\documentclass[11pt]{article}
\usepackage[utf8]{inputenc}
\usepackage[T1]{fontenc}

\usepackage{color}
\usepackage{latexsym}
\usepackage{amssymb}
\usepackage{graphicx}
\usepackage{dsfont} 
\usepackage{amsmath}
\usepackage{eurosym}

\usepackage{algorithm}
\usepackage{algorithmic}
\usepackage[toc,page]{appendix}
\usepackage{hyperref}
\hypersetup{
    colorlinks,
    citecolor=blue,
    filecolor=blue,
    linkcolor=blue,
    urlcolor=blue
}

\newtheorem{Theorem}{Theorem}[part]
\newtheorem{Definition}{Definition}[part]
\newtheorem{Proposition}{Proposition}[part]

\newtheorem{Lemma}{Lemma}[part]
\newtheorem{Corollary}{Corollary}[part]
\newtheorem{Remark}{Remark}[part]

\DeclareMathOperator*{\argmax}{argmax}

\def \Sum{\displaystyle\sum}

\def \Frac{\displaystyle\frac}
\def \Inf{\displaystyle\inf}
\def \Sup{\displaystyle\sup}
\def \Lim{\displaystyle\lim}
\def \Liminf{\displaystyle\liminf}
\def \Limsup{\displaystyle\limsup}

\def \Min{\displaystyle\min}
\def \Argmax{\displaystyle\argmax}

\def \trans{^{\scriptscriptstyle{\intercal }}}

\def \R{\mathbb{R}}

\def \E{\mathbb{E}}
\def \F{\mathbb{F}}

\def \P{\mathbb{P}}

\def \S{\mathbb{S}}

\def \Ac{{\cal A}}

\def \Dc{{\cal D}}

\def \Fc{{\cal F}}

\def \Hc{{\cal H}}
\def \Ic{{\cal I}}
\def \Kc{{\cal K}}
\def \Lc{{\cal L}}

\def \Mc{{\cal M}}
\def \Oc{{\cal O}}
\def \Nc{{\cal N}}

\def \Sc{{\cal S}}
\def \Tc{{\cal T}}
\def \Vc{{\cal V}}

\def \Vc{{\cal V}}

\def \eps{\varepsilon}

\def \ep{\hbox{ }\hfill$\Box$}

\def\Dt#1{\Frac{\partial #1}{\partial t}}

\def\reff#1{{\rm(\ref{#1})}}

\def\beqs{\begin{eqnarray*}}
\def\enqs{\end{eqnarray*}}
\def\beq{\begin{eqnarray}}
\def\enq{\end{eqnarray}}

\begin{document}

\title{Nonzero-sum stochastic impulse games with an application in competitive retail energy markets}

\author{René Aïd\footnote{LEDa, Université Paris Dauphine, PSL research university, rene.aid@dauphine.fr. The author thanks for their support the Finance for Energy Markets, the Chaire EDF-CACIB {\em Finance \& Sustainable Development: Quantitative Approach}  and the ANR project EcoREES ANR-19-CE05-0042.}\and Lamia Ben Ajmia\footnote{ENIT-LAMSIN, University of Tunis El Manar, Tunis, Tunisia, lamia.benajmia@enit.utm.tn}\and M'hamed Gaïgi\footnote{ENIT-LAMSIN, University of Tunis El Manar, Tunis, Tunisia, mhamed.gaigi@enit.utm.tn}\and Mohamed Mnif\footnote{ENIT-LAMSIN, University of Tunis El Manar, Tunis, Tunisia, mohamed.mnif@enit.utm.tn}}


\maketitle

\begin{abstract}
We study a nonzero-sum stochastic differential game with both players adopting impulse controls, on a finite time horizon.
The objective of each player is to maximize her total expected discounted profits. The resolution methodology relies
on the connection between Nash equilibrium and the corresponding system of quasi-variational inequalities (QVIs in short).
We prove, by means of the weak dynamic programming principle for the
stochastic differential game, that the value function of each player is a constrained viscosity solution to the
associated QVIs system in the class of linear growth functions.
We also introduce a family of value functions
converging to our value function of each player, and which is characterized as the unique constrained
viscosity solutions of an approximation of our QVIs system. This
convergence result is useful for numerical purpose. We apply a probabilistic numerical scheme which approximates
the solution of the QVIs system to the case of the competition between two electricity retailers. We show how our model reproduces the qualitative behaviour of electricity retail competition.
\end{abstract}

\vspace{7mm}

\noindent {\bf Key words~:}  stochastic impulse games, Nash equilibrium, viscosity solution.

\vspace{2mm}

\noindent {\bf MSC Classification (2000)~:} 93E20, 49L20,49L25,49N70

\clearpage
\tableofcontents

\section{Introduction}

In this paper, we study a general two-players nonzero-sum stochastic differential game with
impulse controls, on a finite horizon by characterising the associated Nash Equilibrium.
Stochastic differential games have been widely studied.
The theory of two-player zero-sum differential games was pioneered by Isaacs \cite{Isaacs}.
Evans and Souganidis \cite{Evans} studied differential games by means of the viscosity theory, characterizing the
upper and the lower value functions as the unique viscosity solutions to the corresponding Hamilton-
Jacobi-Bellman-Isaacs partial differential equations.
In the literature, many authors focused on stochastic differential game with continuous controls.
The case of impulse games has less been considered. Stochastic impulse games (SIGs) are at the intersection
between differential game theory and stochastic impulse control. They are more realistic in modeling finance problems for example.
Aid et al. \cite{ABCCV19} studied a problem of Nash equilibrium
in a general nonzero-sum impulse game for two players. They derived the corresponding system of quasi-variational inequalities (QVIs) and provides an application in the case of a one-dimensional state problem of forex trading competition between two players. Thanks to a verification theorem, they identified reasonable candidates for the intervention and continuation regions of both players and their strategies.
Cosso \cite{cosso} studied a two-players zero-sum stochastic differential game with both players
adopting impulse controls, on a finite time horizon. The HJBI partial differential equation of the game turns out to
be a double-obstacle quasi-variational inequality. The two obstacles
are implicitly given. He proved that the upper and lower value functions coincide.
He also showed by means of the dynamic programming principle for the
stochastic differential game, that the upper and the lower value functions are the unique viscosity solution to the
HJBI equation, and so the game admits a value.

\noindent Our model is inspired by Aid et al. \cite{ABCCV19} together with the electricity retail pricing model of Basei \cite{Basei19}. In our model,  we assume that retailers buy the energy (electricity or gas) in the wholesale market at the same price. The price of the energy is modeled by the Black-Scholes process.
Finally, retailers sell the energy to their final consumers at a fixed price which could be different for each retailer. Both retailers’ objective is to maximize their total expected discounted profits.
Their instantaneous profits are composed of two parts: sale revenue (market share times retail
price) and sourcing cost (market share times wholesale market price).
The interaction between the opposing retailers is implicitly considered through the market share.
 We formulate the problem as a nonzero-sum stochastic impulse control game.
 The impulse control problem of each player is formulated with three state variables (besides time variable) related to
the wholesale price, the price she offers and the proposed price of the other player which evolve in $\R_+^3$.
The strategy of the retailers consist on changing the energy price in discrete time by using impulse control.
Such model reproduces the fierce competition between the two players which characterizes modern deregulated energy markets.

\noindent To tackle our problem, we use the dynamic programming approach. We prove that the
two players can act simultaneously and the order of intervention of the two players is not important since we can permute
between the two impulse operators of players, which is not the case in the model considered by  Aid et al. \cite{ABCCV19}.
We prove a verification theorem in a general setting
which shows that under some regularity assumptions, if we solve the system of QVIs, then we provide the Nash
equilibrium of the impulse game which characterizes completely the value function of each player.
Then we prove that the value functions of the players
are viscosity solutions of the associated system of QVIs.\\

\noindent In our set-up, it is not obvious to prove the continuity of the value functions of the two players, and
so it is then natural to consider the concept of discontinuous viscosity solutions,
which provides by now a well established  method for dealing with stochastic control problems.  More precisely, we need to consider
constrained viscosity solutions since the state process is allowed to be everywhere in $\R_+^3$
and so the system of QVIs is satisfied even in the boundary of the solvency region.
For the comparison theorem, we consider a small variation of our original model by adding a fixed cost paid by player $i$ when player $j$ makes an intervention.
We prove that the value functions of the perturbed problem converge pointwise to the value functions of the original problem.
We provide a suitable strict supersolution of the system of QVIs solution in the spirit of Crandal et al. \cite{Crand}. In the comparison theorem,
we  construct a test function by introducing a penalization function related to the distance from an element of the solvency region to the boundary,
so that the optimum  associated with the strict supersolution  is not attained on the boundary.
We adapt the arguments of Barles \cite{bar94} which needs a smooth boundary. Similar difficulties have been studied in Akian et al. \cite{akisultak01} and Ly Vath et al.~\cite{MLP07}.
It is known that the function distance is continuously differentiable with bounded derivatives in the smooth part of the boundary (see Gilbarg \cite{giltru77}).
In our context, the boundary is smooth except the corner lines. We prove that one can compare a subsolution with a supersolution
to the system of QVIs provided that one can compare them at the terminal date but also on the corner lines
of the solvency boundary.
It is also known that the comparison theorem is crucial to approximate numerically  the system of QVIs solution.
In this work,
we limit our numerical study to present the numerical scheme and some numerical results for the optimal
intervention strategy. We claim that
the value functions could be obtained as the limit of an iterative
procedure where each step is an optimal stopping problem and the reward function is
related to the impulse operator. We use a numerical approximation algorithm based on
quantization procedure instead of finite difference methods to approximate the value
functions, the intervention and the no-intervention regions of each player. The convergence of
our numerical scheme, in particular, the monotonicity, the stability
and the consistency properties will be postponed in a future research.
In Aïd et al. \cite{aid}, the authors tackled a nonzero-sum stochastic impulse games
by implementing a policy iteration algorithm. To the best of our knowledge,
it is the only paper which is interested to approximate numerically the solution of a system of QVIs, and it opens up the
possibility of finding other numerical schemes to approximate the value functions and the associated optimal strategies.\\

\noindent The paper is organized as follows. Section~\ref{sec:model} is devoted to present the model,
to introduce rigorously the stochastic differential game, Nash equilibrium and
the associated system of QVIs. Then, we study the
properties of the value functions.
In section~\ref{sec:valuefct}
, we give some properties of the value function.
In section~\label{sec:veriftheo}, we prove a verification theorem which
gives a characterization of the Nash equilibrium.
In section~\ref{sec:viscosity}, we prove that the two value functions are viscosity solutions to the system of QVIs followed in section~\ref{sec:compare} by a comparison theorem for an auxiliary family of value functions.
Finally, in section~\ref{sec:numeric},  we explain the numerical scheme and we give some numerical results.

\setcounter{equation}{0} \setcounter{Assumption}{0}
\setcounter{Theorem}{0} \setcounter{Proposition}{0}
\setcounter{Corollary}{0} \setcounter{Lemma}{0}
\setcounter{Definition}{0} \setcounter{Remark}{0}

\section{The model}
\label{sec:model}

This section presents the details of the model.
\subsection{The Dynamics of the state process}
Let $(\Omega,\Fc,\P)$ be a probability space equipped with a filtration $\mathbb{F}:=(\Fc_t)_{0\leq t\leq T}$ supporting a one dimensional
Brownian motion $W$ on a finite horizon $[0,T]$, $T<\infty$. We assume that the wholesale energy price is modeled by a Black-Scholes process.
For a fixed $t\in [0,T]$ and $x\in \R_+$, the process $X^{t,x}$ satisfies
\beq\label{BS}
dX_s^{t,x}= X_{s}^{t,x}(\mu ds+\sigma dW_s), \,s\in[t, T],\,X_s^{t,x}=x,\, s\in[0,t], \label{dynX}
\enq
where $\mu\in \mathbb{R}$ and $\sigma>0$.\\
We consider a nonzero-sum impulse game where two players, which are retailers, can intervene on a continuous
time. After buying the energy, they sell it to consumers. The retailers price is modeled by a stochastic process
$Y := (Y^1, Y^2) \in (0,\infty)\times(0,\infty)$ which is a piecewise-constant process.
More precisely, the control strategy of each player, representing her intervention,
is represented by an impulse control strategy $\alpha^i$  $:=$ $(\tau_{n}^{i},\zeta_{n}^{i})_{n\geq 1}$ for $ i\in\{1,2\}$, where for $t\in[0,T]$,
 $t\leq\tau_{1}^{i}$ $\leq$ $\ldots$ $\tau_{n}^{i}$ $\leq$ $\ldots<T$
are $\F$ stopping times representing the intervention times of the player $i$ and $\zeta_{n}^{i}$, $n$ $\geq$ $1$, are $\Fc_{\tau_{n}^{i}}$-measurable random variables.
The sequence $(\tau_{n}^{i},\zeta_{n}^{i})_{n\geq 1}$ may be a priori finite or infinite.
The process $Y^i$ starts from $Y^i_{t^-}=y^i$ and evolves according to~:
\beq
dY^i_s &=& 0, \;\;\;\;\;  \tau_n^{i} \leq s < \tau_{n+1}^{i} \label{dynY0}  \\
Y^i_{\tau_{n+1}^{i}} &=& \Gamma^i(Y^i_{\tau_{n+1}^{i-}},\zeta_{n+1}^{i}),  \label{dynY}
\enq
where $\Gamma^i(y^i,\zeta^i)=y^i e^{\lambda\zeta^i}$ and $\lambda>0$.
We further assume that, for $ i\in\{1,2\}$, $\zeta^i \in [\zeta_{min},\zeta_{max}]$ where $\zeta_{min}<0<\zeta_{max}$ are fixed constants.
This assumption could be seen as a restriction imposed on the players by the regulator to avoid over-shifting prices. Moreover,
with our specific choice of the function $\Gamma^i$, on the one hand we ensure the positivity of the new prices, on the other
hand, the players could increase the price with a positive impulse (i.e. $\zeta>0$) and decrease it with a negative
impulse (i.e. $\zeta<0$).\\
We then naturally introduce solvency  region: $\Sc =  (0,\infty)^3$, and we denote its closure and its boundary by
\beqs
\bar\Sc  =   \R_+^{3},\,\partial \Sc  =  \bar\Sc\setminus \Sc.
\enqs
In the sequel, we explicit the corner lines of the boundary $\partial \Sc$ as follows:
\beqs
\partial^{x} \Sc=\R^+\times\{(0,0)\},\,\partial^{y_1} \Sc=\{0\}\times (0,\infty)\times \{0\},\,\partial^{y_2}\Sc=\{(0,0)\}\times (0,\infty).
\enqs
We denote by
\beqs
D_0=\partial^{x} \Sc \cup \partial^{y_1} \Sc \cup \partial^{y_2}\Sc.
\enqs
Let $Z=(X,Y^1,Y^2)$ the process living in the space $\Sc$.
We denote by $Z^{t,z,\alpha_1,\alpha_2}$ the process $Z$ that starts at $z=(x,y^1,y^2)$ at time $t$
and is controlled by $\alpha^1$ and $\alpha^2$.\\
For the rest of the paper, to simplify the notations
and when there is no ambiguity, we will use $Z^{t,z}$ instead of $Z^{t,z,\alpha^1,\alpha^2}$. Hence
when we use $Z^{t,z}$, one has to check under which control the process $Z$ evolves. In most cases it can be
easily guessed.
\subsection{Value Function}
\noindent The objective of the player is
to maximize the gain function $J^i$ defined as follows:
\beq
J^i(t,z,\alpha^i,\alpha^j) &=&\E [\int_t^{ T  }e^{-\rho^i (s-t)} f^i(Z^{t,z,\alpha^i,\alpha^j}_s)ds-\sum_{t\leq\tau_{k}^i <T}e^{-\rho^i (\tau_{k}^{i}-t)}\phi^i(Y^{i,t,y^i,\alpha^i}_{(\tau_{k}^{i})^-},\zeta_{k}^{i})\nonumber\\
  &+&e^{-\rho^i (T-t)}g^i(Z^{t,z,\alpha^i,\alpha^j}_T)],
\enq
where, for $ i\in\{1,2\}$, $j\neq i$, $\rho^i$ is a discount factor; $f^i:\R_+^3\longrightarrow \R$ is the running payoff, $g^i:\R_+^3\longrightarrow \R$ is the terminal payoff and $\phi^i:\R_+\times [\zeta_{min},\zeta_{max}] \longrightarrow \R$ is the intervention cost function.
In order to have a realistic model, the running payoff and the terminal payoff of player $i$ depend on the price she asks, the opponent pricing choice and the wholesale price proposed in the market. In our model, the retailer margin of player $i$, represented by $y^i-x$, is truncated by a linear function of the retail price. This regulatory constraint forces retailers to drop down their prices when the wholesale price falls.
In our model $f^i$ and $g^i$ take the form:
\beq\label{payoff}
f^i(x,y^i,y^j)=g^i(x,y^i,y^j)= ((y^i-x) \pi^i(y^i,y^j))\wedge K^ix,
\enq
here the function $\pi^i $ represents the market share and is defined by:
\begin{eqnarray}
 \label{eq:phi}
  \pi^i(y^i,y^j)  =
  \begin{cases}
    1 & \quad \text{if } y^i-y^j\leq -\Delta, \\
    -\frac{1}{2\Delta}(y^i-y^j-\Delta)& \quad \text{if } -\Delta<y^i-y^j<\Delta,\\
    0 & \quad \text{if }  y^i-y^j\geq \Delta,
   \end{cases}
\end{eqnarray}
$K^i$ and $\Delta$ are positive constants.
For $y^i=y^j$, the market is split equally between the two retailers.
The following assumption will be needed in our work:\\
$(\Hc)$ (i) There exist positive constants $C_1^{\phi}$ and $C_2^{\phi}$ s.t.
\begin{eqnarray*}
C_1^{\phi}\leq\phi^i(y^i,\zeta)\leq C_2^{\phi} \,\,\forall (y^i,\zeta) \in \R_+^2\times [\zeta_{min},\zeta_{max}].
\end{eqnarray*}
(ii)The function $\phi^i$ is  continuous.
\begin{Remark}\label{rborne}
  Thanks to the choice of $f^i$ and $g^i$, we deduce that for all $ z \in \Sc$,  $f^i(z)\leq K^ix $ and $g^i(z) \leq K^ix$.
  Moreover  $f^i(z)\geq -(C^i+x)$ and $g^i(z)\geq  -(C^i+x)$, where $C^i$ is a positive constant.
\end{Remark}
\begin{Definition}
We say that the couple $(\alpha^{i,*},\alpha^{j,*})$ is a Nash equilibrium, if, for every other couple of strategies $(\alpha^{i},\alpha^{j})$, we have :
$$ J^i(t,z,\alpha^{i,*},\alpha^{j,*})\geq J^i(t,z,\alpha^i,\alpha^{j,*})\quad \hbox{and} \quad J^j(t,z,\alpha^{i,*},\alpha^{j,*})\geq J^j(t,z,\alpha^{i,*},\alpha^j).$$
\end{Definition}
We define the number of intervention times between $t$ and $T$ by:
\beq\label{numb}
\Nc_t(\alpha^{i}):=Card\{n\;: t\leq\tau_{n}^{i} < T\}.
\enq
We need to impose a uniform integrability condition on $\Nc_t(\alpha^{i})$.
\begin{Definition}
 We define $\Dc_t^i$ as follows
  \beq\label{struniform}
  \Dc_t^i:=\{\alpha^i=(\tau_{n}^{i},\zeta_{n}^{i})_{n\geq 1},\mbox{s.t.}\, \Nc_t(\alpha^{i})  \mbox{ is uniformly integrable}\}.
  \enq
\end{Definition}
For technical reason related to the dynamic programming principle, see Remark 5.2 in Bouchard and Touzi \cite{bou},
we shall restrict the admissible strategies to the following set. For $i\in \{1,2\}$
\beq\label{str}
\Ac^i_t:=\{\alpha^i=(\tau^{i}_{n},\zeta^{i}_{n})_{n\geq 1};\; s.t. \; \tau^{i}_{1}\geq t,\,\alpha^i\in \Dc_t^i\mbox{ and }\alpha^i \mbox{ is independent of } \Fc_t\},
\enq
and therefore, we define the value functions as $(v_i)_{i\in\lbrace 1,2\rbrace}$ on $ [0,T]\times \Sc$ by
\beq\label{ValFct} v^i(t,z) := \Sup_{\alpha^i \in \Ac^i_t} J^i(t,z,\alpha^i,\alpha^{j})\quad i\neq j,
\enq
where $\alpha^{j}$ is the best response against  $\alpha^{i}$.
\begin{Remark}
 It is easy to see that if $\alpha^{i*}=(\tau_{n}^{i*},\zeta_{n}^{i*})$ exists, then, we have
 $\E[\Nc_t(\alpha^{i*})]<\infty$. In our case, to obtain the growth property of the value function, we need a uniform integrability condition of the number of intervention times which is stronger itegrability condition.
 In fact, for $i\in\{1,2\}$, $\alpha^i$ an admissible strategy and $\alpha^j$ the best response against $\alpha^i$, from Assumption $(\Hc)$(i),  we have:
\beqs\label{numberinter}
J^i(t,z,\alpha^{i},\alpha^{j}) &\leq& \E [\int_t^{ T }f^i(Z^{t,z}_s)ds -C_1^{\phi}e^{-\rho^i T} \Nc_t(\alpha^{i})+g^i(Z^{t,z}_T)]\nonumber\\
&\leq& \E [\int_t^{ T }K^i X_s^{t,x}ds -C_1^{\phi}e^{-\rho^i T} \Nc_t(\alpha^{i})+K^iX_T^{t,x}] .\\
\enqs
where the second inequality is obtained from Remark \ref{rborne}. If $\E[\Nc_t(\alpha^{i})]=\infty$, then $J^i(t,z,\alpha^{i},\alpha^{j})=\infty$.
Then, from \reff{numberinter}, we deduce that:
\beqs
J^i(t,z,\alpha^{i},\alpha^{j})=-\infty
\enqs
On the other hand, by choosing the null strategy by player $i$, $J^i(t,z,0,\alpha^{j})>-\infty$,
where here $\alpha^{j}$ denotes the best response against the null strategy.
Hence, an optimal strategy can not have infinite expected number of interventions.
\end{Remark}
\subsection{HJB QVIs}
In this subsection, we provide an heuristic derivation of the system of QVIs.
We define the non-local operator called intervention operator ${\cal M}^i$ for $i\in\{1,2\}$ by:
$$\Mc^i h(t,z)=\sup_{\zeta^i \in [\zeta_{min},\zeta_{max}]}\Big\{h(t,x,\Gamma^{i}(y^i,\zeta^i),y^j)-\phi^i(z,\zeta^i)\Big\},$$
for all locally bounded $h: [0,T]\times S \longrightarrow {\mathbb R}$.
It is known that if the function $h$ is upper semicontinuous, by selection measurable theorem, there exists
a Borel-measurable function $\hat \zeta^{i*}: [0,T]\times S \longrightarrow {\mathbb R}$ such that
$$\Mc^i h(t,z)=h(t,x,\Gamma^{i}(y^i, \hat \zeta^{i*}(t,y^i)),y^j)-\phi^i(z,\hat \zeta^{i*}(t,y^i)),$$
(See for example Proposition 7.33, chapter 7, in \cite{BeSh78}).
Similarly, we define the second intervention operator ${\cal H}^i$ for $i\in\{1,2\}$ by
$$\Hc^i h(t,z)=h(t,x,y^i,\Gamma^j(y^j,\zeta^j)),$$
where $\zeta^j=\Argmax_{\zeta \in [\zeta_{min},\zeta_{max}]}\Big\{h(t,y^i,\Gamma^{j}(y^j,\zeta))-\phi^j(z,\zeta)\Big\}$.\\
Thanks to the following Proposition, the two players could act together. In fact, we can permute between $\Mc^i$ and $\Hc^i$ for $i\in \{1,2\}$
and so we get rid of an assumption in Cosso \cite{cosso} and Aid et al. \cite {ABCCV19}. Such result is obtained as a consequence from the structure
of gain of the two players. In our model, when player $i$ acts an intervention, she pays a cost $\phi^i(y^i,\zeta^i)$. The permutation between the two intervention operators holds.
\begin{Proposition}
For all $(t,z)\in [0,T]\times \Sc$, we have $$\Mc^i\Hc^iv^i(t,z)=\Hc^i\Mc^iv^i(t,z),\forall i,j \in \lbrace 1,2 \rbrace,$$
\end{Proposition}
\underline{\emph{Proof:}} From the definition of $\Mc^i$ and $\Hc^i$ we have
\begin{align*}
\Mc^i\Hc^iv^i(t,z)=&\sup_{\zeta^i \in [\zeta_{min},\zeta_{max}]}\Big\{\Hc^iv^i(t,x,\Gamma^{i}(y^i,\zeta^i),y^j)-\phi^i(y^i,\zeta^i)\Big\},\\
=&\sup_{\zeta^i \in [\zeta_{min},\zeta_{max}]}\Big\{v^i(t,x,\Gamma^{i}(y^i,\zeta^i),\Gamma^{j}(y^j,\zeta^j))-\phi^i(y^i,\zeta^i)\Big\},\\
=&\Mc^iv^i(t,x,y^i,\Gamma^{j}(y^j,\zeta^j))\\
=&\Hc^i\Mc^iv^i(t,x,y^i,y^j)\\
=&\Hc^i\Mc^iv^i(t,z)
\end{align*}
\ep\\
\noindent We define $\overline{\Ic}^i$ as the intervention region of player $j$ and $\Ic^i=\bar \Sc\setminus \overline{\Ic}^i$.
The associated QVIs of our control problem are as follows:
\beq\label{HJBQVI}
\min\{v^i- \Mc^i\Hc^iv^i, v^i-\Hc^iv^i\}&=&0 \quad \hbox{in} \quad\overline{\Ic}^i\label{HIj}\\
\min\{-\Dt{v^i}-\Lc v^i - f^i+\rho^i v^i, v^i-\Mc^iv^i\}&=&0 \quad \hbox{in} \quad \Ic^i \label{HI}\\
v^i&=& g^i\quad \hbox{in } \quad \lbrace T \rbrace\times \bar \Sc \label{ct},
\enq
where
\beq
\Lc v^i(t,z)&=&\mu x\frac{\partial v^i }{\partial x}(t,z) +\frac{\sigma^2 x^2}{2}\frac{\partial^2 v^i }{\partial x^2}(t,z),\label{HJBQVI2}\\
\Mc^iv^i(t,z)&=&\sup_{\zeta^i \in [\zeta_{min},\zeta_{max}]}\Big\{v^i(t,x,\Gamma^{i}(y^i,\zeta^i),y^j)-\phi^i(z,\zeta^i)\Big\},\\
\Hc^iv^i(t,z)&=&v^i(t,x,y^i,\Gamma^j(y^j,\zeta^j)),
\enq
and $ \zeta^j=\Argmax_{\zeta \in [\zeta_{min},\zeta_{max}]}\Big\{v^j(t,y^i,\Gamma^{j}(y^j,\zeta))-\phi^j(z,\zeta)\Big\}$.\\
\begin{Remark}
$\Mc^i v^i(t,z)$ represents the value function
for player $i$ when she makes the optimal impulse in order to increase her value function.
$\Hc^i v^i(t,z)$ represents the value function for player $i$, when player $j$ takes the best intervention.
\end{Remark}
\begin{Remark}
 The system of variational inequalities \reff{HIj}-\reff{HI} is satisfied up to the boundary of $\bar \Sc$.
 Since the dynamics of the process $Z$ is allowed to go everywhere in $\bar \Sc$. We recall that
 the wholesale price $X$ is modeled by a geometric Brownian motion and the process $(Y^i,Y^j)$ are non-negative piecewise constant and at the intervention times,
 the price offered by each player is modeled by the exponential function.
 \end{Remark}
 
 
 \setcounter{equation}{0} \setcounter{Assumption}{0}
\setcounter{Theorem}{0} \setcounter{Proposition}{0}
\setcounter{Corollary}{0} \setcounter{Lemma}{0}
\setcounter{Definition}{0} \setcounter{Remark}{0}

\section{Properties of the value function}
\label{sec:valuefct}

\subsection{Bound on the value functions and Boundary properties}
We start this section with the following lemma which gives a sharp upper bound on the value function $v^i$ for $i\in \{1,2\}$.
\begin{Proposition}
We fix $i\in\{1,2\}$. Under Assumption $(\Hc)$, for all $(t,z)\in [0,T]\times \bar \Sc $, we have
\beq\label{growth}
|v^i(t,z)|\leq C(1+x),
\enq
where $C$ is a positive constant independent of t and z.
\end{Proposition}
\underline{\emph{Proof:}}
We consider the no impulse strategy for player i denoted by $0$ and $\alpha^j$ the best response against the null strategy,
from the definition of $f^i$ and $g^i$, Remark \ref{rborne} and the positivity of the process $Y^i$, we obtain
\beq \label{lower}
v^i(t,z) &\geq& J^i(t,z,0,\alpha^{j})\nonumber\\
&=&\E [\int_t^{ T }e^{-\rho^i(s-t)}f^i(Z^{t,z,0,\alpha^{j}}_s)ds +e^{-\rho^i(T-t)}g^i(Z^{t,z,0,\alpha^{j}}_T)]\nonumber\\
&\geq & -\E [\int_t^{ T }e^{-\rho^i(s-t)}X_s^{t,x}ds +e^{-\rho^i(T-t)}X_T^{t,x}]\nonumber\\
&\geq &-C(1+x),
\enq
where the last inequality is obtained from \reff{BS}, and $C$ is a generic positive constant which could change from line to line.\\
Let $(\alpha^{i},\alpha^{j})$ be an admissible strategy. Using Assumption $(\Hc)(i)$, we have
\beqs
& &J^i(t,z,\alpha^{i},\alpha^{j}) \\
&\leq& \E [\int_t^{ T }e^{-\rho^i(s-t)}f^i(Z^{t,z}_s)ds +e^{-\rho^i(T-t)}g^i(Z^{t,z}_T)]\\
&\leq& \E [\int_t^{ T } K^iX^{t,x}_sds+  K^iX^{t,x}_T]\\
&\leq &\E [\int_t^{ T }  K^i xe^{\mu (s-t)}ds + K^i xe^{\mu (T-t)}].
\enqs
Since $(\alpha^{i},\alpha^{j})$ is arbitrary, we obtain
\beq\label{upper}
v^i(t,z)\leq C x.
\enq
From inequalities \reff{lower} and \reff{upper}, we deduce that (\ref{growth}) is proved.\\
\ep\\
We turn to the behaviour of the values functions on the corner lines of the boundary $\partial \Sc$.
The following proposition determines the value function on the corner lines of the boundary $\partial \Sc$.
\begin{Proposition}\label{cl}
For all $t\in [0,T)$, We have
\beq
v^i(t,z)&=&0\,\,\mbox{ for all } z\in \partial^{y^j}{\cal S}\label{vfc1}\\
v^i(t,z)&=&0\,\,\mbox{ for all } z\in \partial^{y^i}{\cal S}\label{vfc2}\\
v^i(t,z)&=&-\frac{x}{2}(\frac{e^{(\mu-\rho^i)(T-t)}-1}{\mu-\rho^i}+e^{(\mu-\rho^i) (T-t)})\,\,\mbox{ for all } z\in \partial^{x}{\cal S} \label{vfc3}
\enq
\end{Proposition}
\underline{\emph{Proof:}}
$\star$ From the definition of \reff{payoff}, we have $f^j(0,0,y^j)=g^j(0,0,y^j)=0$ for all $y^j>0$.
If player $j$ intervenes, her payoff remains unchanged.
On the other hand, as $\Gamma^i(0,\zeta^i)=0$, any intervention of player $i$ will not change the proposed retailer price, and so it not optimal for the two players to intervene. For the price in wholesale market, it remains equal to zero on $[t,T)$. Then, we deduce \reff{vfc1}.\\
  $\star$ Since $f^i(0,y^i,0)=g^i(0,y^i,0)=0$ for all $y^i>0$, the payoff of player $i$ remains equal to 0.
  As above and since $y^j=0$, any intervention of player $j$ will not change her proposed retailer price. This shows that \reff{vfc2} holds.\\
$\star$ Since $\pi^i(0,0)=\frac{1}{2}$, and for any intervention from one of the two players, their retailer prices remain equal to 0, the market share of each player remains unchanged. From the definition of the value function $v^i$, we deduce that
\beqs
v^i(t,x,0,0)=-\E[\int_t^Te^{-\rho^i (s-t)}\frac{X_s^{t,x}}{2}ds+e^{-\rho^i (T-t)}\frac{X_T^{t,x}}{2}].
\enqs
From the definition of the wholesale price, we have $\E[X_s^{t,x}]=xe^{\mu (s-t)}$. This shows that
\beqs
v^i(t,x,0,0)=-\frac{x}{2}\big(\int_t^Te^{(\mu-\rho^i) (s-t)}ds+e^{(\mu-\rho^i) (T-t)}\big),
\enqs
and so \reff{vfc3} is proved.
\ep\\

\subsection{Terminal condition}

We determine the right terminal condition of the value functions.  We recall that in our model, at the terminal date, each player can not do an intervention.
We obtain such characterization by considering the lower (resp. upper) semicontinuous envelope of $v^i$, for $i\in \{1,2\}$. For all $z\in \Sc$, we set
\beqs
\bar v^{i}(T,z)  :=   \limsup_{\tiny\begin{array}{c}
        (t,z') \rightarrow (T,z) \\
        t < T, z'\in\Sc
       \end{array}} v (t,z'), & &
\underline v^i(T,z)  :=   \liminf_{\tiny\begin{array}{c}
        (t,z') \rightarrow (T,z) \\
        t < T, z'\in\Sc
       \end{array}} v (t,z')
\enqs

\begin{Proposition} \label{lemterm}
We have:
\beqs
{\underline v}^i(T,z) \; = \; {\overline v}^{i}(T,z) \;=\;   g^i(z) , \;\;\; \forall z \in \bar \Sc.
\enqs

\end{Proposition}
\underline{\emph{Proof:}} 1) Fix some $z$ $\in$ $\bar\Sc$ and consider some sequence $(t_m,z_m)_m$ $\in$ $[0,T)\times\Sc$
converging to $(T,z)$ and s.t. $\Lim_{m\rightarrow \infty} v^i(t_m,z_m)$ $=$ $\underline v^i(T,z)$.
By taking the no impulse control strategy for player $i$ noted by 0, and $\alpha^j$ the best response against this strategy by player $j$
on $[t_m,T]$, we have
\beq\label{ct1}
v^i(t_m,z_m) & \geq & \E [\int_{t_m}^{ T  }e^{-\rho^i (s-t_m)}f^i(Z^{t_m,z_m,0,\alpha^j}_s)ds +e^{-\rho^i (T-t_m)}g^i(Z^{t_m,z_m,0,\alpha^j}_T)].
\enq
where
\beqs
Z^{t_m,z_m,0,\alpha^j}=(X^{t_m,x_m},y^i_m,Y^{j,t_m,y^{j}_m,\alpha^j}).
\enqs
We fix $\epsilon >0$. As $(t_m,z_m)_m$ converges to $(T,z)$ when $m$ goes to infinity, there exists $m_0$, such that for
$m\geq m_0$, we have $|z_m-z|\leq \epsilon$.
We fix $m\geq m_0$. From \reff{BS}, we have for $s\leq t_m$, $X^{t_m,x_m}_s=x$ and for $s\geq t_m$, $X^{t_m,x_m}_s=x_me^{(\mu-\frac{\sigma^2}{2})(s-t_m)+\sigma(W_s-W_{t_m})}$, and so $\Lim_{m\rightarrow \infty}X^{t_m,x_m}_s= x$ $dt\otimes d\P$ a.e.
For the jump term $(Y^{j,t_m,y^j_m,\alpha^j}_{s})_s$, we have $Y^{j,t_m,y^j_m,\alpha^j}_s=y_j^m$ for $s\leq t_m$ and for  $s\geq t_m$, we have
\beqs
Y^{j,t_m,y^j_m}_s= y_j^m\prod_{n=0}^{\Nc_{t_m}^j(\hat \alpha)}e^{\lambda \zeta_n^j},
\enqs
where, we recall that $\Nc_{t_m}^j(\alpha^j)$ is the number of interventions of player $j$ between $t_m$ and $T$.
As it is not allowed to make an intervention at the terminal date,
then  $\Lim_{m\rightarrow \infty}\Nc_{t_m}^j( \alpha^j)=0$ a.s. which shows that $\Lim_{m\rightarrow \infty}Y^{j,t_m,y^j_m,\alpha^j}_s=y^j$ $dt\otimes d\P$ a.e. We fix $p>1$, from Remark \ref{rborne}, for $h^i=f^i,g^i$, we have $|h^i(z)|\leq C^i(1+x)$, which implies
\beq\label{majfg1}
 & &\E [\left(\int_{t_m}^{ T  }e^{-\rho^i (s-t_m)}f^i(Z^{t_m,z_m,0,\alpha^j}_s)ds
   +e^{-\rho^i (T-t_m)}g^i(Z^{t_m,z_m,0,\alpha^j}_T)\right)^p]\\
 &\leq& C\left(\E [\left(\int_{t_m}^{ T  }f^i(Z^{t_m,z_m,0,\alpha^j}_s)ds\right)^p]
 +\E [g^i(Z^{t_m,z_m,0,\alpha^j}_T)^p]\right)\nonumber\\
 &\leq& C\E [\int_{t_m}^{ T  }f^i(Z^{t_m,z_m,0,\alpha^j}_s)^pds](T-t_m)^{p-1}
 +\E [g^i(Z^{t_m,z_m,0,\alpha^j}_T)^p])\nonumber\\
 &\leq & C(1+T^{p})(1+\E[\Sup_{s\in [t_m,T]}|X_s^{t_m,x_m}|^{p}]),\nonumber
 \enq
where the second inequality is obtained by using H\"{o}lder inequality and $C$ is a generic constant independent of $m$.
 From the definition of the wholesale price and for $m\geq m_0$, we have
 \beq\label{majx1}
 \E[\Sup_{s\in [t_m,T]}|X_s^{t_m,x_m}|^{p}]\leq C|x_m|^{p}\leq C(1+|x|^{p}).
 \enq
 From \reff{majfg1}-\reff{majx1}, we deduce
 \beqs
 \E [\left(\int_{t_m}^{ T  }e^{-\rho^i (s-t_m)}f^i(Z^{t_m,z_m,0,\alpha^j}_s)ds
   +e^{-\rho^i (T-t_m)}g^i(Z^{t_m,z_m,0,\alpha^j}_T)\right)^p] \leq C(1+|x|^{p}),
 \enqs
 where $C$ is a positive constant independent of $m$. This shows the boundedness of \\
 $\int_{t_m}^{ T  }e^{-\rho^i (s-t_m)}f^i(Z^{t_m,z_m,0,\alpha^j}_s)ds+e^{-\rho^i (T-t_m)}g^i(Z^{t_m,z_m,0,\alpha^j}_T)$ in $L^p(\P)$ for $p>1$,
 which implies the uniform integrability of
 $\left(\int_{t_m}^{ T  }e^{-\rho^i (s-t_m)}f^i(Z^{t_m,z_m,0,\alpha^j}_s)ds+e^{-\rho^i (T-t_m)}g^i(Z^{t_m,z_m,0,\alpha^j}_T)\right)_m$. It yields that
 \beq\label{cteminale1}
{\underline v}^i(T,z)=\Lim_{m\longrightarrow \infty}v^i(t_m,z_m) \geq g^i(z).
\enq
2) Fix some $z$ $\in$ $\bar \Sc$ and consider some sequence $(t_m,z_m)_m$ $\in$ $[0,T)\times\Sc$
converging to $(T,z)$ and s.t. $\lim_{m\rightarrow \infty} v^i(t_m,z_m)$ $=$ $\bar v^{i}(T,z)$. For any $m$, one can find
$\hat\alpha^{m,i}$ $=$ $(\hat\tau_k^{m,i},\hat\zeta_k^{m,i})_k$ $\in$ $\Ac^i_{t_m}$ and $\hat\alpha^{m,j}$ $=$ $(\hat\tau_k^{m,j},\hat\zeta_k^{m,j})_k$ $\in$ $\Ac^j_{t_m}$
the best response against the strategy $\hat\alpha^{m,i}$ s.t.
\beq \label{intervsup}
v(t_m,z_m)\leq
   \E\left[\int_{t_m}^{ T  }f^i(\hat Z^m_s)ds
  -\sum_{t_m\leq\tau_{k}^{i}<T}e^{-\rho^i (\tau_k^i-t_m)}\phi^i(\hat Y^{i,m}_{(\tau_{k}^{m,i})^-},\zeta_{k}^{i,m})+ g^i(\hat Z_T^m) \right] + \frac{1}{m},\nonumber
\enq
where $\hat Z^m$ $=$ $(X^{t_m,x_m},\hat Y^{i,m},\hat Y^{j,m})$ denotes the state process controlled by $\hat\alpha^{m,i}$ and $\hat\alpha^{m,j}$ .
Arguing as above, we have
$\Lim_{m\rightarrow \infty}\Nc_{t_m}^i( \hat\alpha^{m,i})=0$ a.s. which shows that
$\Lim_{m\rightarrow \infty}Y^{i,t_m,y^{i,m},\hat\alpha^{m,i}}_s=y^i$ $dt\otimes d\P$ a.e. For $p>1$ and $m\geq m_0$, we have :
\beqs\label{majfg2}
 \E [\left(\int_{t_m}^{ T  }f^i(\hat Z^m_s)ds
   +g^i(\hat Z^m_T)\right)^p]
&\leq & C(1+(T-t_m)^{p})(1+\E[\Sup_{s\in [t_m,T]}|X_s^{t_m,x_m}|^{p}]),\nonumber\\
&\leq& C(1+|x|^{2p}),\nonumber
 \enqs
 where $C$ is a positive constant independent of $m$. This shows the boundedness of $\int_{t_m}^{ T  }f^i(\hat Z^m_s)ds
   +g^i(\hat Z^m_T)$ in $L^p(\mathbb{P})$ for $p>1$,
 which implies the uniform integrability of $\left(\int_{t_m}^{ T  }f^i(\hat Z^m_s)ds
   +g^i(\hat Z^m_T)\right)_m$. It yields that
 \beq\label{cteminale2}
{\underline v}^i(T,z)=\Lim_{m\longrightarrow \infty}v^i(t_m,z_m) \leq g^i(z).
\enq
The result follows from \reff{cteminale1} and \reff{cteminale2}.
\ep


\setcounter{equation}{0} \setcounter{Assumption}{0}
\setcounter{Theorem}{0} \setcounter{Proposition}{0}
\setcounter{Corollary}{0} \setcounter{Lemma}{0}
\setcounter{Definition}{0} \setcounter{Remark}{0}

\section{Verification theorem}
\label{sec:veriftheo}

We provide in this section a verification theorem for the problem formalized previously. It says that if a
smooth function $v^i$ is solution to the system of QVIs, then it is associated to a Nash equilibrium.
We construct the optimal strategy of each player in an inductive way.

\begin{Proposition} Let $v^i$ be a function from $[0,T]\times\bar \Sc$ to $\mathbb{R}$, with $i \in \lbrace 1,2\rbrace.$ We suppose that Assumption $\Hc$ holds and,
 for $i \in \lbrace 1,2\rbrace$ we have:
\begin{itemize}
\item $v^i  \hbox{ is a solution to } \reff{HIj}-\eqref{ct};$
\item $ v^i \in  C^{1,2}([0,T)\times \bar\Sc)$ satisfying the growth condition \reff{growth}.
\end{itemize}
Let $(t,z) \in [0,T]\times\bar \Sc$ and define the continuation region:
\begin{eqnarray*}
  D^i:=\{(s,\vartheta)\in[t,T]\times {\cal S};\,v^i(s,\vartheta)>\Mc^iv^i(s,\vartheta)\mbox{ or }
  v^i(s,\vartheta)>\Mc^i\Hc^iv^i(s,\vartheta)\}.
\end{eqnarray*}
Define the impulse control $\alpha^{i*}=(\tau_{n}^{i*},\zeta_{n}^{i*})_{n\geq 1}$,
as follows: $\tau_{0}^{i*}=t^-$ and inductively
\begin{eqnarray}
\tau_{n+1}^{i*}&:=&\inf\{ s>\tau_{n}^{i*}\,,(s,Z^{t,z,\alpha^{i*},\alpha^{j*}}_s)\notin D^i\}\label{str1}\\
\zeta_{n+1}^{i*}&:=&\hat \zeta^{i*}(\tau_{n+1}^{i*},Z^{t,z,\alpha^{i*},\alpha^{j*}}_{\tau_{n+1}^{i*-}})\label{str2}
\end{eqnarray}
Assume that $(\alpha^{i,*},\alpha^{j,*}) \in \mathcal{A}_t^i\times\mathcal{A}_t^j $ (i.e. an admissible couple of strategies).\\
Then, $(\alpha^{i,*},\alpha^{j,*})$ is a Nash equilibrium and $v^i(t,z) = J^i(t,z , \alpha^{i,*} , \alpha^{j,*})$ for $i \in \lbrace 1,2 \rbrace$, $j\neq i$.
\end{Proposition}
\noindent\underline{\emph{Proof:}}
We fix $(t,z) \in [0,T]\times\bar \Sc$. We have to prove that
$$ v^i(t,z)=J^i(t,z;\alpha^{i*},\alpha^{j*}),\,\mbox{ for }i \in  \lbrace 1,2\rbrace,\,j\neq i~~~~~~~~ v^i(t,z) \geq J^i(t,z;\alpha^{i},\alpha^{j}).$$
for every $\alpha^i\in \Ac^i_t$ and $\alpha^{j}\in \Ac^j_t$ the best response against $\alpha^{i}$.\\
\textbf{Step 1:} We prove $v^i(t,z) \geq J^i(t,z;\alpha^{i},\alpha^{j}).$\\
We will use the following shortened notation:
$$ Z=Z^{t,z;\alpha^{i},\alpha^{j}}, ~~~~~~ \alpha^i=(\tau^{i}_{n},\zeta^i_{n})_{n\geq 1}. $$
For each $r > 0$ and $n \in \mathbb{N},$ we set
$\tau_r = \Inf\lbrace s > t : X_s \in B(x; r) \rbrace$ is the exit time from the ball with radius $r$ and center x.\\
We denote by $IT^c$ (resp. $IT^d$) the set of intervention times of the two players which coincide (resp. differ).
We apply It\^{o}'s formula to the process
$e^{-\rho^i s}v^i(s, Z_{s} )$
over the interval $[t,\tau_r\wedge T ]$. We take the conditional expectations, the stochastic integral vanishes. Then by taking the expectation, we get:
\beq\label{lam}
&&e^{-\rho^i t}v^i(t,z)\\
&=&
\mathbb{E} \Big[ \int^{\tau_r\wedge T }_{t}e^{-\rho^i u}\left( -\frac{\partial{v^i}}{\partial{t}}(u,Z_u)- \mathcal{L} v^i(u,Z_u) +\rho^i v^i(u,Z_u)\right)du \nonumber\\
  &- &\sum_{t \leq {\tau_{n}^i}<{\tau_r\wedge T }}e^{-\rho^i \tau_{n}^i} (v^i({\tau_{n}^i},Z_{\tau_{n}^{i}})-{v^i}(\tau_{n}^{i},Z_{\tau_{n}^{i-}}))1_{IT^d}(\tau_{n}^i) \nonumber \\
  &-&\sum_{t \leq{\tau_{n}^j}<{\tau_r\wedge T }} e^{-\rho^i \tau_{n}^j}(v^i({\tau_{n}^{j}},Z_{\tau_{n}^{j}})-{v^i}(\tau_{n}^{j},Z_{\tau_{n}^{j-}}))1_{IT^d}(\tau_{n}^j)  \nonumber\\
  &- &\sum_{t \leq{\tau_{n}^i}<{\tau_r\wedge T }} e^{-\rho^i \tau_{n}^i}(v^i({\tau_{n}^{i}},Z_{\tau_{n}^{i}})-{v^i}(\tau_{n}^{i},Z_{\tau_{n}^{i-}}))1_{IT^c}(\tau_{n}^i)
  + e^{-\rho^i( \tau_r\wedge T)}v^i(\tau_r\wedge T ,Z_{\tau_r\wedge T })\Big]. \nonumber
\enq
We now estimate each term in the right-hand side of \eqref{lam}. As for the first term, two cases are possible:
$(s,Z_s)\in {\cal I}^i$ i.e. the QVI (\ref{HI}) holds or  $(s,Z_s)\in \bar{\cal I}^i$ i.e. the second player makes impulses
which could be simultaneously with the first player.
As $\alpha^i\in \Dc^i_t$, we are interested only in the strategies whose number of interventions is finite almost surely, it follows that:
\beq \label{V1Z}
-\frac{\partial{v^i}}{\partial{t}}(s,Z_s)-\mathcal{L} v^i(s,Z_s)+\rho^iv^i(s,Z_s)\geq f^i(Z_s), \mbox{ for all }s\in[t,T]\, \P\,a.s.
\enq
Let us now consider the second term: by the definition of $\mathcal{M}^i v^i$, we have:
\beq \label{QSD}
v^i(\tau_{n}^{i},Z_{\tau_{n}^{i-}})& \geq& \Mc^iv^i(\tau_n^{i},X_{\tau_n^{i-}},Y^i_{\tau_n^{i-}}, Y^j_{\tau_n^{i-}}) \nonumber\\
&=&\sup_{\zeta_n^{i} \in [\zeta_{min},\zeta_{max}]}\Big\{v^i(\tau_n^{i},X_{\tau_n^{i-}},\Gamma^{i}(Y^i_{\tau^n_{i-}},\zeta_n^{i-}), Y^j_{\tau_n^{i-}})-\phi^i(Y^i_{(\tau_n^{i-})},\zeta_n^i)\Big\} \nonumber\\
&\geq& {v^i(\tau_n^{i},X_{\tau_n^{i-}},\Gamma^{i}(Y^i_{\tau_n^{i-}},\zeta_n^{i}),Y^i_{\tau_n^{i-}} )-\phi^i(Y^i_{\tau_n^{i-}},\zeta_{n}^{i})} \nonumber\\
&= &{v^i(\tau_n^{i},Z_{\tau_n^{i}})-\phi^i(Z_{\tau_n^{i-}},\zeta_n^{i})}.
\enq
For the third term, we have $(\Mc^jv^j-v^j)(\tau_n^{j},Z_{\tau_n^{j-}})=0$, hence,
by the definition of $\mathcal{H}^iv^i$, we have:
\beq \label{QSD2}
v^i(\tau_n^{j},Z_{(\tau_n^{j-})})& =& \Hc^iv^i(\tau_n^{j},X_{\tau_n^{j-}},Y^i_{\tau_n^{j-}},Y^j_{\tau_n^{j-}}) \nonumber\\
&=&v^i(\tau_n^{j},X_{\tau_n^{j-}},Y^i_{\tau_n^{j-}},\Gamma^{j}(Y^j_{\tau_n^{j-}},\zeta_n^{j})) \nonumber\\
&= &v^i(\tau_n^{j},Z_{\tau_n^{j}}).
\enq
We focus on the fourth term which is related to the simultaneous interventions case of the two players i.e for time interventions which are in $IT^c$ , we have:
\beq \label{QSD3}
& &v^i(\tau_{n}^{i},Z_{\tau_{n}^{i-}})\nonumber\\
&=&\Hc^iv^i(\tau_{n}^{i},Z_{\tau_{n}^{i-}})\nonumber\\
& =& \Mc^i\Hc^iv^i(\tau_n^{i},X_{\tau_n^{i-}},Y^i_{\tau_n^{i-}}, Y^j_{\tau_n^{i-}}) \nonumber\\
&=&\sup_{\zeta_n^{i} \in [\zeta_{min},\zeta_{max}]}\Big\{v^i(\tau_n^{i},X_{\tau_n^{i-}},\Gamma^{i}(Y^i_{\tau_n^{i-}},\zeta_n^{i-}), \Gamma^{j}(Y^j_{\tau_n^{j-}},\zeta_n^{j}))
-\phi^i(Y^i_{(\tau_n^{i-})},\zeta_n^i)\Big\} \nonumber\\
&\geq& {v^i(\tau_n^{i},X_{\tau_n^{i-}},\Gamma^{i}(Y^i_{\tau_n^{i-}},\zeta_n^{i}),\Gamma^{i}(Y^i_{\tau_n^{i-}},\zeta_n^{j}) )-\phi^i(Y^i_{\tau_n^{i-}},\zeta_{n}^{i})} \nonumber\\
&= &{v^i(\tau_n^{i},Z_{\tau_n^{i}})-\phi^i(Y^i_{\tau_n^{i-}},\zeta_n^{i})}.
\enq
By \eqref{lam}, and the inequalities \eqref{V1Z}-\eqref{QSD2}, it follows that:
\beqs
v^i(t,z)
&\geq& \mathbb{E} \big[ \int^{\tau_r\wedge T }_{t} e^{-\rho^i (s-t)}f^i(Z_u) du
  -\sum_{t \leq{\tau_{n}^{i}}<{\tau_r\wedge T }} e^{-\rho^i (\tau_{n}^i-t)}\phi^i(Y^i_{\tau_n^{i-}},\zeta_n^{i})\\
  &+&e^{-\rho^i (\tau_r\wedge T-t)}{v^i}(\tau_r\wedge T ,Z_{\tau_r\wedge T })   \big].
\enqs
As $f^i(z)\leq K^ix$ (see Remark \ref{rborne}), $\phi^i$ is bounded ( See Assumption $\Hc(i)$),  $v^i$ satisfies linear growth (see inequality \reff{growth}), and since $\alpha^i\in \Ac^i_t$,
we can pass to the limit as $r\rightarrow \infty $ and use the dominated
convergence theorem. it yields that:
\beqs
v^i(t,z)&\geq& J^i(t,z,{\alpha}^i,\alpha^{j}).
\enqs
 \noindent\textbf{Step 2:} $ v^i(t,z) = J^i(t,z,\alpha^{i*},\alpha^{j*} )$, where $(\alpha^{i*},\alpha^{j*})$
 are defined by (\ref{str1})-(\ref{str2})\\
We argue as in Step 1, but here all the inequalities are equalities by the properties of $\alpha^{i*}.$\\
\ep


\setcounter{equation}{0} \setcounter{Assumption}{0}
\setcounter{Theorem}{0} \setcounter{Proposition}{0}
\setcounter{Corollary}{0} \setcounter{Lemma}{0}
\setcounter{Definition}{0} \setcounter{Remark}{0}

\section{Viscosity characterization}
\label{sec:viscosity}
It is well known that the theory of viscosity solutions is a powerful tool to characterize the value function as a solution in a
weaker sense of the associated Hamilton Jacobi Bellman equation.
The couple of value functions are not known to be continuous a priori and so we shall work with
the notion of discontinuous viscosity solutions. For a locally bounded function $u$ on $[0,T]\times\bar \Sc$,
which is the case of the couple of value functions $(v^1,v^2)$, we denote by $\underline{u}$ (resp.
$\overline{u}$) the lower semi-continuous (LSC) (resp. upper semi-continuous (USC)) envelope of $u$. We recall that in general,
$\underline{u}$ $\leq$ $u$ $\leq$ $\overline{u}$, and that $u$ is LSC iff $u$ $=$ $\underline{u}$, $u$ is UCS iff $u$ $=$ $\overline{u}$, and  $u$ is continuous iff
$\underline{u}$ $=$ $\overline{u}$ ($=$ $u$). We denote by $LSC([0,T)\times\bar \Sc)$ (resp. $USC([0,T)\times\bar \Sc)$)  the set of lsc (resp. usc) functions
on $[0,T)\times\bar \Sc$.\\
We work with the suitable notion of constrained viscosity solutions.  The use of constrained viscosity solutions was initiated by \cite{soner86} for first-order equations and applied in stochastic control problems arising in optimal investment problems in \cite{zar88}. The subsolution viscosity is satisfied in $[0, T) \times \bar \Sc$ and the supersolution viscosity is satisfied in $[0, T) \times \Sc$.

\begin{Definition}
(i)  Let $\Oc \subset \bar \Sc$. A locally bounded function $u^i$ on $[0, T] \times \bar \Sc$ is a viscosity
subsolution of \reff{HIj}-\reff{HI} in $[0, T) \times \Oc$ if for all $(\bar{t},\bar{z})\in [0, T) \times \Oc$ and
$\varphi^i\in C^{1,2}([0, T) \times \bar \Sc)$ s.t. $(\overline{u}^i-\varphi^i)(\bar{t},\bar{z}) = 0$
and $(\bar{t},\bar{z})$ is a maximum
of $\overline{u}^i-\varphi^i$  on $[0, T) \times \Oc$, we have
\beq\label{Visc}
\overline{u}^i-\max\{\Mc^i\Hc^i\overline{u}^i,\Hc^i\overline{u}^i\}&\leq&0 \quad \hbox{in} \quad\overline{\Ic}^i\\
\min\{-\Dt{\varphi^i}-\Lc \varphi^i +\rho^i \varphi^i- f^i, \overline{u}^i-\Mc^i\overline{u}^i\}&\leq&0 \quad \hbox{in} \quad \Ic^i
\enq
(ii) Let $\Oc \subset  \Sc$. A locally bounded function $u^i$ on $[0, T] \times \bar \Sc$ is a viscosity
supersolution  of \reff{HIj}-\reff{HI}  in $[0, T) \times \Oc$ if for all $(\bar{t},\bar{z})\in [0, T) \times \Oc$ and
$\varphi^i\in C^{1,2}([0, T) \times \bar \Sc)$ s.t. $(\underline{u}^i-\varphi^i)(\bar{t},\bar{z}) = 0$
and $(\bar{t},\bar{z})$ is a minimum of $\underline{u}^i-\varphi^i$ on $[0, T) \times \Oc$, we have
\beq\label{Visc}
\underline{u}^i-\max\{\Mc^i\Hc^i\underline{u}^i,\Hc^i\underline{u}^i\}&\geq &0 \quad \hbox{in} \quad\overline{\Ic}^i\\
\min\{-\Dt{\varphi^i}-\Lc \varphi^i +\rho^i \varphi^i- f^i, \underline{u}^i-\Mc^i{\underline u}^i\}&\geq&0 \quad \hbox{in} \quad \Ic^i
\enq
(iii) A function $u^i$ on $[0, T] \times \bar \Sc$ is a constrained viscosity solution of \reff{HIj}-\reff{HI} in $[0, T) \times \bar \Sc$
 if it is a viscosity subsolution of \reff{HIj}-\reff{HI} in $[0, T) \times \bar \Sc$ and a viscosity supersolution of \reff{HIj}-\reff{HI} in $[0,T) \times \Sc$.
\end{Definition}

\begin{Remark}
{\rm
There is an equivalent formulation of viscosity solutions, which is useful for proving uniqueness results, see \cite{Crand}~:

\vspace{1mm}

\noindent (i) Let $\Oc$ $\subset$ $\bar\Sc$. A function $u^i$ $\in$  $USC([0,T)\times\bar\Sc)$ is viscosity subsolution
of \reff{HIj}-\reff{HI} in $[0,T)\times\Oc$ if
  \beq
  \min&\{&\overline{u}^i-\Mc^i\Hc^i\overline{u}^i,\overline{u}^i-\Hc^i\overline{u}^i\}\leq 0 \, \hbox{in}\,[0,T)\times\Oc\cap\overline{\Ic}^i,\label{viscosubjet1}\\
  \min&\{&- s_0 - \mu x s_1  - \frac{\sigma^2 x^2 }{2} M_{11}+\rho^iu^i-f^i,  u^i - \Mc^i u^i\} \leq  0\,   \hbox{in}\,[0,T)\times\Oc\cap{\Ic}^i \label{viscosubjet2},
\enq
where  $(s_0,s=(s_k)_{1\leq k\leq 3},M=(M_{k_1k_2})_{1\leq k_1,k_2\leq 3})$ $\in$
$\bar J^{2,+} u^i(t,z)$.

  \vspace{1mm}
\noindent (ii)Let $\Oc$ $\subset$ $\Sc$. A function $u^i$ $\in$  $USC([0,T)\times\bar\Sc)$ is viscosity supersolution
of  \reff{HIj}-\reff{HI} in $[0,T)\times\Oc$ if
  \beq
&  \min&\{\overline{u}^i-\Mc^i\Hc^i\overline{u}^i,\overline{u}^i-\Hc^i\overline{u}^i\}\geq 0 \, \hbox{in} \,[0,T)\times\Oc\cap\overline{\Ic}^i, \label{viscosuperjet1} \\
  &\min&\{- s_0 - \mu x s_1  - \frac{\sigma^2 x^2}{2} M_{11}+\rho^iu^i-f^i,u^i - \Mc^i u^i\}  \geq  0  \,  \hbox{in}\, [0,T)\times\Oc\cap{\Ic}^i,  \label{viscosuperjet2}
\enq
where $(s_0,s=(s_k)_{1\leq k\leq 3},M=(M_{k_1k_2})_{1\leq k_1,k_2\leq 3})$ $\in$
$\bar J^{2,-} u^i(t,z)$.
  \vspace{1mm}

\noindent (iii) A locally bounded function $u^i$  on $[0,T)\times\bar\Sc$ is a constrained viscosity solution to
\reff{HIj}-\reff{HI}  if  $\bar u^i$ satisfies \reff{viscosubjet1}-\reff{viscosubjet2} in $[0,T)\times\bar\Sc$, where $(s_0,s,M)$ $\in$ $\bar J^{2,+} \bar u^i(t,z)$,
and $\underline u^i$ satisfies \reff{viscosuperjet1}-\reff{viscosuperjet2} in $[0,T)\times\Sc$, where $(s_0,s,M)$ $\in$ $\bar J^{2,-} \underline u^i(t,z)$.

\vspace{1mm} Here $J^{2,+} u(t,z)$ is the parabolic second order superjet defined by~:
\beqs
& & J^{2,+} u^i(t,z) \; = \;  \{ (s_0,s,M) \in \R\times\R^3\times\S^3~:  \\
 & &  \;\;\;  \limsup_{\tiny\begin{array}{c}
      (t',z')\rightarrow (t,z) \\
       (t',z') \in [0,T)\times\Sc\end{array}}
\frac{u^i(t',z') - u^i(t,z) - s_0 (t'-t) - s.(z'-z) - \frac{1}{2}(z'-z).M(z'-z)}{|t'-t|+|z'-z|^2} \; \leq \;  0  \},
\enqs
$\S^3$ is  the set of symmetric $3\times 3$ matrices, $\bar J^{2,+} \bar u^i(t,z)$ is its closure~:
\beqs
\bar J^{2,+} u^i(t,x) &=& \left\{ (s_0,s,M) \; = \; \lim_{m\rightarrow\infty} (s_0^m,s^m,M^m) \;\;\;
\mbox { with } (s_0^m,s^m,M^m) \in J^{2,+} u^i(t_m,z_m) \right. \\
& & \;\;\; \left. \mbox { and } \; \lim_{m\rightarrow\infty} (t_m,z_m,u^i(t_m,z_m)) \; = \; (t,z,u^i(t,z)) \right\},
\enqs
and  $J^{2,-} u^i(t,z)$ $=$ $-J^{2,+}(- u^i)(t,z)$, $\bar J^{2,-} u(t,z)$ $=$ $-\bar J^{2,+}(- u^i)(t,z)$.

}
\end{Remark}
The dynamic programming principle (DPP) is the key tool to characterize the couple of value functions $(v^i,v^j)$ defined by \reff{ValFct} as a viscosity solution of
the system \reff{HIj}-\reff{HI}.
We will use the weak version of the dynamic programming principle (WDPP) developed by Bouchard and Touzi \cite{bou} (See Theorem 3.5 and Corollary 3.6).
In fact, the DPP is intuitive and its proof requires that the value function is measurable on the first stage and the delicate measurable selection theorem holds.
The WDPP is a powerful tool to avoid these technical points. Under the assumptions of independence, causality,
stability under concatenation and consistency with deterministic initial data,
which could be easily checked in our case, we have the following result.
\begin{Proposition}
Let Assumption $(\Hc)$ holds. Then,\\
(i)  for all  $(\alpha^i,\alpha^j)\in \Ac^i_0\times \Ac^j_0$, we have  $J^i(.,.,\alpha^i,\alpha^j)\in LSC([0,T)\times\bar \Sc)$. \\
(ii) for any stopping time $\theta\in\Tc_{t,T}$ (set of stopping times valued in $[t,T]$) we have the following inequalities:
\beq\label{DPP} v^i(t,z) &\leq &\Sup_{\alpha^i \in \Ac^i_t} \E [\int_t^{\theta }e^{-\rho^i (s-t)}f^i(Z^{t,z,\alpha^i,\alpha^{j}}_s)ds\nonumber\\
& & -\sum_{t\leq \tau_{k}^{i}<\theta}e^{-\rho^i (\tau_k^i-t)}\phi^i(Y^{i,t,y^i,\alpha^i}_{(\tau_{k}^{i})^-},\zeta_{k}^{i}) +e^{-\rho^i (\theta-t)}{\bar v}^i(\theta,Z^{t,z,\alpha^i,\alpha^{j}}_{\theta})],
\enq
\beq\label{DPP2}
v^i(t,z) &\geq& \Sup_{\alpha^i \in \Ac^i_t} \E [\int_t^{\theta }e^{-\rho^i (s-t)}f^i(Z^{t,z,\alpha^i,\alpha^{j}}_s)ds\nonumber\\
  & & -\sum_{t\leq \tau_{k}^{i}<\theta}e^{-\rho^i (\tau_k^i-t)}\phi_i(Y^{i,t,y^i,\alpha^i}_{(\tau_{k}^{i})^-},\zeta_{k}^{i})+e^{-\rho^i (\theta-t)}{\underline v}^i(\theta,Z^{t,z,\alpha^i,\alpha^{j}}_{\theta})].
    \enq
\end{Proposition}
\underline{\emph{Proof:}}
(i) Let $(\alpha^i,\alpha^j) \in \Ac^i_0\times \Ac^j_0$.
From Remark \ref{rborne} and \reff{growth}, the random variable $ \int_{t_n}^{T }e^{-\rho^i (s-t_n)}f^i(Z^{t_n,z_n,\alpha^i,\alpha^{j}}_s)ds+e^{-\rho^i (T-t_n)}{g}^i(Z^{t_n,z_n,\alpha^i,\alpha^{j}}_{T})$ is bounded from below by an integrable random variable. By applying Fatou's lemma, we obtain
\beqs
\Liminf_{\stackrel{(t_n,z_n)\rightarrow (t,z)}{t_n\in [0,T),\,z_n\in(0,+\infty)^3}}&\E& [\int_{t_n}^{T }e^{-\rho^i (s-t_n)}f^i(Z^{t_n,z_n,\alpha^i,\alpha^{j}}_s)ds+e^{-\rho^i (T-t_n)}{g}^i(Z^{t_n,z_n,\alpha^i,\alpha^{j}}_{T})]\\
 \geq &\E& [\int_{t}^{T }e^{-\rho^i (s-t)}f^i(Z^{t,z,\alpha^i,\alpha^{j}}_s)ds+e^{-\rho^i (T-t)}{g}^i(Z^{t,z,\alpha^i,\alpha^{j}}_{T})].
\enqs
Under Assumption  $(\Hc)(i)$, and since $\alpha^i\in \Dc_0^i$, the sequence
$\left(- \Sum_{t_n\leq \tau_{k}^{i}< T}e^{-\rho^i ( \tau_k^i-t_n)}\phi_i(Y^{i,t_n,y_n^i,\alpha^i}_{(\tau_{k}^{i})^-},\zeta_{k}^{i})\right)_n$ is uniformly integrable, and so we obtain
\beqs
 \Lim_{n\rightarrow \infty}\E[- \sum_{t_n\leq \tau_{k}^{i}< T}e^{-\rho^i (\tau_k^i-t_n)}\phi^i(Y^{i,t_n,y_n^i,\alpha^i}_{(\tau_{k}^{i})^-},\zeta_{k}^{i})]
=\E[- \sum_{t\leq \tau_{k}^{i}< T}e^{-\rho^i (\tau_k^i-t)}\phi_i(Y^{i,t,y^i,\alpha^i}_{(\tau_{k}^{i})^-},\zeta_{k}^{i})].
\enqs
It yields that
$$\underline J^i(t,z,\alpha^i,\alpha^j)= \Liminf_{\stackrel{(t_n,z_n)\rightarrow (t,z)}{t_n\in [0,T),\,z_n\in(0,+\infty)^3}}J^i(t_n,z_n,\alpha^i,\alpha^j)\geq J^i(t,z,\alpha^i,\alpha^j).$$
  This shows that $J^i(.,.,\alpha^i,\alpha^j)\in LSC([0,T)\times\bar \Sc)$.\\
 (ii)    Since  the strategy $\alpha^i$ is independent of $\Fc_t$, and following \cite{bou} (See their Remark 5.2),
one can prove
\beq v^i(t,z) := \Sup_{\alpha^i \in \Ac^i_0} J^i(t,z,\alpha^i,\alpha^{j})\quad i\neq j,\nonumber
\enq
which means that one could get rid of the dependence of time $t$ in the set $\Ac^i_t$. \\
The proof of inequality \reff{DPP} is then similar to  \cite{bou}.\\
We fix $R>0$ and a control $(\alpha^i,\alpha^j)\in \Ac^i_0\times \Ac^j_0$. We consider the stopping time
$\tau^R:=\inf\{s\geq t\mbox{ s.t. }|Z^{t,z,\alpha^i,\alpha^j}_s-z|>R\}$.
Then $Z^{t,z,\alpha^i,\alpha^j}1_{[t,\theta\wedge \tau^R)}$ is
  $\L^{\infty}$ bounded. By Corollary 3.6 in \cite{bou}, we have
  \beqs
  v^i(t,z) &\geq & \E [\int_t^{\theta \wedge \tau^R}
    e^{-\rho^i (s-t)}f^i(Z^{t,z,\alpha^i,\alpha^{j}}_s)ds
  - \Sum_{t\leq \tau_{k}^{i}<\theta\wedge \tau^R}e^{-\rho^i ( \tau^i_k-t)}\phi^i(Y^{i,t,y^i,\alpha^i}_{(\tau_{k}^{i})^-},\zeta_{k}^{i})\\
  & & +e^{-\rho^i (\theta\wedge \tau^R-t)}{\underline v}^i(\theta\wedge \tau^R,Z^{t,z,\alpha^i,\alpha^{j}}_{\theta\wedge \tau^R})].
\enqs
Using Assumption $(\Hc)(i)$, Remark \ref{rborne}, inequality \reff{growth}, and  since $(\alpha^i,\alpha^j)\in \Dc_0^i\times \Dc_0^j$,
Fatou's lemma implies
\beqs
& \Liminf_{R\longrightarrow \infty}& \E [\int_t^{\theta \wedge \tau^R}
   e^{-\rho^i (s-t)} f^i(Z^{t,z,\alpha^i,\alpha^{j}}_s)ds
  - \Sum_{t\leq \tau_{k}^{i}<\theta\wedge \tau^R}e^{-\rho^i ( \tau^i_k-t)}\phi^i(Y^{i,t,y^i,\alpha^i}_{(\tau_{k}^{i})^-},\zeta_{k}^{i})\\
  & &+e^{-\rho^i (\theta\wedge \tau^R-t)}{\underline v}^i(\theta\wedge \tau^R,Z^{t,z,\alpha^i,\alpha^{j}}_{\theta\wedge \tau^R})]\\
&\geq&\E [\int_t^{\theta}
    e^{-\rho^i (s-t)}f^i(Z^{t,z,\alpha^i,\alpha^{j}}_s)ds
  - \Sum_{t\leq \tau_{k}^{i}<\theta}e^{-\rho^i (\tau_k^i-t)}\phi^i(Y^{i,t,y^i,\alpha^i}_{(\tau_{k}^{i})^-},\zeta_{k}^{i})\\
& &  +e^{-\rho^i (\theta-t)}{\underline v}^i(\theta,Z^{t,z,\alpha^i,\alpha^{j}}_{\theta})].
\enqs
As $(\alpha^i,\alpha^j)$  is chosen arbitrary, we obtain \reff{DPP2}.
\ep\\ \\
In the following lemma, we give some auxiliary results on the intervention operators $\Mc^i$ and $\Hc^i$ which will be useful later.
 \begin{Lemma}\label{Lemmatechnical}
    Let $u^i$ be a locally bounded function on $[0,T)\times\bar\Sc$.

\noindent (i) If $u^i$ is lsc, then $\Mc^i u^i$ is lsc.

\noindent (ii) If $u^i$ is usc, then $\Mc^i u^i$ is usc.

\noindent (iii) If $u^i$ is lsc (resp. usc), then $\Hc^i u^i$ is lsc (resp. usc).
 \end{Lemma}
 \underline{\emph{Proof:}}
 We prove only the first and the second statement.\\
\noindent (i)  Fix some $(t,z)$ $\in$ $[0,T)\times\bar\Sc$ and let $(t_n,z_n)$ be  a sequence in $[0,T)\times\bar\Sc$ converging to $(t,z)$ and s.t. $\Mc^i u^i(t_n,z_n)$ converges to
${\underline {\Mc^i u^i}}(t,z)$. Then, using also the lower semi-continuity of $u^i$, we have:
\beqs
\Mc^i u^i(t,z) & = & \sup_{\zeta\in [\zeta_{min},\zeta_{max}]} u^i(t,x,\Gamma^i(y^i,\zeta),y^j))  \;  \leq \;  \sup_{\zeta \in [\zeta_{min},\zeta_{max}]}
\liminf_{n\rightarrow\infty} u^i(t_n,x_n,\Gamma^i(y_n^i,\zeta),y_n^j) \\
&\leq & \liminf_{n\rightarrow\infty}  \sup_{\zeta\in [\zeta_{min},\zeta_{max}]}  u^i(t_n,x_n,\Gamma^i(y_n^i,\zeta),y_n^j)  \; \leq \;
\lim_{n\rightarrow\infty}  \Mc^i u^i(t_n,z_n) \; = \; \underline{\Mc^i u^i}(t,z),
\enqs
which shows the lower semi-continuity of $\Mc^i u^i$.

\noindent (ii) Fix some $(t,z)$ $\in$ $[0,T)\times \Sc$ and let $(t_n,z_n)_{n\geq 1}$ be a sequence in $[0,T)\times \Sc$ converging to
$(t,z)$ when $n$ goes to infinity. Since $u^i$ is usc,
for each $n$ $\geq$ $1$, by selection measurable theorem, there exists a sequence $(\hat\zeta_n)_{n\geq 1}$ with $\hat\zeta_n$ $\in$ $[\zeta_{min},\zeta_{max}]$ such that:
\beqs
\Mc^i  u^i(t_n,z_n) &=& u^i(t_n,x_n,\Gamma^i(y^i_n,\hat\zeta_n),y_n^j) , \;\;\; \forall n \geq 1.
\enqs
As $[\zeta_{min},\zeta_{max}] $ is a compact, the sequence $(\hat\zeta_n)_{n\geq 1}$  converges, up to a subsequence, to some
$\hat\zeta$ $\in$ $[\zeta_{min},\zeta_{max}]$.   Therefore, we get~:
\beqs
\Mc^i u^i(t,z) \; \geq \; u^i(t,x,\Gamma^i(y^i,\hat\zeta),y^j)  & \geq & \limsup_{n\rightarrow\infty}  u^i(t_n,x_n,\Gamma^i(y^i_n,\hat\zeta_n),y^j_n)  \;
= \;  \limsup_{n\rightarrow\infty} \Mc^i u^i(t_n,z_n),
\enqs
which proves that $\Mc^i u^i$ is usc.\\
\noindent (iii) The proof is similar to (i) and (ii).\\
 \ep\\
The following technical lemma is needed to prove that the value functions $(v^1,v^2)$ are viscosity solutions of the associated  QVIs system.
\begin{Lemma}\label{LV2}
For all $(t,z) \in {\Ic}^i$ and $i\in\{1,2\}$ we have:
\beq
\underline v^i(t,z)-\Mc^i\underline v^i(t,z)&\geq&0.\label{ineq1}
\enq
For all $(t,z) \in \overline{\Ic}^i$ and $i\in\{1,2\}$ we have:
\beq
\min\{\underline v^i(t,z)-\Mc^i\Hc^i\underline{v}^i(t,z),\underline v^i(t,z)-\Hc^i\underline{v}^i(t,z)\}&\geq& 0,\label{ineq3}\\
\min\{\overline v^i(t,z)-\Mc^i\Hc^i\overline{v}^i(t,z),\overline v^i(t,z) -\Hc^i\overline{v}^i(t,z)\}&\leq& 0.\label{ineq4}
\enq
\end{Lemma}
\underline{\emph{Proof:}}
We fix $i\in\{1,2\}.$ Let $\alpha^i$ an admissible control associated to the player $i$ where she chooses
to intervene immediately at time $t$ with an arbitrary size $\zeta^i \in [\zeta_{min},\zeta_{max}]$.
By applying the WDPP \reff{DPP2}, for all $(t,z) \in [0,T)\times \Sc$, we have that:
\beqs
v^i(t,z) &\geq& \underline v^i(t,x,\Gamma^i(y^i,\zeta^i),y^j)- \phi^i(z,\zeta^i).
\enqs
As $\phi^i$ is continuous (See Assumption $(\Hc)(ii)$), the right hand side of the last inequality is lsc, and so we deduce that:
\beqs
\underline v^i(t,z) &\geq& \underline v^i(t,x,\Gamma^i(y^i,\zeta^i),y^j)- \phi^i(z,\zeta^i).
\enqs
From the arbitrariness of $\zeta^i$, we deduce
\beqs
\underline v^i(t,z) &\geq& \Mc^i \underline v^i(t,z),
\enqs
and so the inequality (\ref{ineq1}) is proved.\\
We fix $(t,z) \in \overline{\Ic}^i$. Two cases are possible: the player $j$ decides to intervene or the two players $i$ and $j$ intervene simultaneously.
In the first case, we denote by $\zeta^j$ the size of intervention of player $j$.
Then, by applying the WDPP \reff{DPP2}, we have
\beqs
v^i(t,z)\geq \underline v^i(t,x,y^i,\Gamma^j(y^j,\zeta^j))=\Hc^i \underline v^i(t,z).
\enqs
As the right hand side of the last inequality is lsc, we deduce that
\beq\label{under1}
\underline v^i(t,z)\geq \Hc^i \underline v^i(t,z).
\enq
In the second case, if the players $i$ and $j$ decide to intervene in $(t,z)$, then, by applying the WDPP \reff{DPP2},
for all $(t,z) \in [0,T)\times \Sc$, we have that
\beqs
v^i(t,z) &\geq& \underline v^i(t,x,\Gamma^i(y^i,\zeta^i),\Gamma^j(y^j,\zeta^j))- \phi^i(z,\zeta^i).
\enqs
Arguing as above, we deduce
\beq\label{under2}
\underline v^i(t,z) &\geq& \Mc^i \Hc^i \underline v^i(t,z).
\enq
From \reff{under1} and \reff{under2}, we deduce that \reff{ineq3}. It remains to prove \reff{ineq4}.\\
From the WDPP \reff{DPP}, if player $j$ decides to intervene, we obtain
\beqs\label{over1}
{ v}^i(t,z)\leq \Hc^i \overline{v}^i(t,z).
\enqs
As the right hand side of the obove inequality is usc, then
\beq\label{over1}
{\overline v}^i(t,z) \leq \Hc^i \overline{v}^i(t,z).
\enq
If the two players intervene simultanously, then  by using \reff{DPP} and
arguing as above, we obtain
\beq\label{over2}
\overline v^i(t,z) \leq {\Mc}^i{ \Hc}^i\overline v^i(t,z).
\enq
From \reff{over1} and \reff{over2}, we deduce \reff{ineq4}.
\ep\\
\noindent Our main result of the section is the following.
\begin{Theorem}\label{Visc}
  The function $v^i, \; i\in\{1,2\}$ defined by \reff{ValFct} is a constrained viscosity solution
  of \reff{HIj}-\reff{HI} in $[0, T) \times \Sc$.
\end{Theorem}
\underline{\emph{{Proof of supersolution property on $[0,T)\times \Sc$:}}}\\
We fix $i\in\{1,2\}$. Let $(\bar{t},\bar{z}) \in [0,T)\times \Sc$ and $\varphi^i\in C^{1,2}([0, T) \times \bar \Sc)$ s.t. $(\underline{v}^i-\varphi^i)(\bar{t},\bar{z}) = 0$
    and $(\bar{t},\bar{z})$ is a minimum of $\underline{v}^i-\varphi^i$ on $[0, T) \times \Sc$.
      Two cases are possible: $(\bar{t},\bar{z}) \in { \bar \Ic}^i$ or  $(\bar{t},\bar{z}) \in { \Ic}^i$. If $(\bar{t},\bar{z}) \in {\bar \Ic}^i$, then
      from Lemma \ref{LV2}, we have $\Min\{\underline v^i(\bar t,\bar z)-\Mc^i\Hc^i \underline v^i(\bar t,\bar z) ,\underline v^i(\bar t,\bar z)-\Hc^i \underline v^i(\bar t,\bar z)\}\geq0$ (See inequality (\ref{ineq3})).
      If $(\bar{t},\bar{z}) \in { \Ic}^i$, then we have $\underline v^i(\bar t,\bar z)-\Mc^i\underline v^i(\bar t,\bar z)\geq0$  (See inequality (\ref{ineq1})). It remains to show that
\beq\label{inegsursolinter}
-\Dt{\varphi^i}(\bar{t},\bar{z})-\Lc \varphi^i(\bar{t},\bar{z}) +\rho^i\varphi^i(\bar{t},\bar{z})- f^i(\bar{z}) \geq0,
\enq
From the definition of $\underline{v}^i$, there exists a sequence $(t_m,z_m)_{m\geq 1} \in \Sc$
s.t. $(t_m,z_m)$ and $v^i(t_m,z_m)$ converge respectively to $(\bar{t},\bar{z})$ and
 $\underline{v}^i(\bar{t},\bar{z})$ as  $m$ goes to infinity. By continuity of $\varphi^i$, we also have that
$\gamma_m$ $:=$ $v^i(t_m,z_m)-\varphi^i(t_m,z_m)$ converges to $0$ as $m$ goes to infinity.
We denote by $\tau_{1,m}^{j*}$ the first optimal interversion of player $j$ when the state process $Z$ starts from $z_m$ at time $t_m$.
Two cases are possible $\tau_{1,m}^{j*}=t_m$ or $\tau_{1,m}^{j*}>t_m$. If the first case occurs, then
the player $j$ makes an intervention immediately. By applying the WDPP \reff{DPP2},
$v^i(t_m,z_m)\geq \Hc^i\underline v^i(t_m,z_m)$.
As $t_m$ converges to $\bar t$ as $m$ goes to infinity, and since $\Hc^i\underline v^i$ is lsc (See Lemma \ref{Lemmatechnical} (iii)), we obtain $\underline v^i(\bar t,\bar z)\geq \Hc^i\underline v^i(\bar t,\bar z)$. This means that $(\bar{t},\bar{z}) \in {\bar \Ic}^i$, which is false, and so we must have $\tau_{1,m}^{j*}>t_m$.
As $\bar t<T$ and $\bar z\in \Sc$,  for $m$ large enough, one could have $t_m$ $<$ $T$ and there exists $\delta$ $>$ $0$ such
$B(z_m,\delta/2):=\{z \mbox{ s.t } |z-z_m|\leq \frac{\delta}{2}\}\subset \Sc$.
Let us then consider the admissible control in $\Ac^i_{t_m}$ with no impulse until the first exit time $\tau_m$, from the ball $B(z_m,\delta/2)$
before $T\wedge\tau_{1,m}^{j*}$
of the associated state process $Z_s$, defined by:
\beqs
\tau_m &:=& \inf\left\{ s \geq t_m~: |Z_s^{t_m,z_m} - z_m| \geq \delta/2 \right\} \wedge T  \wedge\tau_{1,m}^{j*}.
\enqs
Consider also a strictly positive sequence $(h_m)_m$
s.t. $h_m$ and $\gamma_m/h_m$ converge to zero as $m$ goes to infinity. By using the WDPP \reff{DPP2}
for $v_i(t_m,z_m)$ and $\hat\tau_m$ $:=$ $\tau_m\wedge (t_m+h_m)$, we get~:
\beqs
v^i(t_m,z_m)& = & \gamma_m+\varphi^i(t_m,z_m) \\
&\geq&  \E[\int_{t_m}^{\hat\tau_m }e^{-\rho^i(s-t_m)}f^i(Z^{t_m,z_m}_s)ds
  +e^{-\rho^i(\hat\tau_m-t_m)}\underline v^i(\hat\tau_m,Z_{\hat\tau_m}^{t_m,z_m})]\\
&\geq&  \E[\int_{t_m}^{\hat\tau_m }e^{-\rho^i(s-t_m)}f^i(Z^{t_m,z_m}_s)ds
  +e^{-\rho^i(\hat\tau_m-t_m)}\varphi^i(\hat\tau_m,Z_{\hat\tau_m}^{t_m,z_m})].
\enqs
Applying It\^o's formula to $\varphi^i(s,Z_s^{t_m,z_m})$ between
$t_m$ and $\hat\tau_m$ and noting that the integrant of the stochastic integral term is bounded, we obtain by taking expectation~:
\beq \label{inegsursolinter2}
\frac{\gamma_m}{h_m}
+ \E\left[\frac{1}{h_m} \int_{t_m}^{\hat\tau_m}e^{-\rho^i(s-t_m)}\left( - \Dt{\varphi^i} - \Lc \varphi^i+\rho^i\varphi^i-f^i\right)(s,Z_s^{t_m,z_m} ) ds \right]
& \geq & 0.
\enq
By continuity a.s. of $Z_s^{t_m,z_m}$, we have for $m$ large enough, $\hat\tau_m$ $=$ $t_m+h_m$, and so by the mean-value theorem,
the random variable inside the expectation in \reff{inegsursolinter2} converges a.s. to
$(- \Dt{\varphi^i} - \Lc \varphi^i+\rho^i\varphi^i-f^i)(\bar t,\bar z)$ as $m$ goes to infinity.  Since this random variable is also bounded
by a constant independent of $m$, we conclude by the dominated convergence theorem and we obtain \reff{inegsursolinter}.\\
\noindent\underline{\emph{{Proof of subsolution property on $[0,T)\times \bar \Sc$:}}}\\
We fix $i\in\{1,2\}$ and $(\bar t,\bar z) \in[0,T)\times\bar \Sc$. If $(\bar{t},\bar{z}) \in {\bar \Ic}^i$, then
from Lemma \ref{LV2}, we have $\Min\{\overline v^i(\bar t,\bar z)-\Mc^i\Hc^i \overline v^i(\bar t,\bar z),\overline v^i(\bar t,\bar z)-\Hc^i \overline v^i(\bar t,\bar z)\} \leq 0$
(See inequality (\ref{ineq4})).
 If $(\bar{t},\bar{z}) \in {\Ic}^i$, we consider $\varphi^i$ $\in$ $C^{1,2}([0,T)\times\bar \Sc)$
s.t. $\overline{v}^i(\bar t,\bar z)$ $=$ $\varphi^i(\bar t,\bar z)$ and $\varphi^i$ $\geq$ $\overline{v}^i$ on $\bar\Sc$.
If $\overline{v}^i(\bar t,\bar z)$ $\leq$ $\Mc^i \overline{v}^i(\bar t,\bar z)$ then the subsolution inequality  holds
trivially.
Consider now the case where $\overline{v}^i(\bar t,\bar z)$ $>$ $\Mc^i \overline{v}^i(\bar t,\bar z)$ and argue by contradiction
by assuming on the contrary that
\beqs
\bar \eta \; := \; -\Dt{\varphi^i}(\bar t,\bar z) - \Lc \varphi^i(\bar t,\bar z)+\rho^i\varphi^i(\bar t,\bar z)-f^i(\bar z) & > & 0.
\enqs
By continuity of $f^i$ ( see equation\reff{payoff}), $\varphi^i$ and its derivatives, and since  $(\bar t,\bar z) \in [0,T)\times \bar \Sc$, there exists
some $\delta_0$ $>$ $0$ s.t. $\bar t+\delta_0$ $<$ $T$ and for all $0<\delta\leq\delta_0$ and for all $ (t,z) \in ((\bar t-\delta)_+,\bar t+\delta)\times (B(\bar z,\delta)\cap \bar \Sc)$, we have:
\beq  \label{inegphidelta}
-\Dt{\varphi^i}(t,z) - \Lc  \varphi^i(t,z)+\rho^i\varphi^i(t, z)-f^i(z) > \frac{\bar \eta}{2}.
\enq
From the definition of $\overline{v}_i$, there exists a sequence $(t_m,z_m)_{m\geq 1}$ $\in$
$((\bar t-\delta/2)_+,\bar t+\delta/2)\times B(\bar z,\delta/2)\cap \bar \Sc$ s.t.
$(t_m,z_m)$ and $v^i(t_m,z_m)$ converge respectively to $(\bar{t},\bar{z})$ and
 $\overline{v}_i(\bar{t},\bar{z})$ as  $m$ goes to infinity. By continuity of $\varphi^i$, we also have that
$\gamma_m$ $:=$ $v^i(t_m,z_m)-\varphi^i(t_m,z_m)$ converges to $0$ as $m$ goes to infinity.
By the WDPP \reff{DPP}, given $m$ $\geq$ $1$,  for any stopping time $\tau$ valued in
$[t_m,T]$, we have
\beq \label{inegdpp1}
v^i(t_m,z_m) & \leq & \Sup_{\alpha^i \in \Ac^i_{t_m}} \E [\int_{t_m}^{\tau}e^{-\rho^i(s-t_m)}f^i(Z^{t_m,z,\alpha^i,\alpha^{j}}_s)ds
  \nonumber\\
  &-& \sum_{t_m\leq \tau_{k}^{i}<\tau}e^{-\rho^i(\tau_{k}^{i}-t_m)}\phi^i(Y^{i,t_m,y^i_m,\alpha^i}_{(\tau_{k}^{i})^-},\zeta_{k}^{i})+e^{-\rho^i(\tau-t_m)}{\bar v}^i(\tau,Z^{t_m,z_m,\alpha^i,\alpha^{j}}_{\tau})]\nonumber\\
&\leq& \E[\int_{t_m}^{\tau }e^{-\rho^i(s-t_m)}f^i(\hat Z_s^m)ds-
  \sum_{t_m\leq \tau_{m,k}^{i}<\tau}e^{-\rho^i(\tau_{m,k}^{i}-t_m)}\phi^i(\hat Y^{i,m}_{(\tau_{m,k}^{i})^-},\zeta_{m,k}^{i}) \nonumber\\
  &+&e^{-\rho^i(\tau-t_m)}{\bar v}^i(\tau,\hat Z_\tau^m)] + \frac{1}{m}.
\enq
Here $\hat Z^m$ is the state process, starting from $z_m$ at $t_m$, and controlled by
$(\alpha_m^{i},\alpha_m^{j})$ an $\frac{1}{m}-\mbox{optimal strategy}$.

We define $\bar\tau_m$ $:$ $=$ $\tau_{m,1}^{i}\wedge\tau_\delta^m\wedge\tau_{m,1}^j$ where,
\beqs
\tau_\delta^m &=& \inf\left\{ s \geq t_m~: \hat Z_s^m \notin {\mathring B}(z_m,\delta/2)\cap \bar \Sc \right\} \wedge (t_m+\delta/2)
\enqs
is the first exit time before $t_m+\delta/2$  of  $\hat Z^m$ from the intersection between open ball $\mathring B(z_m,\delta/2)$ and $\bar \Sc$.
We claim that $\tau_{m,1}^{i}>t_m$. In fact, if $\tau_{m,1}^{i}=t_m$, the player $i$ makes an intervention immediately, then
\beqs
v^i(t_m,z_m)  \leq {\bar v}^i(t_m,x_m,\Gamma^i(y^i_m,\zeta_{m,k}^{i}),y^j_m) -\phi^i(y^i_m,\zeta_{m,k}^{i})+ \frac{1}{m}.
\enqs
Sending m to infinity and taking the supremum over all admissible interventions,
we obtain ${\bar v}^i(\bar t,\bar z)\leq \Mc^i{\bar v}^i(\bar t,\bar z)$ which is false in this particular case.\\
We claim that  $\tau_{m,1}^{j}>t_m$. In fact if $\tau_{m,1}^{j}=t_m$, the player $j$ makes an intervention immediately, then
\beqs
v^j(t_m,z_m)  \leq  {\bar v}^j(t_m,x_m,y^i_m,\Gamma^j(y^j_m,\zeta_{m,k}^{j})) -\phi^j(y^j_m,\zeta_{m,k}^{j})+\frac{1}{m}.
\enqs
Sending m to infinity, we deduce that ${\bar v}^j(\bar t,\bar z)\leq \Mc^j{\bar v}^j(\bar t,\bar z)$ and so $(\bar t,\bar z)\in {\bar \Ic}^i$ which is false. As the process $\hat Z^m$ is c\`adl\`ag, we have $\tau_\delta^m >t_m$ $P$ a.s. Then we deduce that $\bar\tau^m>t_m$ $P$ a.s.
As the filtration $(\Fc_t)_{0\leq t\leq T}$ is a Brownian filtration, the stopping time $\bar\tau^m$ is predictable and
so there exists a sequence of stopping times $(\tau_{n,m})_n$ s.t. $t_m<\bar\tau_{n,m} <\bar\tau^m$, $\bar \tau_{n,m} \uparrow \bar\tau^m$,  when $n$ goes to infinity.
We fix $n$ and $m$, we choose $\tau=\bar\tau_{n,m}$ in \reff{inegdpp1}. Then, we get:
\beq
v^i(t_m,z_m) & \leq &
\E[\int_{t_m}^{\bar\tau_{n,m} }e^{-\rho^i(s-t_m)}f^i(\hat Z_s^m)ds\\
& &+ e^{-\rho^i(\bar\tau_{n,m}-t_m)}v^i(\bar\tau_{n,m}, \hat Z^m_{\bar\tau_{n,m}})] +\frac{1}{m}\nonumber.  \label{intersoussol}
\enq
Now, since $\overline{v}_i$ $\leq$ $\varphi_i$, we obtain:
\beqs
\varphi^i(t_m,z_m) + \gamma_m  & \leq & \E[\int_{t_m}^{\bar\tau_{n,m}}e^{-\rho^i(s-t_m)}f^i(\hat Z_s^m)ds\\
 & & +e^{-\rho^i(\bar\tau_{n,m}-t_m)}\varphi^i(\bar\tau_{n,m}, \hat Z_{\bar\tau_{n,m}}) ]+\frac{1}{m}.
\enqs
By applying It\^o's formula to $\varphi^i(s,\hat Z_s^m)$ between $t_m$ and $\bar\tau_{n,m}$, and using  \reff{inegphidelta}, we then get:
\beqs
\gamma_m & \leq & \E\left[\int_{t_m}^{\bar \tau_{n,m}} e^{-\rho^i(s-t_m)}\left(\Dt{\varphi^i} + \Lc\varphi^i-\rho^i\varphi^i+f^i\right)(s,\hat Z_s^m) ds\right]\\
 &  \leq & - \frac{\bar \eta}{2\rho^i} (1-\E[e^{-\rho^i(\bar\tau_{n,m} - t_m)}]).
\enqs
This implies that: $\lim_{n,m\rightarrow\infty} \E[e^{-\rho^i\bar\tau_{n,m}}] =  e^{-\rho^i \bar t} $, and
\beqs \label{limbartaum}
\lim_{m\rightarrow\infty} \E[e^{-\rho^i\bar\tau_{m}}]  =  e^{-\rho^i \bar t},
\enqs
where the last equality is deduced by using monotone convergence theorem. Then, along a subsequence, one could have $\Lim_{m\rightarrow\infty} \bar \tau_{m}=\bar t$ a.s.
Applying the WDPP \reff{DPP} for $\tau=\bar \tau_m$, we obtain:
\beq
& &v^i(t_m,z_m) \nonumber\\
& \leq&
\E[\int_{t_m}^{\bar\tau_m }e^{-\rho^i(s-t_m)}f^i(\hat Z_s^m)ds]+\E[e^{-\rho^i(\bar\tau_m-t_m)} \bar v^i(\bar\tau_m, \hat Z^m_{\bar\tau_m})1_{\tau_\delta^m<\tau_{m,1}^i\wedge\tau_{m,1}^j}]  \nonumber\\
&+&\E[e^{-\rho^i(\bar\tau_m-t_m)} \bar v^i (\bar\tau_m, \hat X^{m}_{\bar\tau_m},\hat Y^{i,m}_{\bar\tau_m},\Gamma^j(\hat Y^{j,m}_{\bar\tau_m^{-}},\zeta_{m,1}^{j} ) )  1_{\tau_{m,1}^j\leq \tau_{m,1}^i\wedge \tau_\delta^m,\,\tau_{m,1}^j\neq \tau_{m,1}^i} ]\nonumber \\
&+&\E[e^{-\rho^i(\bar\tau_m-t_m)} (\bar v^i(\bar\tau_m, \hat X^{m}_{\bar\tau_m},\Gamma^i(\hat Y^{i,m}_{\bar\tau_m^{-}},\zeta_{m,1}^i),\Gamma^j(\hat Y^{j,m}_{\bar\tau_m^{-}},\zeta_{m,1}^{j}))\nonumber\\
 & -&\phi^i(\hat Z^m_{(\tau_{m,1}^i)^-},\zeta_{m,1}^i)) 1_{\tau_{m,1}^i=\tau_{m,1}^j\leq \tau_\delta^m} ]\nonumber \\
& + &  \E[e^{-\rho^i(\bar\tau_m-t_m)}( \bar v^i(\bar\tau_m,\hat X^{m}_{\bar\tau_m}, \Gamma^i(\hat Z^m_{\bar\tau^{m,-}},\zeta_{m,1}^i),\hat Y^{j,m}_{\bar\tau_m} -\phi^i(\hat Z^m_{(\tau_{m,1}^i)^-},\zeta_{m,1}^i)\nonumber)\\
& &  1_{ \tau_{m,1}^i \leq \tau_\delta^m \wedge \tau_{m,1}^j, \tau_{m,1}^i \neq \tau_{m,1}^j } ] +\frac{1}{m}
\nonumber  \\
&\leq & \E[\int_{t_m}^{\bar\tau_m }e^{-\rho^i(\bar\tau_m-t_m)}f^i(\hat Z_s^m)ds] +\E[e^{-\rho^i(\bar\tau_m-t_m)} \bar v^i(\bar\tau_m, \hat Z^m_{\bar\tau_m})
  1_{\tau_\delta^m<\tau_{m,1}^i\wedge\tau_{m,1}^j}] \nonumber\\
&+&\E[e^{-\rho^i(\bar\tau_m-t_m)}\Hc^i  \bar v^i(\bar\tau_m, \hat Z^m_{\bar\tau_m^{-}})1_{\tau_{m,1}^j\leq \tau_{m,1}^i\wedge \tau_\delta^m,\,\tau_{m,1}^j\neq \tau_{m,1}^i}]\nonumber\\
&+&
\E[e^{-\rho^i(\bar\tau_m-t_m)}
  \Mc^i\Hc^i  \bar v^i(\bar\tau_m, \hat Z^m_{\bar\tau_m^{-}}) 1_{\tau_{m,1}^i=\tau_{m,1}^j\leq \tau_\delta^m}]\nonumber\\
&+& \E[e^{-\rho^i(\bar\tau_m-t_m)}\Mc^i  \bar v^i(\bar\tau_m, \hat Z^m_{\bar\tau_m^{-}})1_{\tau_{m,1}^i\leq \tau_\delta^m\wedge\tau_{m,1}^j,\tau_{m,1}^i\neq\tau_{m,1}^j}]+\frac{1}{m}
 ,  \label{intersoussol}
\enq
which implies,
\beq\label{contrasuper}
v^i(t_m,z_m) & \leq & \sup_{\tiny\begin{array}{c}
      |t'-\bar t| < \delta \\ |z'-\bar z| < \delta
\end{array}} f^i(z')\frac{(1-\E[e^{-\rho^i(\bar\tau_{m} - t_m)}])}{\rho^i} \\
&+&
       \sup_{\tiny\begin{array}{c}
      |t'-\bar t| < \delta \\ |z'-\bar z| < \delta
       \end{array}} \bar v^i(t',z') \E[e^{-\rho^i(\bar\tau_m-t_m)}1_{\tau_\delta^m<\tau_{m,1}^i\wedge\tau_{m,1}^j}] \nonumber \\
       &+&
       \sup_{\tiny\begin{array}{c}
      |t'-\bar t| < \delta \\ |z'-\bar z| < \delta
       \end{array}} \Hc^i \bar v^i(t',z') \E[e^{-\rho^i(\bar\tau_m-t_m)}1_{\tau_{m,1}^j\leq\tau_{m,1}^i\wedge \tau_\delta^m,\,\tau_{m,1}^j\neq\tau_{m,1}^i}]\nonumber\\
        &+&
       \sup_{\tiny\begin{array}{c}
      |t'-\bar t| < \delta \\ |z'-\bar z| < \delta
       \end{array}} \Mc^i\Hc^i \bar v^i(t',z') \E[e^{-\rho^i(\bar\tau_m-t_m)}1_{\tau_{m,1}^j=\tau_{m,1}^i\leq\tau_\delta^m}]\nonumber\\
        &+ &\sup_{\tiny\begin{array}{c}
      |t'-\bar t| < \delta \\ |z'-\bar z| < \delta
       \end{array}} \Mc^i \bar v^i(t',z')    \E[e^{-\rho^i(\bar\tau_m-t_m)}1_{\tau_{m,1}^i\leq \tau_\delta^m\wedge\tau_{m,1}^j,\,\tau_{m,1}^i\neq\tau_{m,1}^j}]+\frac{1}{m} \nonumber.
\enq
Since
\beq\label{limtau}
\Lim_{m\rightarrow\infty}\bar \tau_m =\Lim_{m\rightarrow\infty} \tau_{m,1}^j\wedge\tau_{m,1}^i\wedge \tau_\delta^m=\bar t\,\,\mbox{ a.s.},
\enq
by the dominated convergence theorem, we obtain:
\beq\label{arg1}
\Lim_{m\rightarrow\infty}\E[e^{-\rho^i(\bar\tau_{m} - t_m)}]=1
\enq
For the second term in the right hand side of \reff{contrasuper}, we define $A_m:=\{\tau_\delta^m<\tau_{m,1}^i\wedge\tau_{m,1}^j\}$ and $A:=\cap_{m\geq 0}\cup_{n\geq m} A_n$.
We assume that $\P(A)>0$ and we consider $\omega\in A$.
As $\Lim_{m\rightarrow \infty}\tau_\delta^m(\omega)=\bar t$ and  $\Lim_{m\rightarrow \infty}x_m= \bar x$, then
\beqs
\Lim_{m\rightarrow \infty}\hat X_{\tau_\delta^m(\omega)}^m(\omega)=\Lim_{m\rightarrow \infty}x_me^{(\mu-\frac{\sigma^2}{2})(\tau_\delta^m-t_m)+\sigma(W_{\tau_\delta^m}(\omega)-W_{t_m}(\omega))}=\bar x.
\enqs
As the processes $\hat Y^{i,m}$ and $\hat Y^{j,m}$ are constant on $[t_m,\tau_\delta^m(\omega)]$ and $\Lim_{m\rightarrow \infty}(y^i_m,y^j_m)= (\bar y^i,\bar y^j)$, then
$\Lim_{m\rightarrow \infty}\hat Y_{\tau_\delta^m(\omega)}^{i,m}(\omega)=\Lim_{m\rightarrow \infty}y^i_m=\bar y^i$,
and $\Lim_{m\rightarrow \infty}\hat Y_{\tau_\delta^m(\omega)}^{j,m}(\omega)=\Lim_{m\rightarrow \infty}y^j_m=\bar y^j$. This shows that
$\Lim_{m\rightarrow \infty}\hat Z_{\tau_\delta^m(\omega)}^m(\omega)=\bar z$. On the other hand, from the definition of $\tau_\delta^m$, we have $\hat Z_{\tau_\delta^m(\omega)}^m(\omega)\notin B(z_m,\frac{\delta}{2})$ i.e $|\hat Z_{\tau_\delta^m(\omega)}^m(\omega)-z_m|>\frac{\delta}{2}$.
Sending $m$ to infinity, we obtain $0>\frac{\delta}{2}$ which is false, and so $\P(A)=0$. On the other hand, $0\leq \P(A_m)\leq \P(\cup_{n\geq m} A_n)$, and
$\P(\cup_{n\geq m} A_n)\downarrow P(A)$ whem $m$ goes to infinity, then $\P(A_m)$ goes to zero when $m$ goes to infinity. By the dominated convergence theorem, we obtain:
\beq\label{arg2}
\Lim_{m\rightarrow\infty}\E[e^{-\rho^i(\bar\tau_m-t_m)}1_{\tau_\delta^m<\tau_{m,1}^i\wedge\tau_{m,1}^j}]=0.
\enq
For the third term in the right hand side of \reff{contrasuper},
we consider the set $B_m:=\{\tau_{m,1}^j\leq\tau_{m,1}^i\wedge\tau_\delta^m,\,\tau_{m,1}^j\neq\tau_{m,1}^i\}$.
We define $B:=\cap_{m\geq 0}\cup_{n\geq m} B_n$. We assume that $\P(B)>0$ and we consider $\omega\in B$.
The player $j$ is the first to make an intervention. As the process $\hat Z_.^m$ is right continuous, then $\tau_{m,1}^j(\omega)>\bar t$.
Since $\Lim_{m\rightarrow \infty}\tau_{m,1}^j(\omega)=\bar t$, and $\hat Z_.^m$ is right continuous
then $\Lim_{m\rightarrow \infty}\hat Z_{\tau_{m,1}^{j-}(\omega)}^m(\omega)=\bar z$. It yields that $(\bar t,\bar z)\in \bar \Ic^i$ which is false.
Then, we have $\P(B)=0$. Arguing as above, by the dominated convergence theorem, we obtain:
\beq\label{arg3}
\Lim_{m\rightarrow\infty}\E[e^{-\rho^i(\bar\tau_m-t_m)}1_{\tau_{m,1}^j<\tau_{m,1}^i\wedge \tau_\delta^m}]=0.
\enq
For the fourth  term in the right hand side of \reff{contrasuper},
we consider the set $C_m:=\{\tau_{m,1}^j=\tau_{m,1}^i\leq\tau_\delta^m\}$.
We define $C:=\cap_{m\geq 0}\cup_{n\geq m} C_n$. Arguing as above, one could prove that $\P(C)=0$ and by the dominated convergence theorem, we obtain:
\beq\label{arg4}
\Lim_{m\rightarrow\infty}\E[e^{-\rho^i(\bar\tau_m-t_m)}1_{\tau_{m,1}^j=\tau_{m,1}^i\leq\tau_\delta^m}]=0.
\enq
Sending $m$ to infinity, from \reff{contrasuper},\reff{arg1}-\reff{arg4} we deduce that
\beqs
\overline{v}^i(\bar t,\bar z) & \leq &
\sup_{\tiny\begin{array}{c}
      |t'-\bar t| < \delta \\ |z'-\bar z| < \delta
       \end{array}} \Mc^i \bar v^i(t',z')
\enqs
Sending $\delta$ to $0$, by Lemma \ref{Lemmatechnical}, we have
\beqs
\overline{v}_i(\bar t,\bar z) &\leq &  \Mc_i \overline{v}_i(\bar t,\bar z).
\enqs
As $\overline{v}^i(\bar t,\bar z)$ $>$ $\Mc^i \overline{v}^i(\bar t,\bar z)$, we obtain a contradiction.
This shows that the subsolution property is satisfied.
\ep


\setcounter{equation}{0} \setcounter{Assumption}{0}
\setcounter{Theorem}{0} \setcounter{Proposition}{0}
\setcounter{Corollary}{0} \setcounter{Lemma}{0}
\setcounter{Definition}{0} \setcounter{Remark}{0}

\section{Comparison principle of viscosity solutions}
\label{sec:compare}

We turn to uniqueness result by proving a comparison principle for
discontinuous viscosity solutions to the QVIs \reff{HIj}-\reff{HI}.
The comparison theorem is based on the Ishii technique \cite{ishiimpul} to produce a strict supersolution and we adapt arguments to handle the nonlocal operators.
The idea is to build a test function so that the maximum associated with the strict supersolution is not attained on the boundary.
Furthermmore, to construct a strict supersolution, we consider an auxiliary problem by  adding a fixed cost $\kappa$ $>$ $0$ to the intervention operator $\Hc^i$.
This means that for an intervention of size $\zeta^j$ made by player $j$, the operator $\Hc^i_\kappa$ is defined by:
\beqs
\Hc_\kappa^i h(t,z):=h(t,x,y^i,\Gamma^j(y^j,\zeta^j))-\kappa,
\enqs
for all locally bounded function $h: [0,T]\times S \longrightarrow {\mathbb R}$.
We define the value function of the perturbed problem $v^i_\kappa$ by:
\beqs
v^i_\kappa(t,z) := \Sup_{\alpha^i \in \Ac^i_t}J_\kappa^i(t,z,\alpha^i,\alpha^j),
\enqs
where
\beqs
J_\kappa^i(t,z,\alpha^i,\alpha^j)&:=&\E [\int_t^{ T  }e^{-\rho^i (s-t)} f^i(Z^{t,z,\alpha^i,\alpha^j}_s)ds
  -\sum_{t\leq \tau_{k}^i< T}e^{-\rho^i (\tau_{k}^{i}-t)}\phi^i(Z^{t,z,\alpha^i,\alpha^j}_{(\tau_{k}^{i})^-},\zeta_{k}^{i})\\
  &-&\kappa \sum_{t\leq \tau_{k}^j<T}e^{-\rho^i (\tau_{k}^{j}-t)}
+e^{-\rho^i (T-t)}g^i(Z^{t,z,\alpha^i,\alpha^j}_T)],\,(t,z)\in [0,T]\times \bar \Sc
\enqs
One can prove that the value function $v^i_\kappa$ satisfies the growth condition \reff{growth}, and is a constrained viscosity solution in $[0,T)\times \bar \Sc$ of the following system of QVIs:
\beq
\min\{v^i_\kappa-\Mc^i\Hc_\kappa^i v^i_\kappa,v^i_\kappa-\Hc_\kappa^i v^i_\kappa\}&=&0 \quad \hbox{in} \quad\overline{\Ic}^i\label{HIjp}\\
\min\{-\Dt{v^i_\kappa}-\Lc v^i_\kappa +\rho^i v^i_\kappa- f^i, v^i_\kappa-\Mc^iv^i_\kappa\}&=&0 \quad \hbox{in} \quad \Ic^i,\label{HIp}
\enq
with terminal condition:
\beq
v^i_\kappa(T,z)&=& g^i(z)\quad \hbox{in } \quad \lbrace T \rbrace\times \bar \Sc \label{ctlam},
\enq
and boundary conditions:
\beq\label{clkappa}
v_\kappa^i(t,z)=
\begin{cases}
     0& \quad \text{if } (t,z)\in [0,T)\times\partial^{y^1}{\cal S}\cup \partial^{y^2}{\cal S}, \\
     -\frac{x}{2}(\frac{e^{(\mu-\rho^i)(T-t)}-1}{\mu-\rho^i}+e^{(\mu-\rho^i) (T-t)}) & \quad \text{if }  (t,z)\in [0,T)\times\partial^{x}{\cal S}.
  \end{cases}
\enq
\begin{Remark}
The proof of viscosity property for $v^i_\kappa$, $i\in \{1,2\}$,  follows the same lines of arguments as for $v^i$, $i\in \{1,2\}$.
For the growth condition \reff{growth},
similar arguments holds as for the lower bound (see inequality \reff{lower}) and for the upper bound (see inequality \reff{upper}), since the number of interventions of each player is uniformly integrable and so satisfies $\Sup_{\alpha^i\in \Dc_0^i}\E[\Nc_T(\alpha^{i})]<\infty$.
\end{Remark}
The first result of this section is a convergence result.
\begin{Proposition}
For $i\in \{1,2\}$, the sequence $(v^i_\kappa)_\kappa$ is nonincreasing, and converges  pointwise on
$[0,T]\times \bar\Sc $ towards  $v^i$  as $\kappa$ goes to zero.
\end{Proposition}
\noindent\underline{\emph{Proof:}}
\noindent   Notice that for any $0<\kappa_1\leq\kappa_2$, $(t,z)\in [0,T)\times \bar\Sc$ and $(\alpha^i,\alpha^j)\in \Ac_t^i\times \Ac_t^j$, we have
\beqs
J_{\kappa_1}^i(t,z,\alpha^i,\alpha^j)\geq J_{\kappa_2}^i(t,z,\alpha^i,\alpha^j).
\enqs
As $\alpha^i$ is arbitrary and $\alpha^j$ is the best response againt $\alpha^i$, we deduce that $v_{\kappa_1}^i(t,z)\geq v_{\kappa_2}^i(t,z)$.
 This shows that  the sequence $(v^i_\kappa)$ is nonincreasing, and is upper-bounded by the  value function $v^i$, so that
\beq \label{liminfvi}
\lim_{\kappa \downarrow 0} v_\kappa^i (t,z) &\leq& v^i(t,z), \;\;\; \forall (t,z) \in [0,T]\times\bar\Sc.
\enq
For the converse inequality, we fix some point $(t,z)$ $\in$ $[0,T]\times \bar\Sc$. From the representation \reff{ValFct} of $v^i(t,z)$,
there exists for any  $n$ $\geq$ $1$, an $1/n$-optimal control denoted by $\alpha_n^{i}$ et $\alpha_n^{j}$ the best response against $\alpha_n^{i}$ such that:
\beq\label{1noptimal}
J^i(t,z,\alpha_n^i,\alpha_n^j)\geq v^i(t,z)-\frac{1}{n}.
\enq
On the other hand, we have $J^i(t,z,\alpha_n^i,\alpha_n^j)=J_\kappa^i(t,z,\alpha_n^i,\alpha_n^j)+\kappa \E [\sum_{t\leq \tau_{k}^j<T}e^{-\rho^i (\tau_{k}^{j}-t)}]$. From \reff{1noptimal},
we deduce that:
\beqs
J_\kappa^i(t,z,\alpha_n^i,\alpha_n^j)+\kappa \E [\sum_{t\leq \tau_{k}^j<T}e^{-\rho^i (\tau_{k}^{j}-t)}]\geq v^i(t,z)-\frac{1}{n}.
\enqs
Taking the supremum over all admissible strategies, we obtain:
\beqs
v_\kappa^i(t,z)+\kappa \Sup_{\alpha^{j}\in \Ac^j_t}\E [\sum_{t\leq \tau_{k}^j<T}e^{-\rho^i (\tau_{k}^{j}-t)}]\geq v^i(t,z)-\frac{1}{n}
\enqs
By Sending $\kappa$ to zero and $n$ to infinity, and since $\alpha^j\in \Dc_t^j$, we obtain
\beq\label{limsupvi}
\lim_{\kappa \downarrow 0} v_\kappa^i (t,z) &\geq& v^i(t,z).
\enq
The result follows from \reff{1noptimal} and \reff{limsupvi}.
\ep\\
Thanks to the following proposition, we produce a suitable strict viscosity supersolution, which will allow us to compare a viscosity supersolution and a viscosity subsolution on $[0,T)\times \Sc$.
\begin{Proposition}
Let
$w_\kappa^i,\; i \in\{1,2\}$ be a lower semi-continuous viscosity supersolution on $[0,T)\times\Sc$ of the QVIs \reff{HIjp}-\reff{HIp}. For all $(t,z)\in [0,T)\times \Sc$, we consider
\beqs
h^i(t,z)=\left(A^i+B^ix^2+C^i \log(y^i+1)+D^i \log(y^j+1)\right)e^{\rho(T-t)},
\enqs
where
\beq\label{param}
&\rho>(2\mu+\sigma^2-\rho^i)_+,\, A^i>\frac{K^i}{\rho^i+\rho},\,B^i>\frac{K^i}{\rho^i+\rho-2\mu-\sigma^2},\nonumber\\
&\frac{C^{\phi}_1}{\lambda\zeta_{max}}e^{-\rho T}>C^i>0,\mbox{ and }\frac{\kappa}{\lambda\zeta_{max}}e^{-\rho T}>D^i>0.
\enq
Then, for each $\gamma \in (0,1)$, $w_\kappa^{i,\gamma}=(1 - \gamma)w_\kappa^i + \gamma g^i $ is a  viscosity supersolution of
\begin{eqnarray} \label{strictsuper}
\left\{
\begin{array}{rlll}
\min\{w_\kappa^{i,\gamma}-\Mc^i\Hc^i_\kappa w_\kappa^{i,\gamma},w_\kappa^{i,\gamma}-\Hc^i_\kappa w_\kappa^{i,\gamma}\}&\geq&\gamma\eta^{'} \quad \hbox{in} \quad\overline{\Ic}^i,\\
\min \{-\Dt{w_\kappa^{i,\gamma}}-\Lc  w_\kappa^{i,\gamma} +\rho^i w_\kappa^{i,\gamma}- f^i, w_\kappa^{i,\gamma}-\Mc^i w_\kappa^{i,\gamma}\} &\geq& \gamma\eta  \quad \hbox{in}
\quad \mathcal{I}^i,\\
\end{array}
\right.
\end{eqnarray}
where $\eta$ and $\eta^{'}$ are positive constants.
\end{Proposition}
\noindent\underline{\emph{Proof:}}

Let  $(t,z) \in \Ic^i$ and $\gamma \in (0,1).$
We fix $\zeta^i \in  [\zeta_{min},\zeta_{max}]$, then
\beq
& &h^i(t,z)-\left(  h^i(t,x,\Gamma^i(z^i,\zeta^i),y^j)-\phi^i(y^i,\zeta^i) \right) \\
&=&h^i(t,z)-\left(  A^i+B^ix^2+C^i \log(y^i e^{\lambda \zeta^i}+1)+D^i \log(y^j+1)\right)e^{\rho(T-t)}  +\phi^i(y^i,\zeta^i) \nonumber \\
&\geq& h^i(t,z)-\left(  A^i+B^ix^2+C^i \log( (y^i e^{\lambda\zeta_{max}}+e^{\lambda\zeta_{max}})+D^i \log(y^j+1)\right)e^{\rho(T-t)}  +C^{\phi}_1\nonumber\\
&\geq& h^i(t,z)-\left(  A^i+B^ix^2+C^i \log(  e^{\lambda\zeta_{max}}(y^i+1))+D^i \log(y^j+1)\right)e^{\rho(T-t)} +C^{\phi}_1\nonumber\\
&\geq& h^i(t,z)-\left(  A^i+B^ix^2+C^i \log (y^i+1)+C^i\log (e^{\lambda\zeta_{max}})+D^i \log(y^j+1)\right)e^{\rho(T-t)} +C^{\phi}_1\nonumber\\
&\geq& h^i(t,z)-h^i(t,z)- C^i \lambda \zeta_{max}e^{\rho(T-t)}+C^{\phi}_1\nonumber\\
&\geq& - C^i \lambda \zeta_{max}e^{\rho(T-t)}+C^{\phi}_1,\nonumber
\enq
where the first inequality is obtained since $\phi^i $ satisfies $\Hc_1$ and $\zeta_{max}$ is positive.
From the choice of $C^i$ (See \reff{param}) and from the arbitrariness of $\zeta^i$, we deduce that:
\begin{align}
  h^i(t,z)-\mathcal{M}^ih^i(t,z)\geq  -C^i {\lambda \zeta_{max}}e^{\rho T}+C^{\phi}_1:=\eta_1>0. \label{ABB}
\end{align}
Moreover, as $w^i$ is lsc and a supersolution of \reff{HI},  we have
\beq
 w_\kappa^i-\Mc^i w_\kappa^i \geq 0.  \label{BBBB}
\enq
Combining \reff{ABB}, \reff{BBBB} and using the convexity of the operator$\Mc$,  we obtain
\beq\label{premieremaj}
 w_\kappa^{i,\gamma}(t,z)-\Mc^i w_\kappa^{i,\gamma}(t,z)&\geq&(1 - \gamma)(w_\kappa^i(t,z)-\mathcal{M}^iw_\kappa^i(t,z))+ \gamma(h^i(t,z)-\Mc^i h^i(t,z))\nonumber\\
 &\geq &\gamma \eta_1>0
 \enq
From the definition of $h^i$, for all $(t,z)\in  \Ic^i$, we have
\begin{align*}
&  -\frac{\partial h^i}{\partial t} (t,z)- \mathcal{L}h^i (t,z) +\rho^i h^i (t,z)-f^i(z) \\
  &= (\rho^i+\rho) h^i(t,z)-\mu  x \frac{\partial{h^i}}{\partial{x}}(t,z)-\frac{\sigma^2}{2} x^2 \frac{\partial^2{h^i}}{\partial{x}^2}(t,z)-f^i(z)\\
&= (\rho^i+\rho) h^i(t,z)-(2\mu+\sigma^2) B^ix^2 e^{\rho(T-t)}-f^i(z)\\
  &\geq  (\rho^i+\rho) (A^i+B^ix^2+C^i log(y^i+1)+D^i \log(y^j+1))e^{\rho(T-t)}-(2\mu+\sigma^2) B^ix^2 e^{\rho(T-t)}-K^ix\\
  &\geq x^2\left(B^i e^{\rho(T-t)} ( \rho^i+\rho  -2\mu -\sigma^2 )-K^i\right)+\left((\rho^i+\rho) A^i e^{\rho(T-t)}-K^i\right)\\
  &+C^i e^{\rho(T-t)} (\rho^i+\rho) log(y^i+1)+D^i e^{\rho(T-t)} (\rho^i+\rho) log(y^j+1)\\
  &\geq x^2\left(B^i  ( \rho^i+\rho  -2\mu -\sigma^2 )-K^i\right)+\left((\rho^i+\rho) A^i -K^i\right)\\
  &+C^i  (\rho^i+\rho) log(y^i+1)+D^i (\rho^i+\rho) log(y^j+1),
\end{align*}
where the first and the second inequalities are obtained since $f^i(z)\leq K^i x\leq K^i(1+x^2)$.
Thanks to the choice of $A^i$, $B^i$ and $\rho$  (See \reff{param}),
 we have
 \begin{align}
 -\frac{\partial h^i}{\partial t}(t,z) - \mathcal{L} h^i(t,z) -\rho^i h^i(t,z)-f^i(z) \geq (\rho^i+\rho) A^i-K^i=:\eta_2 >0.
 \end{align}
Thus we obtain
\beq
\label{CCCC}
-\Dt{w_\kappa^{i,\gamma}}-\Lc w_\kappa^{i,\gamma} +\rho^i w_\kappa^{i,\gamma}- f^i \geq \gamma\eta_2 \mbox{ in the viscosity sense.}
\enq
Combining \reff{premieremaj} and \reff{CCCC}, we have
\beq\label{visstricsuper}
min\lbrace -\Dt{w_\kappa^{i,\gamma}}-\Lc w_\kappa^{i,\gamma} +\rho^i  w_\kappa^{i,\gamma}- f^i, w_\kappa^{i, \gamma}-\Mc^i w_\kappa^{i,\gamma}\rbrace \geq \gamma \eta\quad \hbox{in} \quad {\Ic}^i,
\enq
where $\eta:=\min \lbrace  \eta_1, \eta_2 \rbrace    > 0$ and  the QVI \reff{visstricsuper} should be interpreted in the viscosity sense.\\
Let  $(t,z) \in \bar \Ic^i$ and $\gamma \in (0,1).$ From the definition of $h^i$, we have:
\beqs
& &h^i(t,z) - \Hc_\kappa^ih^{i}(t,z)\\
&=&h^i(t,z)-\left(  h^i(t,x,y^i,\Gamma^j(z^j,\zeta^j))-\kappa \right) \\
&=&h^i(t,z)-\left(  A^i+B^ix^2+C^i \log(y^i +1)+D^i \log(y^j e^{\lambda \zeta^j}+1)\right)e^{\rho(T-t)} + \kappa\nonumber
\enqs
From the choice of $D^i$ (See \reff{param}), we deduce that:
\begin{align}\label{syst2membre1}
  h^i(t,z)-\Hc_\kappa^ih^i(t,z)\geq  -D^i {\lambda \zeta_{max}}e^{\rho T}+\kappa:=\eta^{'}_1>0.
\end{align}
On the other hand, for $(t,z) \in \bar \Ic^i$, we have:
\beqs
& &h^i(t,z) - \Mc^i\Hc_\kappa^ih^{i}(t,z)\\
&= &h^i(t,z)-\left(  h^i(t,x,\Gamma^i(z^i,\zeta^i),\Gamma^j(z^j,\zeta^j))-\phi^i(y^i,\zeta^i)-\kappa \right) \\
&=&h^i(t,z)-\left(  A^i+B^ix^2+C^i \log(y^i e^{\lambda \zeta^i}+1)+D^i \log(y^j e^{\lambda \zeta^j}+1)\right)e^{\rho(T-t)}+\phi^i(y^i,\zeta^i)+\kappa \nonumber \\
&\geq& - C^i \lambda \zeta_{max}e^{\rho(T-t)}+C^{\phi}_1- D^i \lambda \zeta_{max}e^{\rho(T-t)}+\kappa,\nonumber
\enqs
From the choices of $C^i$ and  $D^i$ (See \reff{param}), we deduce that:
\begin{align}\label{syst2membre2}
  h^i(t,z)-\Mc^i\Hc_\kappa^ih^i(t,z)\geq  -C^i {\lambda \zeta_{max}}e^{\rho T}+C^{\phi}_1-D^i {\lambda \zeta_{max}}e^{\rho T}+\kappa:=\eta^{'}_2>0.
\end{align}
From \reff{syst2membre1} and \reff{syst2membre2}, we deduce that:
\beqs
\Min\{h^i-\Mc^i\Hc_\kappa^ih^{i},h^i-\Hc_\kappa^ih^{i}\}\geq \eta^{'}:=\Min\{\eta^{'}_1,\eta^{'}_2\}\quad \hbox{in} \quad\overline{\Ic}^i.
\enqs
As $w_\kappa^i$ is lsc and a viscosity supersolution of \reff{HIj}, then
\beqs
 \Min\{w_\kappa^i-\Mc^i\Hc_\kappa^iw_\kappa^{i},w_\kappa^i-\Hc_\kappa^iw_\kappa^{i}\}\geq 0\quad \hbox{in} \quad\overline{\Ic}^i.
\enqs
We deduce that:
\beq\label{visstricsuper2}
 \Min\{w_\kappa^{i,\gamma} -\Mc^i\Hc_\kappa^iw_\kappa^{i,\gamma},w_\kappa^{i,\gamma}\Hc_\kappa^iw_\kappa^{i,\gamma}\}\geq \gamma \eta^{'}\quad \hbox{in} \quad\overline{\Ic}^i.
\enq
The result follows from \reff{visstricsuper} and \reff{visstricsuper2}.
\ep\\
\\
The following result is a comparison theorem. It states that we can compare a viscosity subsolution to \reff{HIjp}-\reff{HIp} in $[0,T)\times \Sc$ and a viscosity supersolution to \reff{HIjp}-\reff{HIp} in $[0,T)\times\Sc$, provided that we can compare them at the terminal date but also at the corner lines $D_0$ of the solvency boundary. We have to handle with the difficulties coming from the boundary. We adapt the argument in Barles \cite{bar94} (See Theorem 4.5) which needs a smooth boundary.
Akian et al. \cite{akisultak01} proved a comparison theorem where the Hamilton Jacobi Bellman Variational inequality is satisfied up the boundary.
The non regularity at some points of the boundary is studied in Ly Vath et al. \cite{MLP07}.
They proved a comparison theorem on the solvency region except the non regular part of the boundary. In our case, we can not compare a subsolution and a supersolution at the corner lines.
\begin{Theorem} \label{theocompa}
We fix $i\in \{1,2\}$. Suppose that $u_\kappa^i\in USC([0,T)\times\bar \Sc)$ is a viscosity subsolution to  \reff{HIjp}-\reff{HIp} in $[0,T)\times \bar \Sc$ and
    $w_\kappa^i\in LSC([0,T)\times\bar \Sc)$ is a viscosity supersolution to  \reff{HIjp}-\reff{HIp} in $[0,T)\times \Sc$. We assume that $u^i_\kappa$ and $w^i_\kappa$ satisfy the growth condition \reff{growth} and

\begin{eqnarray}\label{corners}
  u_\kappa^i(t,z)\leq
  \Liminf_{(t^{'},z^{'})\longrightarrow(t,z)}
  w_\kappa^i(t^{'},z^{'}),
  \,\,\forall (t,z)\in [0,T)\times D_0,
\end{eqnarray}
\begin{eqnarray}\label{closuredomain}
u_\kappa^i(T,z):=\Limsup_{\tiny\begin{array}{c}
        (t,z') \rightarrow (T,z)  \\
        t < T, z'\in\Sc
       \end{array}}u_\kappa^i(t^{'},z^{'}) \leq w_\kappa^i(T,z):=\Liminf_{\tiny\begin{array}{c}
        (t,z') \rightarrow (T,z)  \\
        t < T, z'\in\Sc
       \end{array}}w_\kappa^i(t^{'},z^{'})\,\,\forall z\in \bar \Sc.
\end{eqnarray}
Then,
\beqs
u_\kappa^i\leq w_\kappa^{i}\,\,\mbox{ on } [0,T]\times\cal S.
\enqs
\end{Theorem}
\noindent\underline{\emph{Proof:}}
In order to prove the comparison principle, it suffices to show that for all $\gamma  \in (0,1):$
\beqs
\sup_{z \in \Sc}(u_\kappa^i-w_\kappa^{i,\gamma}) \leq 0,
\enqs
where $w_\kappa^{i,\gamma}$ is a strict supersolution satisfying \reff{strictsuper}.
We obtain the required result by letting $\gamma$ to 0. We argue by contradiction. We suppose
that
\beq
\mu^i:=\sup_{z \in \Sc}(u_\kappa^i-w_\kappa^{i,\gamma})> 0.
\enq
We define $w^i_\kappa$ for $i\in \{1,2\}$ on $[0,T)\times\partial\Sc$ by~:
\beq \label{defwbarsc}
w_\kappa^i(t,z) &=& \liminf_{\tiny\begin{array}{c}
      (t',z')\rightarrow (t,z) \\
      (t',z') \in \bar [0,T)\times\Sc
      \end{array} } w_\kappa^i(t',z'), \;\;\; \forall (t,z) \in [0,T)\times\partial\Sc.
    \enq
As $u_\kappa^i$ is usc and  $w_\kappa^{i,\gamma}$ is lsc, then $u_\kappa^i- w_\kappa^{i,\gamma}$ is usc.
By the choice of the strict supersolution, we have $\Lim_{|x|,|y^i|,|y^j|\longrightarrow \infty} u^i(t,z) - w^{i,\gamma}(t,z)=-\infty$.
Hence, the set $\Argmax_{[0,T]\times \bar \Sc }(u_\kappa^i- w_\kappa^{i,\gamma})$ is non empty. As
$u_\kappa^i(T,z) - w_\kappa^{i,\gamma}(T,z)\leq 0$  and $u_\kappa^i(t,z) - w_\kappa^{i,\gamma}(t,z)\leq 0$ for all $(t,z)\in [0,T)\times D_0$,
there exists an open set $\Kc$ with compact closure $\bar \Kc$  s.t.
\beqs
\Argmax_{[0,T]\times \bar \Sc }(u_\kappa^i- w_\kappa^{i,\gamma})\subset [0,T)\times ( (\bar\Sc\setminus D_0) \cap \Kc).
\enqs
We choose $(t_0,z_0)\in [0,T)\times (\bar\Sc\setminus D_0 \cap \Kc)$ s.t.
\beq\label{teta}
\mu^i = u_\kappa^i(t_0,z_0) - w_\kappa^{i,\gamma}(t_0,z_0).
\enq
We distinguish two cases:\\
\noindent $\bullet$ {\bf \it Case 1.}~:$z_0\in \partial \Sc \setminus D_0 \cap \Kc$.
From \reff{defwbarsc}, there exists a sequence $(t_n,z_n)_{n\geq 1}$ in $[0,T)\times(\Sc\cap \Kc)$ converging to
$(t_0,z_0)$  s.t.  $w_\kappa^{i,\gamma}(t_n,z_n)$ tends to $w_\kappa^{i,\gamma}(t_0,z_0)$ when $n$ goes to infinity.  We then set
  $\beta_n$ $:=$ $|t_n-t_0|$, $\eps_n$ $:=$ $|z_n-z_0|$ and consider the function  $\Phi_n^i$
  defined on $[0,T]^2\times(\bar\Sc\cap\bar\Kc)^2$ by~:
\beq
\Phi_n^i(t,t',z,z') &=& u_\kappa^i(t,z) - w_\kappa^{i,\gamma}(t',z') - \varphi^{i}_n(t,t',z,z') \label{defpsi} \\
\varphi_n^i(t,t',z,z') &=&  |t-t_0|^2 +  |z-z_0|^4   +   \frac{|t-t'|^2}{2\beta_n} +  \frac{|z-z'|^2}{2\eps_n} +
\left(\frac{d(z')}{d(z_n)} - 1\right)^4. \nonumber
\enq
Here $d(z)$ denotes the distance from $z$ to $\partial\Sc$.
It is known that for $z_0\notin D_0$, there exists an open neighborhood of $z_0$ in which the distance $d(.)$
is twice differentiable with bounded derivatives. Such property fails when $z_0$ belongs to the corner lines i.e. $z_0\in D_0$ (See \cite{giltru77}). \\
Since $\Phi_n^i$ is usc on the compact set
$[0,T]^2\times(\bar\Sc\cap\bar\Kc)^2$, there exists
$(\hat t_n,\hat t_n^{'},\hat z_n,\hat z_n^{'})$ $\in$ $[0,T]^2\times(\bar\Sc\cap\bar\Kc)^2$
that attains its maximum  on $[0,T]^2\times(\bar\Sc\cap\bar\Kc)^2$. We define
\beqs
\mu_n^i &:= & \sup_{[0,T]^2\times(\bar\Sc\cap \Kc)^2} \Phi_n^i(t,t',z,z') \; = \; \Phi_n^i(\hat t_n,\hat t_n^{'},\hat z_n,\hat z_n^{'}).
\enqs
Moreover, there exists a subsequence, also
denoted $(\hat t_n,\hat t_n^{'},\hat z_n,\hat z_n^{'})_{n\geq 1}$, converging to $(\hat t_0,\hat t_0^{'},\hat z_0,\hat z_0^{'})$
$\in$
$[0,T]^2\times(\bar\Sc\cap\bar\Kc)^2$.
Since  $\Phi_n^i(t_0,t_n,z_0,z_n)$ $\leq$ $\mu^i_n=\Phi_n^i(\hat t_n,\hat t_n^{'},\hat z_n,\hat z_n^{'})$, we have~:
\beq
& & u_\kappa^i(t_0,z_0) - w_\kappa^{i,\gamma}(t_n,z_n) - \frac{1}{2}\left(|t_n-t_0| + |z_n-z_0|\right) \label{interviscouni0} \\
& \leq & \mu^i_n \; = \;  u_\kappa^i(\hat t_n,\hat z_n) - w_\kappa^{i,\gamma}(\hat t_n',\hat z_n') -
\left(|\hat t_n-t_0|^2 + |\hat z_n-z_0|^4 \right) - R_n \label{interviscouni1} \\
& \leq &   u_\kappa^i(\hat t_n,\hat z_n) - w_\kappa^{i,\gamma}(\hat t_n',\hat z_n') -
\left(|\hat t_n-t_0|^2 + |\hat z_n-z_0|^4 \right), \label{interviscouni2}
\enq
where we set
\beqs
R_n &:=& \frac{|\hat t_n-\hat t_n'|^2}{2\beta_n} +  \frac{|\hat z_n-\hat z_n'|^2}{2\eps_n} + \left(\frac{d(\hat z_n')}{d(z_n)} - 1\right)^4.
\enqs
As $u_\kappa^i$, $w_\kappa^{i,\gamma}$ are bounded on $[0,T]\times\bar\Sc\cap\bar\Kc$, then inequality \reff{interviscouni1} implies
the boundedness of
$(R_n)_{n\geq 1}$. This yields that, there exists a subsequence, also denoted $(R_n)_{n\geq 1}$, which is convergent.
As $\beta_n$ and $\eps_n$ go to $0$ when $n$ goes to infinity, we must have
\beq
\hat t_0 \,=\,\hat t_0' \mbox{ and }\hat z_0\,=\,\hat z_0'.
\enq
As $u^i$ is upper semi-continuous and $w_\kappa^{i,\gamma}$ is lower semi-continuous, then, by sending $n$ to infinity into \reff{interviscouni0} and
\reff{interviscouni2}, we obtain
\beqs
\mu^i=u_\kappa^i(t_0,z_0)-w_\kappa^{i,\gamma}(t_0,z_0)\leq u_\kappa^i(\hat t_0,\hat z_0)-w_\kappa^{i,\gamma}(\hat t_0,\hat z_0)-|\hat t_0-t_0|^2-|\hat z_0-z_0|^4\leq \mu^i,
\enqs
where the last inequality is deduced from the definition of $\mu^i$. It yields that:
\beq \label{teta}
\hat t_0 \; =  \hat t_0' \;  = \;  t_0, & &
\hat z_0 \; =  \hat z_0' \;  = \;  z_0.
\enq
Using again \reff{interviscouni0}-\reff{interviscouni1}-\reff{interviscouni2} and sendig $n$ to infinity, we deduce that:
\beq\label{etaieta}
\Lim_{n\rightarrow \infty}\mu^i_n =\mu^i.
\enq
From equality \reff{interviscouni1}, we have $\Lim_{n\rightarrow \infty}\mu^i_n=\mu^i-\Lim_{n\rightarrow \infty} R_n$, then we deduce that:
\beq\label{disti}
\Lim_{n\rightarrow \infty} R_n=0
\enq
In particular, for $n$ large enough , we have $\hat t_n$, $\hat t_n'$ $<$ $T$ (since $t_0$ $<$ $T$). Besides as $\frac{d(\hat z_n')}{d(z_n)}\rightarrow 1$ when $n$ goes to infinity, we have $d(\hat z_n')\geq \frac{d(z_n)}{2}$ for $n$ large enough. It yields that $\hat z_n'\in \Sc$ and consequently $\hat z_n$ $\in$ $\Sc$. \\
Applying Ishii's lemma (see Theorem 8.3 in \cite{Crand}) to
$(\hat t_n,\hat t_n^{'},\hat z_n,\hat z_n^{'})$ $\in$ $[0,T)\times[0,T)\times(\Sc\cap\Vc_0)\times(\Sc\cap\Vc_0)$ that attains the maximum of
$\Phi_n^i$ in \reff{defpsi}, where $\Vc_0$ is a neighborhood of $z_0$.
We get the existence of two $3 \times 3$ symmetric matrices ${M}$ and ${N}$ s.t.:
\beqs
(s_0,s,M)\quad &\in \quad  \bar J^{2,+}u_\kappa^i(\hat t_n,\hat z_n),\\
(q_0,q,N)\quad &\in \quad  \bar J^{2,-}w_\kappa^{i,\gamma}(\hat t_n',\hat z_n'),
\enqs
where
\beq
\label{SQ}
\begin{pmatrix}
M & 0 \\
0 & -N
\end{pmatrix}
\leq D^2_{z,z'}\varphi_{n}^i(\hat t_n,\hat t_n^{'},\hat z_n,\hat z_n^{'})+{\varepsilon_n}
(D^2_{z,z'}\varphi_{n}^i(\hat t_n,\hat t_n^{'},\hat z_n,\hat z_n^{'}))^2 \enq
\beqs
s_0&=& \Dt{\varphi_n^i}(\hat t_n,\hat t_n^{'},\hat z_n,\hat z_n^{'})=2(\hat{t}_n - t_0) + \frac{(\hat{t}_n -\hat{t}^{'}_n)}{\beta_n},\\
q_0&=&-\frac{\partial \varphi_n^i}{\partial t^{'}}(\hat t_n,\hat t_n^{'},\hat z_n,\hat z_n^{'})=\frac{(\hat{t}_n -\hat{t}^{'}_n)}{\beta_n},\\
s&=&(s_k)_{1\leq k\leq 3}
=D_z{\varphi_n^i}
(\hat t_n,\hat t_n^{'},\hat z_n,\hat z_n^{'})
=4(\hat{z}_n - z_0)|\hat{z}_n - z_0|^2 + \frac{(\hat{z}_n - \hat{z}^{'}_n)}{\eps_n},\\
q &=&(q_k)_{1\leq k\leq 3}=-D_{ z^{'}}\varphi_n^i(\hat t_n,\hat t_n^{'},\hat z_n,\hat z_n^{'})=\frac{(\hat{z}_n - \hat{z}^{'}_n)}{\eps_n}
- \frac{4}{d(z_n)}  \left( \frac{d(\hat{z}^{'}_n)}{d(z_n)} - 1 \right)^3D d(\hat{z}^{'}_n),
\enqs
and
\beqs
 &D^2_{z,z'}\varphi_n^i(\hat t_n,\hat t_n',\hat z_n,\hat z_n')=
 \begin{pmatrix}
    \frac{1}{\varepsilon_n}I_3+P_n&-\frac{1}{\varepsilon_n}I_3\\
    -\frac{1}{\varepsilon_n}I_3& \frac{1}{\varepsilon_n}I_3+Q_n

\end{pmatrix},
 \enqs
 where
 \beqs
P_n &=&    4|\hat{z}_n - z_0|^2 I_3 + 8 (\hat{z}_n - z_0)(\hat{z}_n - z_0)\trans, \\
Q_n  &=&  12 \left(\frac{d(\hat{z}'_n)}{d(\hat{z}_n)} - 1 \right)^2 \frac{1}{d(\hat{z}_n)^2}D d(\hat{z}'_n)  D d(\hat{z}'_n)\trans
        + \frac{4}{d(z_n)} \left(\frac{d(\hat{z}'_n)}{d(\hat{z}_n)} - 1 \right)^3 D^2 d(\hat{z}'_n).
        \enqs
        Here $\trans$ denotes the transpose operator.
By writing the viscosity subsolution property of $u_\kappa^i$ and the viscosity supersolution property of $w_\kappa^{i,\gamma},$
we have the following inequalities:
\beq
\label{LK}
&\min& \lbrace -s_0-\mu \hat x_n s_1-\frac{\sigma^2 {\hat x_n}^2}{2}M_{11} - f^i(\hat z_n)+\rho^iu_\kappa^i(\hat t_n,\hat z_n),\\
& &u_\kappa^i(\hat t_n,\hat z_n)-\Mc^i u_\kappa^i(\hat t_n,\hat z_n) \rbrace \leq 0  \quad \hbox{if} \quad(\hat t_n,\hat z_n)\in {\Ic}^i,\nonumber \\
\label{LKK}
&\min& \lbrace -q_0-\mu {\hat x}_n^{'} q_1-\frac{\sigma^2 {\hat x}_n^{'2}}{2}N_{11} - f^i(\hat z_n^{'})+\rho^iu_\kappa^i({\hat t}_n^{'},{\hat z}^{'}_n),\\
& &w_\kappa^{i,\gamma}({\hat t}_n^{'},{\hat z}^{'}_n)-\Mc^i w_\kappa^{i,\gamma}({\hat t}_n^{'},{\hat z}^{'}_n) \rbrace \geq \gamma\eta \quad \hbox{if} \quad({\hat t}_n^{'},{\hat z}^{'}_n)\in {\Ic}^i,\nonumber\\
&\min& \lbrace u_\kappa^i(\hat t_n,\hat z_n)-\Mc^i \Hc_\kappa^iu_\kappa^i(\hat t_n,\hat z_n),\label{LK3}\\
&&u_\kappa^i(\hat t_n,\hat z_n)-\Hc_\kappa^i u_\kappa^i(\hat t_n,\hat z_n) \rbrace \leq 0 \quad \hbox{if}
\quad(\hat t_n,\hat z_n)\in \overline{\Ic}^i,\nonumber \\
&\min&  \lbrace w_\kappa^{i,\gamma}({\hat t}_n^{'},{\hat z}^{'}_n)- \Mc^i\Hc_\kappa^iw_\kappa^{i,\gamma}({\hat t}_n^{'},{\hat z}^{'}_n),\label{LK4}\\
&&w_\kappa^{i,\gamma}({\hat t}_n^{'},{\hat z}^{'}_n)-\Hc_\kappa^iw_\kappa^{i,\gamma} \rbrace \geq \gamma\eta^{'} \quad \hbox{if} \quad({\hat t}_n^{'},{\hat z}^{'}_n)\in \overline{\Ic}^i. \nonumber
\enq
\noindent We then distinguish four subcases:\\
\textbf{Subcase 1:} $u_\kappa^i({\hat t}_n,{\hat z}_n)- \Mc^i u_\kappa^{i}({\hat t}_n,{\hat z}_n)\leq 0$ in \eqref{LK}. \\
From the definition of $(\hat t_n,\hat t_n',\hat z_n,\hat z_n'),$ we have:
\beqs
\mu^i_n
&\leq& u_\kappa^i({\hat t}_n,{\hat z}_n) - w_\kappa^{i,\gamma}({\hat t}_n^{'},{\hat z}^{'}_n)\\
&\leq&   \Mc^i u_\kappa^{i}({\hat t}_n,{\hat z}_n)-\Mc^i w_\kappa^{i,\gamma}({\hat t}_n^{'},{\hat z}_n^{'})-\gamma\eta,
\enqs
where the last inequality is deduced from \reff{LKK}.
Now, letting $n$ going to $\infty,$ and using \reff{etaieta}, we obtain:
 we obtain~:
\beqs
\mu^i  & \leq &  \limsup_{n\rightarrow\infty} \Mc^i u_\kappa^i(\hat t_n,\hat z_n) -  \liminf_{n\rightarrow\infty}  \Mc^i w_\kappa^{i,\gamma}(\hat t_n',\hat z_n')-\gamma\eta \\
& \leq & \Mc^i u_\kappa^i(t_0,z_0) - \Mc^i w_\kappa^{i,\gamma}(t_0,z_0) -\gamma\eta,
\enqs
where we used the upper-semicontinuity  of $\Mc^i u_\kappa^i$ and the lower-semicontinuity of  $\Mc^i w_\kappa^{i,\gamma}$ (see Lemma \ref{Lemmatechnical}).
As $u_\kappa^i$ is usc, there exists $\zeta^i \in [\zeta_{min},\zeta_{max}]$ s.t.\\
$\Mc^i u_\kappa^i(t_{0}, z_{0})$ $=$ $u_\kappa^i(t_{0},x_0,\Gamma^i(y^i_{0},\zeta^i),y^j_0)-\phi^i(y^i_{0},\zeta^i)$.
  We then get
\beqs
\mu^i  & \leq & \Mc^i u_\kappa^i(t_0,z_0) - \Mc^i w_\kappa^{i,\delta}(t_0,z_0) - \gamma\eta  \\
& \leq &  u_\kappa^i(t_{0},x_0,\Gamma^i(y^i_{0},\zeta^i),y^j_0)  -  w_\kappa^{i,\delta}(t_{0},x_0,\Gamma^i(y^i_{0},\zeta^i),y^j_0) - \gamma\eta \;  \leq \;  \mu^i - \gamma\eta,
\enqs
which is obviously a contradiction.\\
\textbf{Subcase 2:} $-s_0-\mu \hat x_n s_1-\frac{\sigma^2 {\hat x_n}^2}{2}M_{11} - f^i(\hat z_n)+\rho^iu_\kappa^i(\hat t_n,\hat z_n)\leq 0$ in \eqref{LK}.
From \eqref{LKK}, we have $-q_0-\mu {\hat x}_n^{'} q_1-\frac{\sigma^2 {\hat x}_n^{'2}}{2}N_{11} - f^i(\hat z_n^{'}) +\rho^iu_\kappa^i({\hat t}_n^{'},{\hat z}^{'}_n)\geq \gamma\eta,$
which implies in this case
\beq
&-&(s_0-q_0)-\mu(\hat x_n s_1 -{\hat x}_n^{'} q_1)-\frac{\sigma^2}{2}(M_{11} {\hat x}^2_n-N_{11} {\hat x}^{'2}_n)
- f^i(\hat z_n)+ f^i(\hat z_n^{'})\nonumber\\
&+&\rho^i(u_\kappa^i(\hat t_n,\hat z_n)-u_\kappa^i({\hat t}_n^{'},{\hat z}^{'}_n))\leq -\gamma\eta. \label{QQQQ}
\enq
We have that
\beqs s_0-q_0=2(\hat t_{n}-t_0).\enqs
Since $\hat t_{n}$ goes to $t_0$ when $n$ goes to infinity, we deduce that $s_0-q_0$  goes to zero when $n$ goes to infinity.
The second term of \reff{QQQQ} is expressed as folows:
\beqs
\hat x_n s_1 -{\hat x}_n^{'} q_1=
4\hat x_n(\hat{x}_n - x_0)|\hat{z}_n - z_0|^2 + 4 \frac{(\hat{x}_n - \hat{x}^{'}_n)^2}{\eps_n}
+ 4 {\hat x}_n^{'} \frac{D_xd(\hat{z}^{'}_n)}{d(\hat{z}_n)} \left( \frac{d(\hat{z}^{'}_n)}{d(z_n)} - 1 \right)^3,
\enqs
where $D_xd(\hat{z}^{'}_n)$ is the first component of $Dd(\hat{z}^{'}_n)$, which is known to be continuous on $\bar \Sc\setminus D_0$.
Using  \reff {teta} and \reff{disti},  we deduce that $\hat x_n s_1 -{\hat x}_n^{'} q_1$ goes to 0 when $n$ goes to infinity.
Moreover, from \eqref{SQ}, we have $$\frac{\sigma^2}{2}(M_{11}{\hat x}_n^2-N_{11}{\hat x}_n^{'2})\leq \Upsilon_n $$
where
\beqs \Upsilon_n=\Lambda_n\left( D_{z,z'}^2 \varphi_n(\hat{t}_n,\hat{t}'_n,\hat{z}_n,\hat{z}'_n)
+ \eps_n \left(D_{z,z'}^2 \varphi_n^i(\hat{t}_n,\hat{t}^{'}_n,\hat{z}_n,\hat{z}^{'}_n)\right) ^2\right)\Lambda_n^{\trans},
\enqs
and
\beqs
\Lambda_n=(\frac{\sigma}{\sqrt{2}}{\hat x}_n,0,0,\frac{\sigma}{\sqrt{2}}{\hat x}^{'}_n,0,0)
\enqs
 After some straightforward calculation, we then get~:
\beqs
\Upsilon_n &=&  3\frac{(\hat{x}_n-\hat{x}'_n)^2}{\eps_n} + \Lambda_n\left(
 \left(
    \begin{array}{cc}
        3 P_n & - P_n-Q_n \\
        - P_n-Q_n & 3 Q_n
    \end{array} \right)
+ \eps_n \left(
    \begin{array}{cc}
        P_n^2 & 0 \\
        0 & Q_n^2
    \end{array} \right)
\right) \Lambda_n^{\trans},
\enqs
which converges also to zero from \reff{teta} and \reff{disti}.
We have that $ \Upsilon_n $ goes to zero when $n$ goes to infinity.
Sending $n$ goes to infinity in inequality \reff{QQQQ}, we obtain the required contradiction: $0  \leq -\gamma\eta< 0.$\\
\textbf{Subcase 3:} $u^i({\hat t}_n,{\hat z}_n)- \Hc_\kappa^i u^{i}({\hat t}_n,{\hat z}_n)\leq 0$ in \reff{LK3}. \\
From the definition of $(\hat t_n,\hat t_n',\hat z_n,\hat z_n'),$ we have
\beqs
\mu^i_n
&\leq& u_\kappa^i({\hat t}_n,{\hat z}_n) - w_\kappa^{i,\gamma}({\hat t}_n^{'},{\hat z}^{'}_n)\\
&\leq&   \Hc_\kappa^i u_\kappa^{i}({\hat t}_n,{\hat z}_n)-\Hc_\kappa^i w_\kappa^{i,\gamma}({\hat t}_n^{'},{\hat z}_n^{'})-\gamma\eta^{'},
\enqs
where the last inequality is deduced from \reff{LK4}.
Now, using the upper-semicontinuity  of $\Hc_\kappa^i u_\kappa^i$ and the lower-semicontinuity of  $\Hc_\kappa^i w_\kappa^{i,\gamma}$ (see Lemma \ref{Lemmatechnical} (iii)), letting $n$ going to $\infty,$ and using \reff{etaieta},
 we obtain~:
\beqs
\mu^i  & \leq &  \limsup_{n\rightarrow\infty} \Hc_\kappa^i u_\kappa^i(\hat t_n,\hat z_n) -  \liminf_{n\rightarrow\infty}  \Hc_\kappa^i w_\kappa^{i,\gamma}(\hat t_n',\hat z_n')-\gamma\eta^{'} \\
& \leq & \Hc_\kappa^i u_\kappa^i(t_0,z_0) - \Hc_\kappa^i w_\kappa^{i,\gamma}(t_0,z_0) -\gamma\eta \;\leq\;  \mu^i - \gamma\eta^{'},
\enqs
which is obviously a contradiction.\\
\textbf{Subcase 4:} $u_\kappa^i({\hat t}_n,{\hat z}_n)- \Mc^i \Hc_\kappa^i u_\kappa^{i}({\hat t}_n,{\hat z}_n)\leq 0$ in \eqref{LK}. \\
From the definition of $(\hat t_n,\hat t_n',\hat z_n,\hat z_n'),$ we have
\beqs
\mu^i_n
&\leq& u_\kappa^i({\hat t}_n,{\hat z}_n) - w_\kappa^{i,\gamma}({\hat t}_n^{'},{\hat z}^{'}_n)\\
&\leq&   \Mc^i \Hc_\kappa^i u_\kappa^{i}({\hat t}_n,{\hat z}_n)-\Mc^i \Hc_\kappa^i w_\kappa^{i,\gamma}({\hat t}_n^{'},{\hat z}_n^{'})-\gamma\eta^{'},
\enqs
where the last inequality is deduced from \reff{LK4}.
Now, letting $n$ going to $\infty,$ and using \reff{etaieta}, we obtain:
\beqs
\mu^i  & \leq &  \limsup_{n\rightarrow\infty} \Mc^i  \Hc_\kappa^iu_\kappa^i(\hat t_n,\hat z_n) -  \liminf_{n\rightarrow\infty}  \Mc^i \Hc_\kappa^i w_\kappa^{i,\gamma}(\hat t_n',\hat z_n')-\gamma\eta \\
& \leq & \Mc^i \Hc_\kappa^i u_\kappa^i(t_0,z_0) - \Mc^i  \Hc_\kappa^iw_\kappa^{i,\gamma}(t_0,z_0) -\gamma\eta^{'},
\enqs
where we used the upper-semicontinuity  of $\Mc^i  \Hc_\kappa^iu_\kappa^i$ and the lower-semicontinuity of  $\Mc^i  \Hc_\kappa^iw_\kappa^{i,\gamma}$ (see Lemma \ref{Lemmatechnical}).
As $u_\kappa^j$ is usc, there exists $\zeta^j \in [\zeta_{min},\zeta_{max}]$ s.t.  $\Mc^j u_\kappa^i(t_{0}, z_{0})$ $=$ \\
$u_\kappa^i(t_{0},x_0,y^i_0,\Gamma^j(y^j_{0},\zeta^j))$,
and then $\Hc_\kappa^i u_\kappa^i(t_0,z_0)=u_\kappa^i(t_{0},x_0,y^i_0,\Gamma^j(y^j_{0},\zeta^j))-\kappa$.
As $u_\kappa^i$ is usc, there exists $\zeta^i \in [\zeta_{min},\zeta_{max}]$ s.t.
$\Mc^i \Hc_\kappa^iu_\kappa^i(t_{0}, z_{0})$ $=$ $u_\kappa^i(t_{0},x_0,\Gamma^i(y^i_{0},\zeta^i),\Gamma^j(y^j_{0},\zeta^j))-\phi(y^i_{0},\zeta^i)-\kappa$.
  We then get
\beqs
\mu^i  & \leq & \Mc^i \Hc_\kappa^i u_\kappa^i(t_0,z_0) - \Mc^i \Hc_\kappa^iw_\kappa^{i,\delta}(t_0,z_0) - \gamma\eta  \\
& \leq &  u_\kappa^i(t_{0},x_0,\Gamma^i(y^i_{0},\zeta^i),\Gamma^j(y^j_{0},\zeta^j))  -  w_\kappa^{i,\delta}(t_{0},x_0,\Gamma^i(y^i_{0},\zeta^i),\Gamma^j(y^j_{0},\zeta^j)) - \gamma\eta^{'} \;  \leq \;  \mu^i - \gamma\eta,
\enqs
which is obviously a contradiction.\\
\noindent $\bullet$ {\it Case 2.}~: $z_0\in  \Sc \cap \Kc$
\noindent We consider the function
\beqs
\Phi_n^i(t,z,z') &=& u^i(t,z) - w^{i,\gamma}(t,z') - \varphi^{i,\gamma}_n(t,z,z') \\
\varphi^{i,\gamma}_n(t,z,z')  &=& |t-t_0|^2 + |z-z_0|^4 + \frac{n}{2} |z-z'|^2,
\enqs
for $n$ $\geq$ $1$,  and to take a maximum $(\tilde t_n,\tilde z_n,\tilde z_n')$ of $\Phi_n^i$.  We then show that
the sequence $(\tilde t_n,\tilde z_n,\tilde z_n')_{n\geq 1}$ converges to $(t_0,z_0,z_0)$ as $n$ goes to infinity and we
apply Ishii's lemma to get the required contradiction.

\ep\\
The condition \reff{corners} in the comparison theorem that must be satisfied is not obvious to check for $v^i$, $i\in \{1,2\}$.
To circumvent this difficulty, we introduce the following function:
\begin{eqnarray}\label{candidat}
{\hat v}^i(t, z):=\E [\int_{t}^{ T }e^{-\rho^i(s-t)}f^i(Z^{t,z}_s)ds +e^{-\rho^i(T-t)}g^i(Z^{t,z}_T)],\,\forall \, (t,z) \in [0,T)\times \bar \Sc,
\end{eqnarray}
where $(Z^{t,z}_s)_{t\leq s\leq T}$ is the state process associated with the no impulse strategy.
The next result shows that the value function $\hat v^i$ is continuous
on the part $[0,T)\times D_0$ of the solvency region.
\begin{Proposition}\label{prophatvi}
For all $i\in \{1,2\}$,  we have:
\begin{eqnarray}\label{continuity}
\lim_{(t^{'},z^{'})\longrightarrow (t,z)} \hat v^i(t^{'},z^{'})= \hat v^i(t,z)\, \mbox{ for all }\, (t,z) \in [0,T)\times D_0,
\end{eqnarray}
and ${\hat v}^i$ is a classical solution
  of  \reff{HIjp}-\reff{HIp} in a neighbourhood of $[0,T)\times D_0$.
\end{Proposition}
\noindent\underline{\emph{Proof:}}
We fix $i\in \{1,2\}$ and $(t,z)\in [0,T)\times D_0$.
We consider $ {\mathring B}((t,z),\delta)\cap \bar\Sc$ a neighbourhood of $(t,z)$, where
${\mathring B}((t,z),\delta):=\{(t^{'}, z^{'})\mbox{ s.t. } |z-z^{'}|+|t-t^{'}|<\delta\}$ and $\delta$ is a positive constant.\\
{\bf Step 1:} We prove the continuity property \reff{continuity}.

From the definition of the process $(Z^{t^{'},z^{'}}_s)_{t^{'}\leq s\leq T}$ with no impulse strategy, we have:
\begin{eqnarray*}
 & Z^{t^{'},z^{'}}_s&=(X^{t^{'},x^{'}}_s,Y^{1,t^{'},y^{'1}}_s,Y^{2,t^{'},y^{'2}}_s)=(X^{t^{'},x^{'}}_s,y^{'1},y^{'2} )
  =(x^{'}e^{(\mu-\frac{\sigma^2}{2})(s-t^{'})+\sigma(W_s-W_{t^{'}})},y^{'1},y^{'2})\\&\longrightarrow&(X^{t,x}_s,y^{1},y^{2} ) dt\otimes d\P,
\mbox{ when }
(t^{'},z^{'})\longrightarrow (t,z).
\end{eqnarray*}
We fix $p>1$. From Remark \ref{rborne}, for $h^i=f^i,g^i$, we have $|h^i(z)|\leq C^i(1+x)$, which implies
\beq\label{majfg}
 & &\E [\left(\int_{t^{'}}^{ T  }f^i(Z^{t^{'},z^{'}}_s)ds
   +g^i(Z^{t^{'},z^{'}}_T)\right)^p]\\
 &\leq& C\left(\E [\left(\int_{t^{'}}^{ T  }f^i(Z^{t^{'},z^{'}}_s)ds\right)^p]
 +\E [g^i(Z^{t^{'},z^{'}}_T)^p]\right)\nonumber\\
 &\leq& C\E [\int_{t^{'}}^{ T  }f^i(Z^{t^{'},z^{'}}_s)^pds](T-t^{'})^{p-1}
 +\E [g^i(Z^{t^{'},z^{'}}_T)^p])\nonumber\\
 &\leq & C(1+T^{p})(1+\E[\Sup_{s\in [t^{'},T]}|X_s^{t^{'},x^{'}}|^{p}]),\nonumber
 \enq
where the second inequality is obtained by using H\"{o}lder inequality and $C$ is a generic constant independent of $(t^{'},z^{'})$.
 From the definition of the wholesale price and for $(t^{'},z^{'})\in {\mathring B}((t,z),\delta)$, we have
 \beq\label{majx}
 \E[\Sup_{s\in [t^{'},T]}|X_s^{t^{'},z^{'}}|^{p}]\leq C|x^{'}|^{p}\leq C(1+|x|^{p}).
 \enq
 From \reff{majfg}-\reff{majx}, we deduce
 \beqs
 \E [\left(\int_{t^{'}}^{ T  }f^i(Z^{t^{'},z^{'}}_s)ds
   +g^i(Z^{t^{'},z^{'}}_T)\right)^p] \leq C(1+|x|^{p}),
 \enqs
 where $C$ is a positive constant which depends on $\delta$. This shows the boundedness of \\
 $\left(\int_{t^{'}}^{ T  }f^i(Z^{t^{'},z^{'}}_s)ds
   +g^i(Z^{t^{'},z^{'}}_T)\right)$ in $L^p(\mathbb{P})$ for $p>1$,
 which implies the uniform integrability of
 $\left(\int_{t^{'}}^{ T  }f^i(Z^{t^{'},z^{'}}_s)ds +g^i(Z^{t^{'},z^{'}}_T)\right)_{t^{'},z^{'}}$.
 It yields that ${ \hat v}^i(t^{'}, z^{'})\longrightarrow { \hat v}^i(t, z)$ when $(t^{'},z^{'})\longrightarrow (t,z)$.\\
{\bf Step 2:} We show that ${\hat v}^i$ is the solution of the IQV
\reff{HIj} on $ {\mathring B}((t,z),\delta^{'})\cap \bar \Sc$ for some $\delta^{'}$. From Equation \reff{candidat},
and since the process $Z^{t^{'},z^{'}}$ is Log-normal, then ${\hat v^i}$ is regular and we have:
\begin{eqnarray*}
-\Dt{\hat v^i}(t^{'}, z^{'})-\Lc {\hat v}^i(t^{'}, z^{'}) +\rho^i {\hat v}^i(t^{'}, z^{'})- f^i(z^{'})=0,
\end{eqnarray*}
on $ {\mathring B}((t,z),\delta)\cap \bar \Sc$. It remains to prove that ${\hat v}^i(t^{'}, z^{'})> \Mc^i{\hat v}^i(t^{'}, z^{'})$
in a neighbourhood of $(t,z)$. We argue by contradiction.
We assume that for all $\delta>0$, there exists $(t^{''}, z^{''})\in {\mathring B}((t,z),\delta)\cap \bar \Sc$
such that  ${\hat v}^i(t^{''}, z^{''})\leq \Mc^i{\hat v}^i(t^{''}, z^{''})$.
Sending $\delta$ towards 0, using Lemma \ref{Lemmatechnical}, and the continuity of $\hat v^i$ w.r.t $(t,z)$, we obtain:
\begin{eqnarray}\label{nointer}
{\hat v^i}(t, z)\leq \Limsup_{(t^{''}, z^{''})\longrightarrow (t, z)}\Mc^i{\hat v^i}(t^{''}, z^{''})\leq \Mc^i{\hat v^i}(t, z).
\end{eqnarray}
On the other hand, as $(t,z)\in [0,T)\times D_0$, we have
\beqs
\hat v^i(t,z)=
\begin{cases}
     0& \quad \text{ if } (t,z)\in [0,T)\times\partial^{y^1}{\cal S}\cup \partial^{y^2}{\cal S}, \\
     -\frac{x}{2}(\frac{e^{(\mu-\rho^i)(T-t)}-1}{\mu-\rho^i}+e^{(\mu-\rho^i) (T-t)}) & \quad \text{ if }  (t,z)\in [0,T)\times\partial^{x}{\cal S}.
  \end{cases}
\enqs
By straight forward computation, we have ${ \hat  v}^i(t, z)>\Mc^i{\hat  v^i}(t, z)$,
which means that there is no intervention in $(t,z)$, and so the contradiction is obtained.
Symmetrically, In a neighbourhood of $(t,z)\in [0,T)\times D_0$, we have ${ \hat  v}^j(t^{'}, z^{'})>\Mc^j{\hat  v^j}(t^{'}, z^{'})$
for $j\neq i$, which means, against the no intervention startegy of player $i$, the best reponse of player $j$ is also not to make an intervention.
This shows that there exists  a neighborhood of $[0,T)\times D_0$ which is included in the continuation region and ${\hat v}^i$ is a regular solution
of \reff{HIp} in this neighbourhood.\\
  \ep
\begin{Proposition}\label{continuituboundary}
(i)For all $i\in \{1,2\}$, we have:
\begin{eqnarray}\label{clthcom3}
 \overline v_\kappa^i(t,z)\leq
  \Liminf_{(t^{'},z^{'})\longrightarrow(t,z)}
  \underline v_\kappa^i(t^{'},z^{'}),
  \,\,\forall (t,z)\in [0,T)\times D_0,
\end {eqnarray}
(ii)For all $i\in \{1,2\}$, we have:
\begin{eqnarray}\label{clthcom4}
 \overline v_\kappa^i(t,z)=
  \underline v_\kappa^i(t,z)=v_\kappa^i(t,z),
  \,\,\forall (t,z)\in [0,T)\times D_0,
\end {eqnarray}
\end{Proposition}
\noindent\underline{\emph{Proof:}} {\bf Step 1:}
We proved that the value function $v_\kappa^i$ is a viscosity subsolution
  of \reff{HIjp}-\reff{HIp} in $[0, T) \times \bar \Sc$, then
for all $({t},{z})\in [0, T) \times D_0$ and
$\varphi^i\in C^{1,2}([0, T) \times \bar \Sc)$ s.t. $(\overline{v}_\kappa^i-\varphi^i)({t},{z}) = 0$
and $({t},{z})$ is a maximum
of $\overline{v}_\kappa^i-\varphi^i$  on $[0, T) \times D_0$, we have
\beqs\label{Visc}
\min\{-\Dt{\varphi^i}({t},{z})-\Lc\varphi^i ({t},{z})+\rho^i \varphi^i({t},{z})- f^i({z}),
\overline{v}_\kappa^i({t},{z})-\Mc^i\overline{v}_\kappa^i({t},{z})\}\leq0 \, \hbox{in} \, [0, T) \times D_0.
  \enqs
 As $(t,z)\in [0, T) \times D_0$, by straight forward computation, we have ${ v}_\kappa^i(t, z)>\Mc^i{  v^i_\kappa}(t, z)$,
  where $v^i_\kappa$ is given by equation \reff{clkappa}, which means that there is no intervention in $(t,z)$. It yields that
\begin{eqnarray}\label{equvar1}
-\Dt{\varphi^i}(t, z)-\Lc \varphi^i(t, z) +\rho^i \varphi^i(t, z)- f^i(z)\leq 0,\,\mbox{ in }[0,T)\times D_0.
\end{eqnarray}
For the terminal condition, we take
\begin{eqnarray}\label{equvar3}
 {\varphi^i}(T, z)\leq g^i(z).
\end{eqnarray}
Tow cases are possible:\\
$\star$ First case: $z\in \partial^{y^j}{\cal S}\cup  \partial^{y^i}{\cal S}$, inequalities \reff{equvar1}-\reff{equvar3} become
\begin{eqnarray*}
-\Dt{\varphi^i}(t, z)+\rho^i \varphi^i(t, z)\leq 0,\,\mbox{ in }[0,T)\times  \partial^{y^j}{\cal S}\cup  \partial^{y^i}{\cal S},
\end{eqnarray*}
and
\begin{eqnarray*}
 {\varphi^i}(T, z)\leq 0=v_\kappa^i(T,z)=\hat v^i(T,z).
\end{eqnarray*}
which implies ${\varphi^i}(t, z)\leq 0= \hat v^i(t,z)$ for all $(t,z)\in [0,T)\times \partial^{y^j}{\cal S}\cup  \partial^{y^i}{\cal S}$\\
$\star$ Second case: $z\in \partial^{x}{\cal S}$, inequalities \reff{equvar1}-\reff{equvar3} become
\begin{eqnarray}\label{equvar4}
-\Dt{\varphi^i}(t, z)-\Lc \varphi^i(t, z)+\rho^i \varphi^i(t, z)-f(t,z)\leq 0,\,\mbox{ in }[0,T)\times \partial^{x}{\cal S}\setminus{(0,0,0)},
\end{eqnarray}
and
\begin{eqnarray}\label{equvar5}
 {\varphi^i}(T, z)\leq 0.
\end{eqnarray}
For the boundary condition, we take
\begin{eqnarray}\label{equvar6}
 {\varphi^i}(t, 0)\leq 0=v_\kappa^i(t,0)=\hat v^i(t,0),\,\mbox{ in }[0,T).
\end{eqnarray}
On the other hand, $\hat v^i$ satisfies \reff{equvar4}-\reff{equvar6} with equalities.
By classical comparison theorem, we deduce that
\begin{eqnarray}\label{clthcom1}
\overline{v}_\kappa^i(t,z)\leq \varphi^i(t,z)\leq  \hat v^i(t,z),\,\mbox{ in }[0,T)\times D_0.
\end {eqnarray}
{\bf Step 2:}
By definition of $\hat v^i$ and  since against the no intervention strategy of player $i$, the best response of player $j$ is also not to make an intervention, we have $ \hat  v^i(t,z)\leq  v_\kappa^i(t,z)$ for all $(t,z)$ in a neighborhood of $[0,T)\times D_0$.
From Proposition \ref{prophatvi}, the function $\hat v^i$ is continuous. It yields that:
\begin{eqnarray*}
 \hat  v^i(t^{'},z^{'})\leq  \underline v_\kappa^i(t^{'},z^{'})\mbox { for all }(t^{'},z^{'})\mbox { in a neighborhood of } (t,z)\in [0,T)\times D_0.
\end {eqnarray*}
Using again the continuity property of $\hat v^i$ (See Proposition \ref{prophatvi}), and since $\underline v_\kappa^i$ is lsc, we obtain:
\begin{eqnarray}\label{clthcom2}
 \hat  v^i(t,z)\leq  \Liminf_{(t^{'},z^{'})\longrightarrow(t,z)}\underline v_\kappa^i(t^{'},z^{'})=\underline v_\kappa^i(t,z) \mbox { for all }(t,z)\in [0,T)\times D_0.
\end {eqnarray}
From inequalities \reff{clthcom1}, \reff{clthcom2} and , we deduce inquality \reff{clthcom3}
and the continuity property of $v_\kappa^i$ in the boundary \reff{clthcom4}.
\ep\\
\\
Finally, combining the previous results, we obtain the following PDE characterization of the value function.
\begin{Corollary}
The value function $v_\kappa^i$ is continuous on $[0,T)\times\Sc$ and is the unique (in $[0,T)\times\Sc$) constrained viscosity solution to
the system of QVIs \reff{HIjp}-\reff{HIp} lying in the class of functions with linear growth in $x$ uniformly in $(t,y^i,y^j)$ and satisfying the boundary condition~:
\beqs
\lim_{(t',z')\rightarrow (t,z)} v_\kappa^i(t',z')=
\begin{cases}
     0& \quad \text{if } (t,z)\in [0,T)\times\partial^{y^1}{\cal S}\cup \partial^{y^2}{\cal S}, \\
     -\frac{x}{2}(\frac{e^{(\mu-\rho^i)(T-t)}-1}{\mu-\rho^i}+e^{(\mu-\rho^i) (T-t)}) & \quad \text{if }  (t,z)\in [0,T)\times\partial^{x}{\cal S} ,
  \end{cases}
\enqs
and the terminal condition
\beqs
v_\kappa^i(T,z)&=& g^i(z), \;\;\; \forall z \in \bar\Sc.
\enqs

\end{Corollary}
{\bf Proof.}
We have $\bar v_\kappa^i$ is an usc viscosity subsolution to \reff{HIjp}-\reff{HIp} in $[0, T) \times \bar \Sc$ and $\underline v^i_\kappa$ is a
lsc viscosity supersolution to \reff{HIjp}-\reff{HIp} in  $[0,T)\times\Sc$. Moreover, by Proposition \ref{continuituboundary} and Proposition \ref{lemterm}, we have
$\overline v_\kappa^i(t,z)$$\leq$
  $\Liminf_{(t^{'},z^{'})\longrightarrow(t,z)}
  \underline v_\kappa^i(t^{'},z^{'})$,
for all $(t,z)$ $\in$ $[0,T)\times D_0$,
    and $\overline v^i_\kappa(T,z)$ $=$ $\underline v_\kappa^i(T,z)$ $=$ $g^i(z)$ for all $z$  $\in$ $\bar\Sc$. Then by Theorem \ref{theocompa}, we deduce
$\bar v_\kappa^i$ $\leq$ $\underline v_\kappa^i$ on $[0,T]\times\Sc$, which proves the continuity of $v_\kappa^i$ on $[0,T)\times\Sc$. On the other hand, suppose that
$w_\kappa^i$ is another constrained viscosity solution to \reff{HIjp}-\reff{HIp} with
\beqs
\lim_{(t',z')\rightarrow (t,z)} w_\kappa^i(t',z')=w_\kappa^i(t,z)=v_\kappa^i(t,z),\;\;\;\mbox{ for all }(t,z)\in [0,T)\times D_0,
\enqs
and
$w_\kappa^i(T,z)$ $=$ $g^i(z)$ for $z$ $\in$ $\bar\Sc$.
Then, $\bar w_\kappa^i(t,z)$ $=$ $\underline v_\kappa^i(t,z)$ $=$ $\bar v_\kappa^i(t,z)$ $=$ $\underline w_\kappa^i(t,z)$ for $(t,z)$ $\in$ $[0,T)\times D_0$ and
  $\bar w_\kappa^i(T,z)$ $=$ $\underline v_\kappa^i(T,z)$ $=$ $\bar v_\kappa^i(T,z)$ $=$ $\underline w_\kappa^i(T,z)$  for $z$ in $\bar\Sc$.
  We  then deduce by Theorem \ref{theocompa} that
  $\bar v_\kappa^i$ $\leq$ $\underline w_\kappa^i$ $\leq$ $\bar w_\kappa^i$ $\leq$ $\underline v_\kappa^i$ on $[0,T]\times\Sc$.  This proves $v_\kappa^i$ $=$ $w_\kappa^i$ on $[0,T]\times\Sc$. \\
\ep

\setcounter{equation}{0} \setcounter{Assumption}{0}
\setcounter{Theorem}{0} \setcounter{Proposition}{0}
\setcounter{Corollary}{0} \setcounter{Lemma}{0}
\setcounter{Definition}{0} \setcounter{Remark}{0}

\section{Numerical illustrations}
\label{sec:numeric}

In this section we provide some numerical results describing the value functions of the players and their
optimal policy.  A forward computation of the value function and the optimal strategy is in our knowledge
impossible due to the high dimension of the state process and the complexity of our model, therefore we
used a numerical scheme based on a quantization technique (see \cite{PPP04}). The convergence of the numerical
solution towards the real solution can be shown using consistency, monotonicity and stability arguments
and will be further investigated in a future work. A detailed description of the numerical algorithm
can be found in the Appendix.\\

Numerical tests are performed on the localized and discretized grid $[0,T]\times[x_{min},..,x_{max}]\times[y_{min}^1,..,y_{max}^1]\times[y_{min}^2,..,y_{max}^2]$. We used the following values for the parameters of the model: $T=1$, $\mu=0$, $\sigma=0.5$, $\zeta_{min}=-2.2$, $\zeta_{max}=1.8$,  $x_{min}$ $=$ $y_{min}^1$ $=$ $y_{min}^2=10$, $x_{max}$ $=$ $y_{max}^1$ $=$ $y_{max}^2$ $=90$, $\lambda=0.1$ and $g^1 =f^1$ and $g^2=f^2$, $\rho_1=\rho_2=0$,  $\phi_1=5$, $\phi_2=2.5$. Besides, the running costs are  $f^1(x,y^1,y_2)=(y^1-x)Q(y^1-y^2)$ and $f^2(x,y^1,y^2)=(y^2-x)Q(y^2-y^1)$ where
\begin{align*}
Q(x)=\mathds{1}_{]-\infty,-\Delta]}-\frac{x-\Delta}{2\Delta}\mathds{1}_{[-\Delta,\Delta]},
\end{align*}
with $\Delta=40$. Further, the terminal payoffs are chosen such that $g^1(x,y^1,y^2)=f^1(x,y^1,y^2)$ and $g^2(x,y^1,y^2)=f^2(x,y^1,y^2)$.

\begin{center}
\begin{figure}[htbp]
    \centering
        \includegraphics[width=15cm]{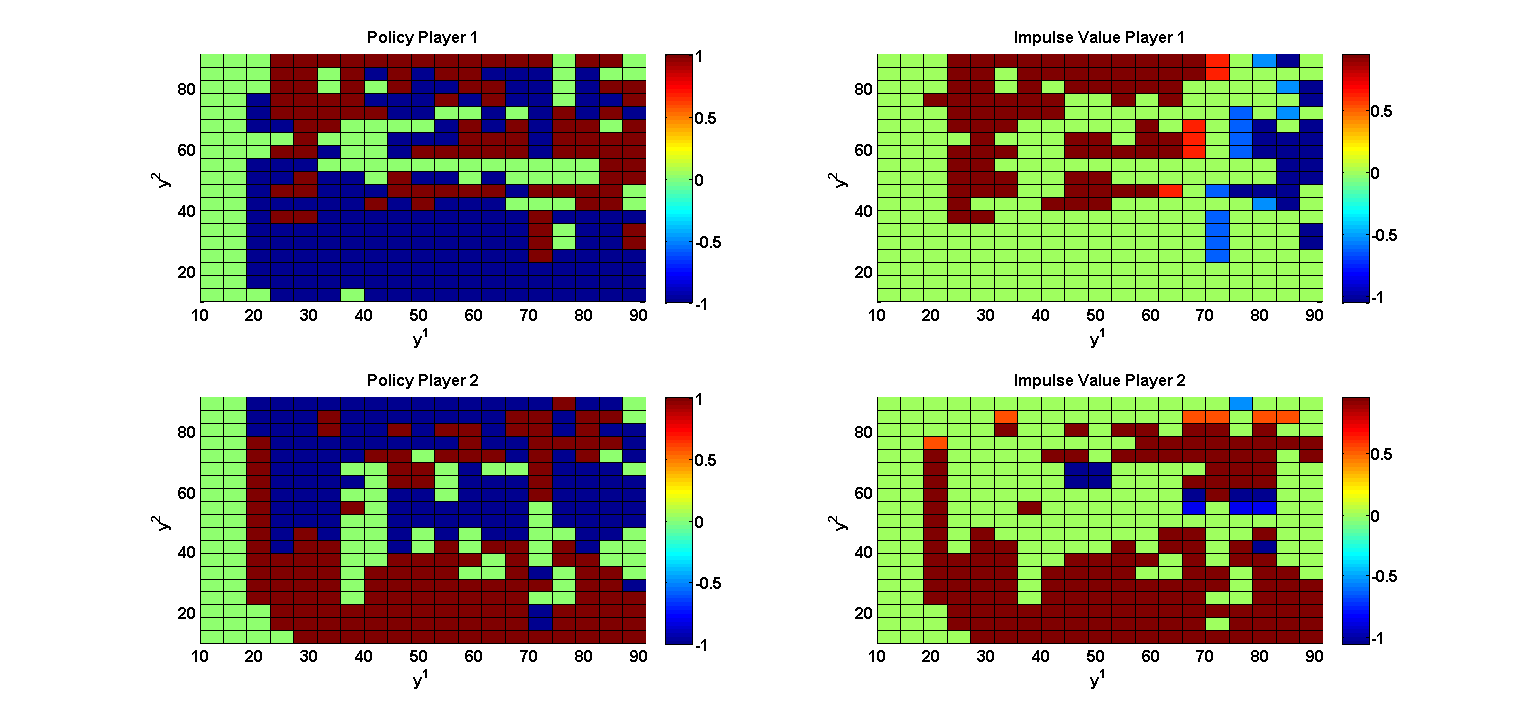}
    \caption{\emph{The optimal policies for a fixed $(t,x) = (,)$ for the first player (First Line) and the second player (Second Line). Color code: red: concerned player intervenes, green: concerned player waits, blue: concerned player endures the intervention of the other player. }}
    \label{PolicyFig}
\end{figure}
\end{center}

First,  the Figure \ref{PolicyFig} presents the optimal transaction policy for the two players, i.e. the different regions of interventions and continuations in the plane ($y^1,y^2)$ for $t=0.5$ and $x=50$~\euro/MWh. The first line (resp. second line) of Figure~\ref{PolicyFig} corresponds to the optimal policy regions and the corresponding interventions of the player 1 (resp. player 2). In the first column we can distinguish, for both of the players, three different regions, represented by three different colors, corresponding to the optimal action given a state $(y^1,y^2)$. Indeed, the blue region represents the states $(y^1,y^2)$ where a player is subject to the intervention of the other player, the green regions represents the states where a player chooses to not intervene and the red region represents the states where the player makes an intervention. The second column represents, whenever a player decides to intervene, the size of the intervention. If the quantity is positive it means that the price is increased and if it is negative it means that the price is lowered.

We can see that, as expected, both the players tend to keep the price spread $|y^1-y^2|$ as low as possible in order to avoid market share losses. In fact, for instance, at the state $(y^1=85,y^2=60)$, player 1 chooses to push down her price to keep an acceptable market share position. On the other hand, at the state $(y^1=30,y^2=70)$, player 1 chooses to push up her price which allows her to make benefits whilst keeping a reasonable market share position.

\begin{center}
\begin{figure}[htbp]
    \centering
        \includegraphics[width=0.5\textwidth]{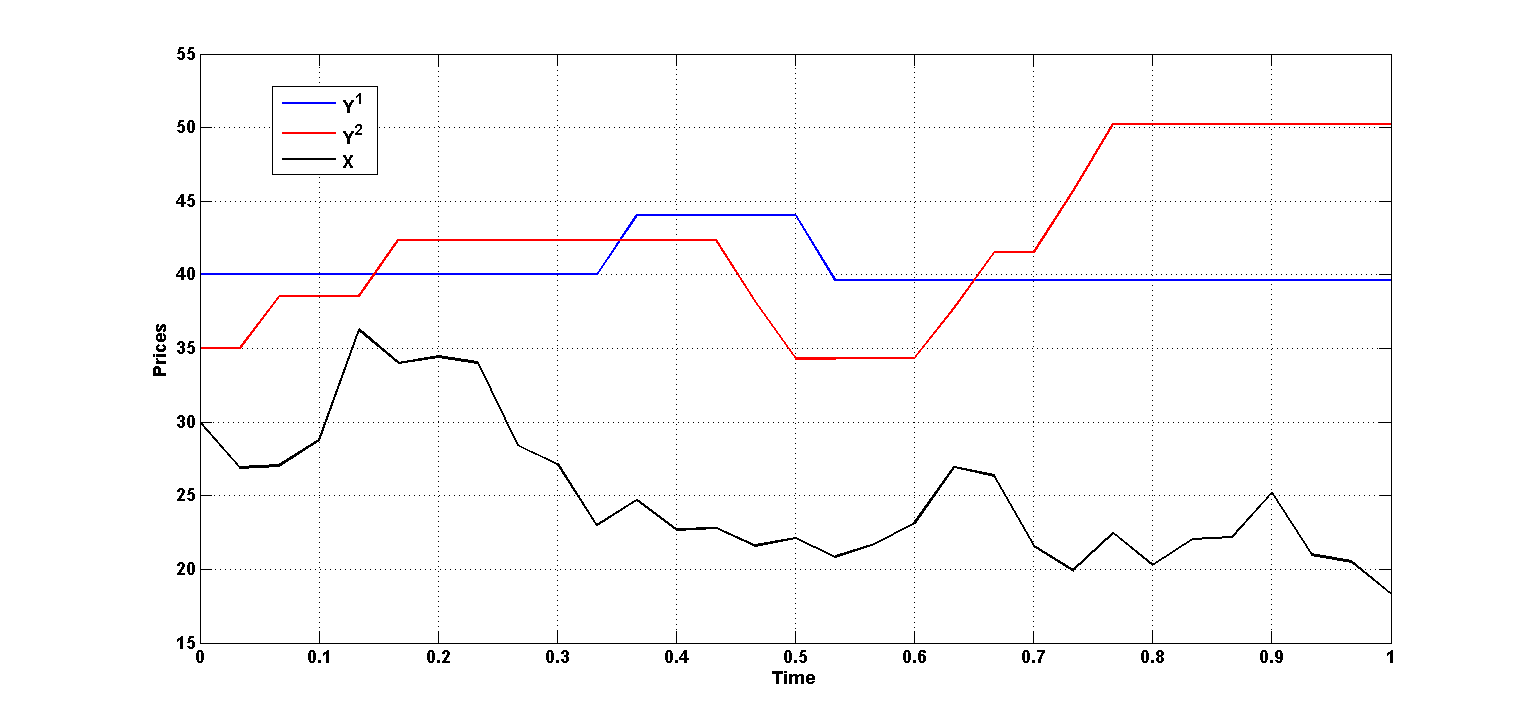}
	\includegraphics[width=0.45\textwidth]{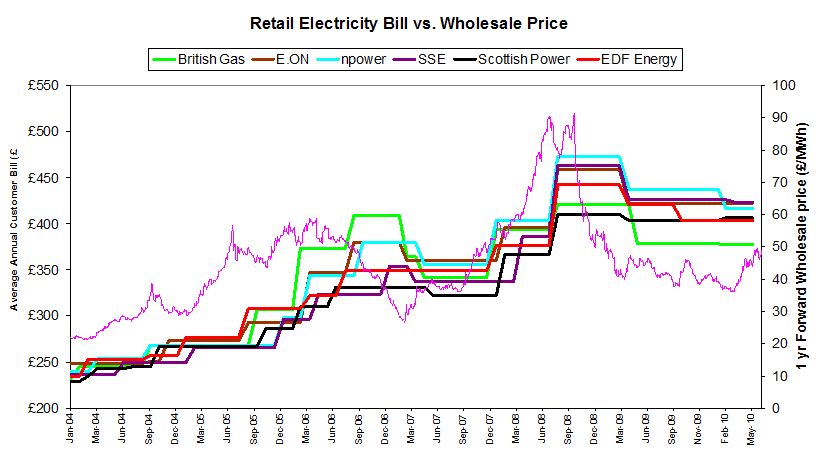}
    \caption{ (Left) One path-scenario of the wholesale market price and the players' retail prices. (Right) Retail electricity bill compared to wholesale price in the UK (source Ofgem). }
    \label{OnePath}
\end{figure}
\end{center}

Second, the Figure~\ref{OnePath} (Left) gives an example of a trajectory of the wholesale electricity price $X$ together with the corresponding retail prices trajectories  $Y^1$ and $Y^2$ of the two players, where the initial state is $(X_0=30,Y^1_0=40,Y^2_0=35)$. As a matter of comparison, Figure~\ref{OnePath} (Right) shows the trajectories of the wholesale price of electricity and retail prices of the six largest energy providers in the UK from January 2004 to March, 2010. We observe several comparable features of the optimal retailers price resulting from our impulse game and the real-life experience. Increases in the wholesale price is not immediately followed by an increase in retail prices. There is a delay given by the optimal time to reach the boundary of the action region. Further, even if our model only involves two players, we observe that they do not intervene at the same time, as it is the case in the UK market example. However, they appear to follow an almost synchronised behaviour: an increase by a first player is mostly to be followed by an increase of the second player and not by a decreases. Further, the optimal trajectories of the retail prices can be increasing while the wholesale price is decreasing (from $0.2$ to $0.3$ for instance), a phenomenon which is also observed in the UK case (from April, 2006 to March 2007, for instance). The optimal trajectories can also decrease, even if these decrease are limited compared to the same reference case of the UK market. Thus, contrary to the belief of the UK energy regulator, the Ofgem\footnote{The British energy regulator launched an inquiry on energy retailers in 2014. The headline findings of the assessment were: {\it (...) Possible tacit co-ordination: The assessment has not found evidence of explicit collusion between suppliers. However, there is evidence of possible tacit coordination reflected in the timing and size of price announcements and new evidence that prices rise faster when costs rise than they reduce when costs fall. Although tacit coordination is not a breach of competition law, it reduces competition and worsens outcomes for consumers}. Published on Ofgem website on June 26th, 2014, at the address: www.ofgem.gov.uk/press-releases/ofgem-refers-energy-market-full-competition-investigation.}, the observed behaviour of almost synchronised increase and decrease of retailers prices might not be the result of a tacit collusion mechanism, but is simply the result of optimal decision in a Nash equilibrium.

\begin{center}
\begin{figure}[htbp]
    \centering
        \includegraphics[width=13cm]{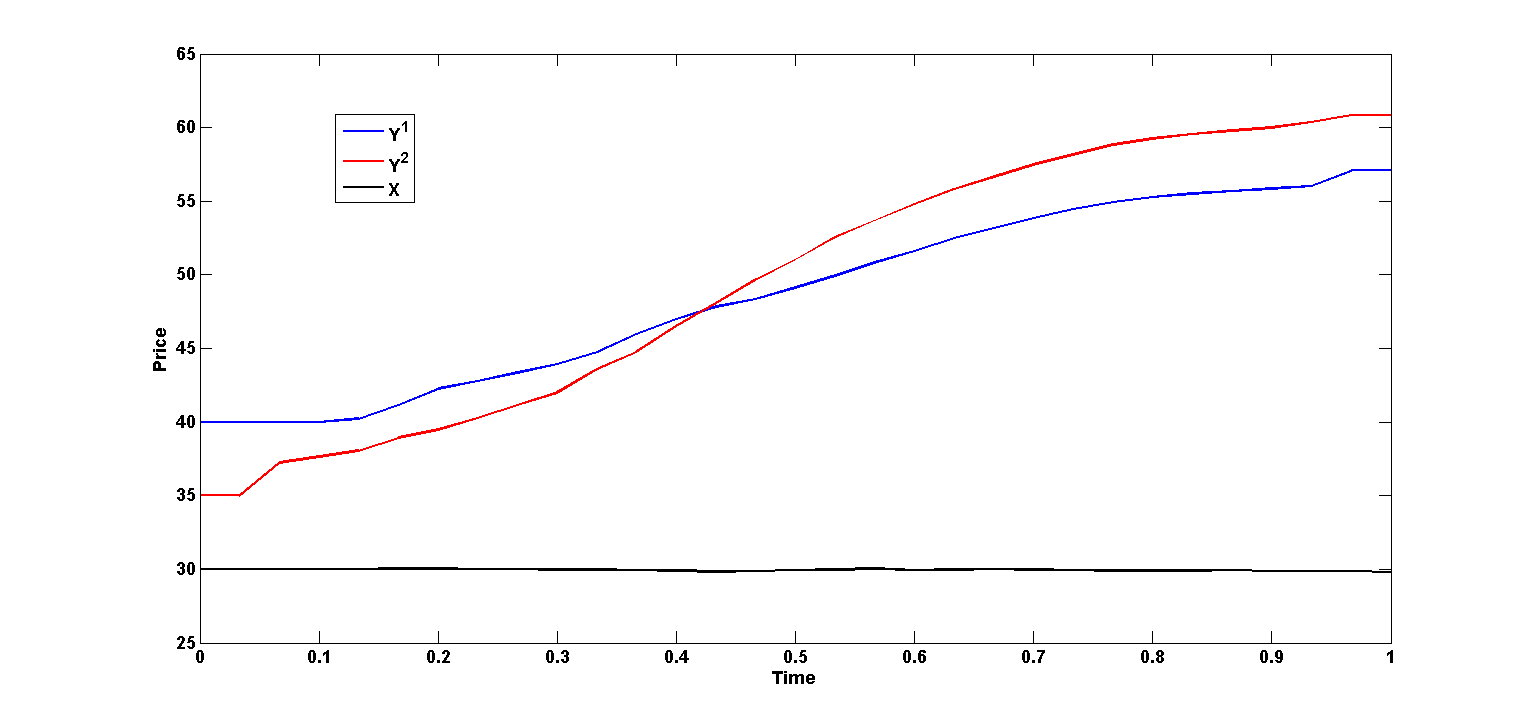}
    \caption{\emph{The average trajectory of the market price and the players' prices.}}
    \label{MeanPath}
\end{figure}
\end{center}

Finally, the Figure \ref{MeanPath} shows the average trajectories of the market price and the players' optimal retail price processes over ten thousand simulated trajectories of $X$, $Y^1$ and $Y^2$ on the horizon $[0,T]$. The initial state is the same as in the Figure~\ref{OnePath}. We notice that, although the wholesale price $X$ is a martingale, the  retail prices offered by the two players are increasing. In addition to this observation, we note that the players have almost the same tendency as they try to keep a balanced market share configuration until the maturity. With our choice of parameters, we observe that player 2 starts with a price lower than the player 1's price and attains the maturity with a higher price. This is because the interventions for the player 2 are less expensive making her more dynamic. We can also observe that, throughout the time period, the price spread between the two players is quite small preventing the market share to be imbalanced.

Our model suggests that the players would rather propose increasing prices to maximize their profit. This result might be surprising as one would expect that the players would stick to the wholesale price tendency and would propose a mean constant prices. But, in our model the market shares are split between the two players only according to the difference in the price they offer: consumers do not have an outside option to switch to another energy and no market entry of a competitor may threaten the two players for practicing increasing prices. The thing we find remarkable in this result is that without setting any potential communication device between the two players, we observe on average a behaviour that looks like tacit collusion.

\section*{Appendix}

In the following, we give a detailed description of the numerical procedure used to compute the value function and the optimal policies associated to the optimal control problem. We recall that we used a numerical scheme based on a quantization technique (see \cite{PPP04}) mixed with an iterative procedure. The convergence of the numerical solution towards the real solution can be shown using consistency, monotonicity and stability arguments and will be further investigated in a future work.\\

For a time step $h > 0$ on the interval $[0, T ],$
 we introduce  a numerical  backward scheme that approximates the solution of the HJB-QVI system via the couple of functions $v^i_h, i=1,2$ through:

 \begin{eqnarray}
\left\{
\begin{array}{rlll}
\mathcal{M}^i v^i_h(t,z)&-  v^i_h(t,z)\leq 0\\
v^i_{h}(t,z)&=max(\mathcal{M}^i\mathcal{H}^i v^i_{h}(t,z),\mathcal{H}^i v^i_{h}(t,z))  \quad \textrm{in} \quad  \overline{\mathcal{I}^i}\\
v^i_h(t,z)&=\max \left[ \mathbb{E}[v^i_h(t+h,Z^{t,z}_{t+h})]+\Sigma_i(t,z), \mathcal{M}^i v^i_h(t,z)\right] \textrm{in} ~~ {\mathcal{I}^i}\\
v^i_h(T,z)&= g^i(z), \quad \textrm{in} ~~\mathcal{S}.
\end{array}
\right.
\end{eqnarray}
 Where
 $$\Sigma_i(t,z)=\int^{t+h}_t f^i(Z_s^{t,z})ds.$$

This approximation scheme seems a
priori implicit due to the nonlocal obstacle terms $\mathcal{M}^i$ and $\mathcal{H}^i$. This is typically the case in impulse
control problems, and the usual way to circumvent this problem is to
iterate the scheme by considering a sequence of optimal stopping problems:
\begin{eqnarray}
\left\{
\begin{array}{rlll}
\mathcal{M}^i v^i_{h,n}(t,z)&-  v^i_{h,n+1}(t,z) \leq 0\\
v^i_{h,n+1}(t,z)&=max(\mathcal{M}^i\mathcal{H}^i v^i_{h,n}(t,z),\mathcal{H}^i v^i_{h,n}(t,z))  \quad \textrm{in} \quad  \overline{\mathcal{I}_i}\\
v^i_{h,n+1}(t,z)&=max \left[ \mathbb{E}[v^i_{h,n+1}(t+h,Z^{t,z}_{t+h})]+\Sigma_i(t,z), \mathcal{M}^i v^i_{h,n}(t,z)\right] ~~ \textrm{in} ~~ {\mathcal{I}^i}\\
v^i_{h,n+1}(T,z)&= g^i(z) \quad \textrm{in} ~~\mathcal{S}.
\end{array}
\right.
\end{eqnarray}

\subsubsection*{Time and Space  discretization}
$\bullet$ Now let us consider the time grid $\mathbb{T}:=\{ t_k = kh,~~ k = 0, .. ,M,~~ h = \frac{T}{M}\}$ and $M\in \mathbb{N} \setminus\lbrace 0\rbrace,$
$ z \in \mathcal{S}$ and starting from  a pair $(v^1_0,v^2_0)$  two fixed vectors. \\
\begin{eqnarray}
\left\{
\begin{array}{rlll}
\mathcal{M}^i v^i_{h,n}(t_k,z)&-  v^i_{h,n+1}(t_k,z)\leq 0\\
v^i_{h,n+1}(t_k,z)&=max(\mathcal{M}^i\mathcal{H}^i v^i_{h,n}(t_k,z),\mathcal{H}^i v^i_{h,n}(t_k,z)) \quad \textrm{in} \quad  \overline{\mathcal{I}^i}\\
v^i_{h,n+1}(t_k,z)&=max \left[ \mathbb{E}[v^i_{h,n+1}(t_{k+1},Z^{t_k,z}_{t_{k+1}})]+\Sigma_i(t_k,z), \mathcal{M}^i v^i_{h,n}(t_k,z)\right] \quad \textrm{in} \quad  {\mathcal{I}^i} \\
v^i_{h,n+1}(T,z)&= g^i(z_j) \quad \textrm{in} ~~\mathcal{S}.\\
\end{array}
\right.
\end{eqnarray}
$\bullet$ Let $\mathbb{X}$  the uniform grid on $[x_{min}, x_{max}]$ of step $dx=\frac{x_{max}-x_{min}}{(N_x-1)},$ where $x_{min}< x_{max} \in (0,+\infty)$ and $N_x >0 $.  For  $j = 0, ... ,N_x,$ we denote $x_j := x_{min}+j dx$.\\
$\bullet$ For $i \in\{1,2\}$, let $\mathbb{Y}_i$  the uniform grid on $[y_{min}^i, y_{max}^i]$ of step $dy_i=\frac{y_{max}^i-y_{min}^i}{(N_y-1)},$
 where $y_{min}^i< y_{max}^i \in (0,+\infty).$  For  $j= 0, ... ,N_y,$ we denote $y_{j}^i := y_{min}^i+jdy_i$.\\
 Let $z_j=(x_j,y_{j}^1,y_{j}^2)\in\mathbb{G}:=\mathbb{X}\times\mathbb{Y}_1\times \mathbb{Y}_2$, we define the following  problem:
 \begin{eqnarray}
\left\{
\begin{array}{rlll}
\mathcal{M}^i v^i_{h,n}(t_k,z_j)&-  v^i_{h,n+1}(t_k,z_j)\leq0\\
 v^i_{h,n+1}(t_k,z_j)&=max(\mathcal{M}^i\mathcal{H}^i v^i_{h,n}(t_k,z_j),\mathcal{H}^i v^i_{h,n}(t_k,z_j))  \quad \textrm{ in } \overline{\mathcal{I}^i}\cap\mathbb{G}\\
v^i_{h,n+1}(t_k,z_j)&=max \left[ \mathbb{E}[v^i_{h,n+1}(t_{k+1},Z^{t_k,z_j}_{t_{k+1}})]+\Sigma_i(t_k,z_j), \mathcal{M}^i v^i_{h,n}(t_k,z_j)\right]\textrm{ in }\mathcal{I}^i\cap\mathbb{G}\\
v^i_{h,n+1}(T,z_j)&= g^i(z_j)  \textrm{ in } \mathcal{S}\cap\mathbb{G}.\\
\end{array}
\right.
\end{eqnarray}
\subsection*{Quantization of the Brownian Motion }
To compute  the conditional expectations arising in the numerical backward scheme, we use the optimal quantization method. The main idea is to use the quantization theory to construct a suitable approximation of the Brownian motion.\\
 It is known that there exists a unique
strong solution for the SDE, $\frac{dX_s^{t,x}}{X_{s}^{t,x}} = \mu ds+\sigma dW_s $. So it suffices to consider a quantization of the Brownian motion itself. \\
Recall that the optimal quantization technique consists in approximating the expectation $\mathbb{E}[f(Z)]$, where $Z$
is a normal distributed variable and $f$ is a given real function, by
\beqs
\mathcal{E}[f(\xi)] &= & \sum_{k\in \xi(\Omega)}f(k)\mathbb{P}(\xi=k)\;.
\enqs
The distribution of the discrete variable $\xi$ is known for a fixed $N:=card(\xi(\Omega))$ and the approximation is optimal as the $L^2$-error
between $\xi$ and $Z$ is of order $1/N$ (see \cite{PPP04}). The optimal grid $\xi(\Omega)$ and the associated weights $\mathbb{P}(\xi=k)$ can be downloaded from
the website: http://www.quantize.maths-fi.com/downloads.\\
Let $N$ denote the number of elementary quantizers used to quantize process $\hat{X}_s$.
We replace ${X}_s$ in  by its quantized random vector $\hat{X}_s,$  the optimal quantization of $X_s$ and we obtain the quantized
dynamic programming backward scheme:
\begin{eqnarray} \label{A}
\left\{
\begin{array}{rlll}
\mathcal{M}^i v^i_{h,n}(t_k,z_j)&-  v^i_{h,n+1}(t_k,z_j)\leq0\\
 v^i_{h,n+1}(t_k,z_j)&=max(\mathcal{M}^i\mathcal{H}^i v^i_{h,n}(t_k,z_j),\mathcal{H}^i v^i_{h,n}(t_k,z_j))  \quad \textrm{ in } \overline{\mathcal{I}^i}\cap\mathbb{G}\\
v^i_{h,n+1}(t_k,z_j)&=max \left[ \mathcal{E}[v^i_{h,n+1}(t_{k+1},Z^{t_k,z_j}_{t_{k+1}}))]+\Sigma_i(t_k,z_j), \mathcal{M}^iv^i_{h,n}(t_k,z_j) \right], \textrm{ in }\mathcal{I}^i\cap\mathbb{G} \\
v^i_{h,n+1}(T,z_j)&= g^i(z_j) \;\textrm{ in } \mathcal{S}\cap\mathbb{G}.\\
\end{array}
\right.
\end{eqnarray}
 Hence, the expectations arising in the  backward scheme are approximated by
 \begin{align*}
  \mathcal{E}[v^i_{h,n+1}(t+h,Z^{t,z}_{t+h})]=\sum_{l=1}^{N}  v^i_{h,n+1}(t+h, x e^{(\mu-\frac{\sigma^2}{2})h+\sigma \sqrt{h} \hat{u_l} },y^i,y^j ) {P}_l,
 \end{align*}
 where $\hat{u}_l$ is the $N$ quantizer of the standard normal distribution.\\
The weight associated to this quantizer is ${P}_l=\mathbb{P}(\hat{U}=\hat{u}_l)$.
 The optimal grid  $\hat{u}_l$ and the associated weights ${P}_l$ are downloaded from
the website: http://www.quantize.maths-fi.com/downloads.\\
Finally, to approximate the integral $\Sigma_i$, we use the rectangle rule and we obtain:
 $$\Sigma_i(t_k,z_j)=\displaystyle{ \int_{t_{k}}^{t_{k+1}}f^i(Z_s^{t_k,z_j} } ds\simeq\displaystyle{ h f^i(z_j)}.$$

\subsection*{Final Numerical Algorithm}
Thus, considering the iterative scheme defined in \eqref{A}, we obtain the following final backward scheme for $(t_k,z_j)\in\mathbb{T}\times\mathbb{G}$:
\begin{algorithm}[H]
\caption{Policy iteration for system of QVIs  (one-player)}
\begin{algorithmic}[1]
\item[$1:$] Set $ \varepsilon >0 $ {(numerical tolerance)}
and $n_{max} \in \mathbb{N}$ (maximum iterations).
\item[$2:$] Pick initial guess: $v^i_{h,0} \in \mathbb{R}.$
\item[$3:$] Let $n=0$ (iteration counter) and $R^0= +\infty$..
\item[$4:$] \textbf{while} ${R}^n >\varepsilon $ and ~$ n\leq n_{max}$ \textbf{do}
\item[$5:$]
$\begin{array}{lcl}
v^i_{h,n+1}(T,z_j)=g^i(z_j) ,\\
v^i_{h,n+1}(t_k,z_j)=max \left[ \mathcal{E}[v^i_{h,n+1}(t_{k+1},Z^{t_k,z_j}_{t_{k+1}})]+\Sigma_i, \mathcal{M}^i v^i_{h,n}(t_k,z_j)\right].\\

\end{array}$
\item[$6:$]Let $R^ {n+1}$ be the largest pointwise residual to the QVI, i.e.\\$R^{n+1}=|| v_{h,n+1}^i-v_{h,n}^i||.$
\item[$7:$]  Let $n=n+1.$
\item[$8:$] \textbf{end while}.
\end{algorithmic}
\end{algorithm}

The final Algorithm is as follows
\begin{algorithm}[H]
\caption{Policy iteration for system of QVIs  (two players)}
\begin{algorithmic}[1]
\item[$1:$] Set $ \varepsilon >0 $ {(numerical tolerance)} $,0 < \alpha < 1,$
$ ~ r^0>0$ (relaxation parameters)
and $n_{max} \in \mathbb{N}$ (maximum iterations).
\item[$2:$] Pick initial guess: $(v^1_{h,0},v^2_{h,0}) \in \mathbb{R} \times \mathbb{R}.$
\item[$3:$] Let $n=0$ (iteration counter) and $R^0= +\infty$.
\item[$4:$] \textbf{while} ${R}^n >\varepsilon $ ~and ~$ n\leq n_{max}$ \textbf{do}
\item[$5:$] \textbf{for} $i=1,2$ (player $i$) \textbf{do}
\item[$6:$]  $ l=3-i$  (player $l.$)
\item[$7:$]  $\mathcal{C}^n_l :=\lbrace\mathcal{M}^l v^l_n -v^l_n<-r^n\rbrace\cap
\mathbb{G}.$
\item[$8:$]  For $t_k\in \mathbb{T}$ and $z_j  \notin  \mathcal{C}^n_l,$ let $v^{n+1}_i (t_k,z_j) =max(\mathcal{M}^i\mathcal{H}^i v^i_{h,n}(t_k,z_j),\mathcal{H}^i v^i_{h,n}(t_k,z_j)).$
\item[$9:$]  For $(t_k,z_j)  \in \mathbb{T}\times \mathcal{C}^n_l,$ solve for $v_{n+1}^i(t_k,z_j) $ by applying Algorithm $1$ to
$\min\lbrace{-\frac{\partial v_{n+1}^i}{\partial t} -\mathcal{L}^i v_{n+1}^i - f^i, v_{n+1}^i-\mathcal{M}^i v_{n+1}^i}\rbrace=0$
\item[$10:$]\textbf{end for.}
\item[$11:$]Let $R^ {n+1}$ be the largest pointwise residual to the system of QVIs, i.e.\\$R^{n+1}=max(|| v_{h,n+1}^1-v_{h,n}^1||,|| v_{h,n+1}^2-v_{h,n}^2||).$
\item[$12:$] $r^{n+1}:=max\lbrace  \alpha R^{n+1},\varepsilon\rbrace$
\item[$13:$]  Let $n=n+1.$
\item[$14:$] \textbf{end while}.
\end{algorithmic}
\end{algorithm}

\begin{small}

\end{small}

\end{document}